\documentclass[10pt,reqno]{amsart}
\ifdefined\pdfsuppresswarningpagegroup
\fi
\ifdefined\pdfminorversion
\fi

\usepackage{graphicx}
\usepackage{array}
\usepackage[figuresright]{rotating}
\usepackage{needspace}
\usepackage{placeins}
\usepackage{amsmath,amssymb,amsthm,mathtools}
\usepackage{enumitem}
\usepackage{xcolor}
\usepackage{hyperref}
\usepackage{bookmark}

\theoremstyle{definition}
\newtheorem{theorem}{Theorem}[section]
\newtheorem{lemma}[theorem]{Lemma}
\newtheorem{proposition}[theorem]{Proposition}
\newtheorem{corollary}[theorem]{Corollary}
\newtheorem*{maintheorema}{Theorem A}
\newtheorem*{maintheoremb}{Theorem B}

\theoremstyle{definition}
\newtheorem{definition}[theorem]{Definition}
\newtheorem{remark}[theorem]{Remark}

\providecommand{\sslash}{\mathbin{/\mkern-6mu/}}
\newcommand{\orbpartial}{\partial_{\mathrm{orb}}}
\newcommand{\SO}{\operatorname{SO}}
\newcommand{\Orth}{\operatorname{O}}
\newcommand{\MCG}{\operatorname{MCG}}
\newcommand{\Id}{\mathrm{id}}
\newcommand{\Cyc}[1]{\mathrm C_{#1}}
\newcommand{\Dih}[1]{\mathrm D_{#1}}

\newcolumntype{L}[1]{>{\raggedright\arraybackslash}p{#1}}

\hypersetup{
colorlinks=true,
linkcolor=blue,
citecolor=blue,
urlcolor=blue,
pdftitle={Orbifold Morse Handles II. Boundary Regions, Gluing Maps, and Surgery Cobordisms},
pdfauthor={Suhyoung Choi and Juseop Lee},
pdfsubject={Boundary regions and gluing maps of four-dimensional orbifold Morse handles},
pdfkeywords={orbifold Morse theory, orbifold handles, spherical orbifolds, surgery cobordism}
}

\title[Orbifold Morse Handles II]{Orbifold Morse Handles II.\\Boundary Regions, Gluing Maps, and Surgery Cobordisms}
\author{Suhyoung Choi}
\address{Department of Mathematical Sciences, KAIST, 291 Daehak-ro Yuseong-gu Daejeon, South Korea}
\email{schoi@math.kaist.ac.kr}
\author{Juseop Lee}
\address{Department of Mathematical Sciences, KAIST, 291 Daehak-ro Yuseong-gu Daejeon, South Korea}
\email{dlwntjq77@kaist.ac.kr}
\date{}

\begin{document}

\begin{abstract}
We study four-dimensional quotient handles defined by faithful representations of finite groups in \(\SO(4)\) and invariant splittings of indices \(1\), \(2\), and \(3\).
We identify their attaching and belt regions as oriented \(3\)-orbifolds described by labeled singular graphs.
We determine the orbifold isotopy class of the boundary gluing map.
Fix an orientation-reversing orbifold embedding of the attaching region into an oriented \(3\)-orbifold.
Assume that the embedding induces isomorphisms on local groups.
The regions, gluing class, and embedding determine the surgery.
For a fixed handle and attaching image, an extension condition on the belt region characterizes diffeomorphisms of the cobordisms fixing a specified common exterior pointwise.
For index-\(2\) handles, this condition is the vanishing of an integer twist in the fixed boundary coordinates.

We use these surgeries to construct a compact connected oriented \(4\)-orbifold cobordism from every closed connected oriented \(3\)-orbifold \(\mathcal O\) to \(S^3\).
A finite sequence of index-\(2\) orbifold handle attachments from \(\mathcal O\) to \(S^3\) exists if and only if \(\mathcal O\) has no tetrahedral, octahedral, or icosahedral vertices.
Let \(v_{\mathbf T}(\mathcal O)\), \(v_{\mathbf O}(\mathcal O)\), and \(v_{\mathbf I}(\mathcal O)\) denote the numbers of these vertices.
Among finite sequences of handle attachments of indices \(1\), \(2\), and \(3\) from \(\mathcal O\) to \(S^3\), the minimum total number of index-\(1\) and index-\(3\) handles is
\[
\left\lceil \frac{v_{\mathbf T}(\mathcal O)}{2}\right\rceil
+
\left\lceil \frac{v_{\mathbf O}(\mathcal O)}{2}\right\rceil
+
\left\lceil \frac{v_{\mathbf I}(\mathcal O)}{2}\right\rceil.
\]
\end{abstract}

\maketitle
{\small\setlength{\baselineskip}{9.5pt}\tableofcontents}

\section{Introduction}\label{sec:introduction}

Part~I~\cite[Definition~1.1 and Theorems~3.15 and~B]{choilee2026part1} defines ordered handles for faithful representations of finite groups in \(\SO(4)\).
It classifies their possible Morse indices and isomorphism classes.
For every ordered handle class of index \(1\), \(2\), or \(3\), we identify the attaching and belt regions.
We determine the orbifold isotopy class of their boundary gluing map.
These calculations describe the surgery induced by an attaching embedding.

We use the definitions and conventions for orbifolds, smooth orbifold maps, and orbifold diffeomorphisms from Part~I~\cite[Section~2]{choilee2026part1}.
Orbifolds are smooth, effective, and orientable unless specified otherwise.
Orbifolds with corners occur in Section~\ref{subsec:corner-smoothing}.
For an orbifold \(\mathcal O\), we write \(\orbpartial\mathcal O\) for its orbifold boundary.
If a finite group \(G\) acts smoothly and effectively on a manifold \(X\), we write \(X\sslash G\) for the quotient orbifold.
The symbol \(\cong\) denotes orbifold diffeomorphism for orbifolds, diffeomorphism for manifolds, and isomorphism for groups.
Orientation is specified when relevant.
For a Euclidean space \(V\), we write \(D(V)\) and \(S(V)\) for its unit disk and unit sphere.
If \(X\) and \(Y\) are smooth \(G\)-manifolds and the diagonal action on \(X\times Y\) is effective, we write \(X\times_GY\) for the quotient by this action.

\begin{definition}\label{def:full-orbifold-embedding}
Let \(P\) and \(Q\) be smooth orbifolds, possibly with boundary.
An \emph{orbifold embedding} \(f\colon P\to Q\) is a smooth orbifold map satisfying the following conditions.
Its underlying map is a topological embedding.
Its compatible local lifts are smooth embeddings equivariant with respect to injective homomorphisms of the local groups.
It is \emph{full} if the induced homomorphism from the local group at \(x\) to the local group at \(f(x)\) is an isomorphism for every \(x\in P\), including \(x\in\partial P\).
At boundary points, the local lifts are smooth up to the boundary.

An \emph{isotopy of full orbifold embeddings} is a family \(f_t\colon P\to Q\), \(t\in[0,1]\), of full orbifold embeddings.
In compatible orbifold charts, the local lifts \(\widetilde f_t\colon U\to V\) define a smooth map \(U\times[0,1]\to V\).
\end{definition}

For general orbifolds, a nonfull embedding can map a point with trivial local group to a point with nontrivial local group~\cite[Theorem~1(2) and Example~11]{borzellino2015suborbifold}.
An orbifold isotopy of the target preserves the isomorphism class of each local group.
If the local group at the image of a fixed source point changes isomorphism class during a family of embeddings, that family cannot extend to an orbifold isotopy of the target.
The proof of Corollary~\ref{cor:attaching-isotopy-cobordisms} uses the fullness assumption.

Throughout this paper, orbifold embeddings and their isotopies are full.

\begin{definition}\label{def:quotient-handle-boundary-regions}
Let \(G\) be a finite group.
Let \(\rho\colon G\to\SO(4)\) be a faithful representation.
Let \(\mathbb R^4=V^-\oplus V^+\) be an orthogonal \(\rho(G)\)-invariant splitting with \(1\leq\dim V^-\leq3\).

The \emph{orbifold handle} associated with this splitting is
\(H_\rho(V^-,V^+)=D(V^-)\times_GD(V^+)\).
Its \emph{index} is \(\dim V^-\).
Its \emph{attaching region} is
\(A_\rho(V^-,V^+)=S(V^-)\times_GD(V^+)\).
Its \emph{belt region} is
\(B_\rho(V^-,V^+)=D(V^-)\times_GS(V^+)\).
Set \(C_\rho(V^-,V^+)=S(V^-)\times_GS(V^+)\).

The inclusions identify \(C_\rho\) with \(\orbpartial A_\rho\) and \(\orbpartial B_\rho\), respectively.
The identity map of \(C_\rho\) therefore induces the orbifold diffeomorphism
\[
\varphi_\rho\colon\orbpartial B_\rho\longrightarrow\orbpartial A_\rho.
\]
\end{definition}

In the notation of Part~I~\cite[Definition~1.1]{choilee2026part1},
\(H_\rho(V^-,V^+)=H_\rho(V^-)\).
When the splitting is understood, we omit \((V^-,V^+)\) from the notation.

After the smoothing of Section~\ref{subsec:corner-smoothing}, radial projection induces an orientation-preserving orbifold diffeomorphism
\[
\orbpartial H_\rho\longrightarrow S^3\sslash G
\]
that identifies \(A_\rho\) and \(B_\rho\) with complementary suborbifolds.

Let \(\mathcal O\) be an oriented \(3\)-orbifold.
Let \(\psi\colon A_\rho\to\operatorname{int}\mathcal O\) be an orientation-reversing orbifold embedding with codimension-zero image.
We attach \(H_\rho\) to \(\mathcal O\times[0,1]\) along \(\mathcal O\times\{1\}\) by \(\psi\).
This attachment induces a surgery on \(\mathcal O\).
Its outgoing boundary is
\[
\left(\mathcal O\setminus\operatorname{int}\psi(A_\rho)\right)
\cup_{\psi|_{\orbpartial A_\rho}\circ\varphi_\rho}B_\rho.
\]
The ordered handle determines \(A_\rho\), \(B_\rho\), and \(\varphi_\rho\).
Together with the attaching embedding \(\psi\), these determine the outgoing orbifold.
The attaching embedding is an additional choice.

\subsection{Main results}\label{subsec:main-results}

Part~I~\cite[Theorem~A]{choilee2026part1} identifies the ordered handle classes with the normalizer orbits of invariant subspaces.

\begin{definition}\label{def:boundary-gluing-models}
Let \(H_\rho\), \(A_\rho\), \(B_\rho\), and \(\varphi_\rho\) be as in Definition~\ref{def:quotient-handle-boundary-regions}.
Let \(A\) and \(B\) be oriented \(3\)-orbifolds.
Let \(f_A\colon A_\rho\to A\) and \(f_B\colon B_\rho\to B\) be orientation-preserving orbifold diffeomorphisms.
The \emph{boundary gluing map} associated with \(f_A\) and \(f_B\) is
\[
\Phi_\rho
=
f_A|_{\orbpartial A_\rho}
\circ\varphi_\rho
\circ
\left(f_B|_{\orbpartial B_\rho}\right)^{-1}
\colon\orbpartial B\longrightarrow\orbpartial A.
\]
\end{definition}

The map \(\Phi_\rho\) is orientation-reversing because \(A_\rho\) and \(B_\rho\) induce opposite orientations on \(C_\rho\).
Tables~\ref{tab:index-two-models}--\ref{tab:polyhedral-regions} use the model identifications from Section~\ref{subsec:gluing-map} and the boundary coordinates from Section~\ref{sec:local-conventions}.

\begin{maintheorema}[Classification of handle surgeries]
Let \(\rho\colon G\to\SO(4)\) be a faithful representation of a nontrivial finite group.
Fix a \(\rho(G)\)-invariant orthogonal splitting \(\mathbb R^4=V^-\oplus V^+\) with \(1\leq\dim V^-\leq3\).
Let \(H_\rho,A_\rho,B_\rho,\varphi_\rho\) be as in Definition~\ref{def:quotient-handle-boundary-regions}.
\begin{enumerate}[label=\textup{(\roman*)}]
\item Tables~\ref{tab:index-two-models}--\ref{tab:polyhedral-regions} identify oriented labeled models for \(A_\rho\) and \(B_\rho\).
With the model identifications and boundary coordinates fixed in Sections~\ref{subsec:rational-tangle-models} and~\ref{subsec:gluing-map}, the tables determine the orbifold isotopy class \([\Phi_\rho]\) of Definition~\ref{def:boundary-gluing-models}.
\item Let \(\mathcal O\) be a closed oriented \(3\)-orbifold.
Let \(\psi_0,\psi_1\colon A_\rho\to\mathcal O\) be orientation-reversing full orbifold embeddings with \(\psi_0(A_\rho)=\psi_1(A_\rho)\).
Form the surgery cobordisms
\[
W_i=(\mathcal O\times[0,1])\cup_{\psi_i}H_\rho,
\qquad i=0,1,
\]
and fix compatible smoothings and a common closed exterior \(E\) as in Definition~\ref{def:relative-handle-attachments}.
There is an orientation-preserving orbifold diffeomorphism \(F\colon W_0\to W_1\) with \(F|_E=\Id_E\) if and only if there is an orientation-preserving orbifold diffeomorphism \(b\colon B_\rho\to B_\rho\) satisfying
\[
\psi_1|_{\orbpartial A_\rho}\circ\varphi_\rho\circ b|_{\orbpartial B_\rho}
=\psi_0|_{\orbpartial A_\rho}\circ\varphi_\rho.
\]
\end{enumerate}
\end{maintheorema}

Theorem~A\textup{(i)} follows from Theorem~\ref{thm:computed-boundary-regions} and Part~I~\cite[Theorems~A and~B and Table~1]{choilee2026part1}.
Theorem~A\textup{(ii)} is proved in Theorem~\ref{thm:relative-handle-attachments}.
For index-\(2\) handles, Corollary~\ref{cor:index-two-framing-relative-exterior} expresses this criterion as the vanishing of an integer twist in the fixed boundary coordinates.
By Corollary~\ref{cor:attaching-isotopy-cobordisms}, the cobordisms obtained from isotopic attaching embeddings are diffeomorphic relative to the incoming boundary.
The attaching images may vary during the isotopy.

The map \(\varphi_\rho\) is induced by the identity on \(C_\rho\).
Its expression in the fixed boundary coordinates need not be the identity.
Thus the orbifold diffeomorphism types of the attaching and belt regions alone do not determine the surgery.
For example, both handlebodies in a genus-one Heegaard splitting of a lens space are solid tori.
The diffeomorphism type of the lens space depends on the meridian slope specified by the gluing map.

We use the description of orientable \(3\)-orbifolds by labeled singular graphs from Section~\ref{subsec:orbifolds}.
Let \(\mathbf T,\mathbf O,\mathbf I<\SO(3)\) denote the tetrahedral, octahedral, and icosahedral rotation groups, respectively.
A vertex with local group \(\mathbf T\), \(\mathbf O\), or \(\mathbf I\) is called tetrahedral, octahedral, or icosahedral, respectively.
The edge-label triples at tetrahedral, octahedral, and icosahedral vertices are \((2,3,3)\), \((2,3,4)\), and \((2,3,5)\), respectively.

Druschel proves that every closed effective oriented \(3\)-orbifold bounds a compact effective oriented \(4\)-orbifold~\cite[Theorem~3.2]{druschel2000cobordism}.
Let \(\mathcal O\) be a closed connected oriented \(3\)-orbifold.
Theorem~\ref{thm:cobordism-to-sphere} yields a compact connected oriented \(4\)-orbifold cobordism from \(\mathcal O\) to \(S^3\) obtained by a finite sequence of orbifold handle attachments.
Every handle in this sequence has index \(1\), \(2\), or \(3\).

\begin{definition}
Let \(\mathcal O\) be a compact \(3\)-orbifold, not necessarily connected.
For \(P\in\{\mathbf T,\mathbf O,\mathbf I\}\), let \(v_P(\mathcal O)\) denote the number of vertices of \(\mathcal O\) with local group \(P\).

Let \(\mathcal O\) be a closed connected oriented \(3\)-orbifold.
Consider all finite sequences of orbifold handle attachments of indices \(1\), \(2\), and \(3\) from \(\mathcal O\) to \(S^3\).
Define \(c_{\mathrm{poly}}(\mathcal O)\) to be the minimum total number of index-\(1\) and index-\(3\) handles among these sequences.
\end{definition}

\begin{maintheoremb}
Let \(\mathcal O\) be a closed connected oriented \(3\)-orbifold.
Then
\[
c_{\mathrm{poly}}(\mathcal O)
=
\left\lceil\frac{v_{\mathbf T}(\mathcal O)}2\right\rceil
+
\left\lceil\frac{v_{\mathbf O}(\mathcal O)}2\right\rceil
+
\left\lceil\frac{v_{\mathbf I}(\mathcal O)}2\right\rceil.
\]

In particular, there exists a compact connected oriented \(4\)-orbifold \(W\) with
\[
\orbpartial W=(-\mathcal O)\sqcup S^3
\]
and a handle decomposition relative to \(\mathcal O\) consisting only of index-\(2\) handles if and only if
\[
v_{\mathbf T}(\mathcal O)
=
v_{\mathbf O}(\mathcal O)
=
v_{\mathbf I}(\mathcal O)
=
0.
\]
\end{maintheoremb}

Theorem~B follows from Theorems~\ref{thm:index-two-characterization} and~\ref{thm:minimum-polyhedral-handles}.

\subsection{Previous work and organization}\label{subsec:previous-work}

In classical Morse theory, crossing a nondegenerate critical value changes a sublevel set by a handle attachment~\cite{milnor1963morse}.
Wasserman proved the corresponding equivariant theorem for compact group actions~\cite[Theorem~4.6]{wasserman1969equivariant}.
Part~I~\cite{choilee2026part1} classifies the ordered handle classes in dimension four.

Boileau, Maillot, and Porti defined Dehn filling for torus boundary components by specifying a curve that becomes a meridian~\cite[Section~2.5]{boileau-maillot-porti2003}.
They also defined Dehn filling in the same way for boundary components whose underlying space is \(S^2\) with four cone points of order \(2\).
The orbifold handles considered here induce surgeries involving dihedral, tetrahedral, octahedral, and icosahedral singularities.
Together with smooth Dehn surgery, these surgeries produce every closed connected oriented \(3\)-orbifold from \(S^3\).
Unlike the manifold case, index-\(2\) handles do not always suffice.
Theorem~B characterizes when index-\(2\) handles suffice and determines the minimum total number of index-\(1\) and index-\(3\) handles otherwise.

Druschel's cobordism construction first removed vertices with icosahedral, octahedral, and tetrahedral local groups.
It then removed vertices with dihedral local groups in pairs and finally removed the cyclic local groups~\cite[proof of Theorem~3.2, pp.~51--53]{druschel2000cobordism}.
Druschel later defined surgery and handle attachment for direct-sum representations~\cite[Definitions~3.3 and~3.4]{druschel2015vertex}.
For \(q\geq2\), Druschel's pairwise removal of vertices of type \((2,2,q)\) corresponds to the index-\(1\) splitting in Proposition~\ref{prop:two-G11p-handles}.
The same representation has an index-\(2\) splitting whose attachment removes these two vertices without changing the underlying \(3\)-manifold.

Mecchia and Seppi determined the Seifert invariants and oriented diffeomorphism types of the spherical quotients that are fibered~\cite{mecchia2015fibered,mecchiaseppi2020diffeomorphism}.
Dunbar described the fibered quotients by rational tangle diagrams and the nonfibered quotients by labeled singular graphs~\cite{dunbar1981fibered,dunbar1988geometric,dunbar1994nonfibering}.
For each splitting chosen in Part~I to represent an orbit~\cite[Theorem~B and Table~1]{choilee2026part1}, we identify the attaching region \(A_\rho\) in the corresponding spherical quotient.
Under the identification in Section~\ref{subsec:corner-smoothing}, \(B_\rho\) is the closure of the complement of \(A_\rho\) in the same quotient.
When the boundary is a torus, we compute the boundary gluing map in meridian--longitude coordinates.
When a boundary component has underlying space \(S^2\) with four cone points of order \(2\), we use coordinates induced by the torus double cover.
When a boundary component has underlying space \(S^2\) with at most three cone points, we use the boundary map induced by the common quotient.

Section~\ref{sec:local-conventions} fixes the orbifold and smoothing conventions.
It also specifies the labeled models and boundary coordinates.
Section~\ref{sec:boundary-regions} identifies the attaching and belt regions and computes their boundary gluing maps.
It also classifies surgeries for a fixed handle and attaching image relative to the common exterior of Definition~\ref{def:relative-handle-attachments}.
Section~\ref{sec:surgery-cobordisms} constructs surgery cobordisms from closed connected oriented \(3\)-orbifolds to \(S^3\).
It characterizes when index-\(2\) handles suffice and determines the minimum total number of index-\(1\) and index-\(3\) handles.
Appendix~\ref{app:spherical-boundary-extension} determines which orientation-preserving orbifold diffeomorphisms of the boundaries of the index-\(1\) and index-\(3\) models extend over these models.
Appendix~\ref{subsec:illustrated-boundary-models} collects the labeled diagrams for the classification tables.

\section{Orbifold preliminaries and boundary coordinates}\label{sec:local-conventions}

In this section, we fix the orbifold conventions and the smoothings used below.
We specify the labeled models and boundary coordinates used in the classification.

Section~\ref{subsec:orbifolds} describes labeled \(3\)-orbifolds.
Section~\ref{subsec:corner-smoothing} applies corner smoothing and identifies the attaching and belt regions in \(S^3\sslash G\).
Section~\ref{subsec:fixed-region-models} specifies the oriented labeled models for these regions.
Section~\ref{subsec:rational-tangle-models} fixes the boundary coordinates and mapping class conventions used to express the boundary gluing maps.

\subsection{Labeled orbifolds}\label{subsec:orbifolds}

We describe orientable \(3\)-orbifolds by labeled singular graphs.

For an orbifold \(\mathcal O\), let \(\Sigma(\mathcal O)\) denote the set of points with nontrivial local group.
We write \(|\mathcal O|\) for its underlying space.

For an integer \(q\geq1\), let \(\Cyc{q}\) denote the cyclic group of order \(q\).
For an integer \(q\geq1\), let \(\Dih{2q}\) denote the rotational dihedral group of order \(2q\).
The group \(\Dih{2}\) is cyclic of order \(2\).
For integers \(q_1,\ldots,q_s\geq2\), let \(S^2(q_1,\ldots,q_s)\) denote the orientable \(2\)-orbifold whose underlying space is \(S^2\) and whose cone points have orders \(q_1,\ldots,q_s\).

\begin{definition}\label{def:labeled-pair}
Let \(M\) be an orientable smooth \(3\)-manifold, possibly with boundary.
A \emph{labeled singular graph} in \(M\) is a properly embedded locally finite graph \(\Gamma\subset M\) whose edges are labeled by integers at least \(2\).
We allow circle components.
Every interior vertex is trivalent.
The univalent vertices are precisely the points of \(\Gamma\cap\partial M\).
Each open edge is smoothly embedded.
A neighborhood of each interior vertex in \((M,\Gamma)\) is diffeomorphic to a \(3\)-ball with three radial segments meeting at its center.
At each univalent vertex, the incident edge is smoothly embedded and meets \(\partial M\) transversely.

The labels of the three edges incident to an interior vertex have one of the following forms.
\[
(2,2,q)\quad(q\geq2),
\qquad
(2,3,3),
\qquad
(2,3,4),
\qquad
(2,3,5).
\]
A vertex of type \((2,2,q)\) is called \emph{dihedral}.
A vertex of type \((2,3,3)\) is called \emph{tetrahedral}.
A vertex of type \((2,3,4)\) is called \emph{octahedral}.
A vertex of type \((2,3,5)\) is called \emph{icosahedral}.
The tetrahedral, octahedral, and icosahedral vertices are called \emph{polyhedral}.
The pair \((M,\Gamma)\) is called a \emph{labeled pair}.
\end{definition}

\begin{proposition}\label{prop:3orb-diffeo-criterion}
Let \(\mathcal O\) be an effective orientable smooth \(3\)-orbifold, possibly with boundary.
The underlying space \(|\mathcal O|\) admits an orientable smooth \(3\)-manifold structure for which \(\Sigma(\mathcal O)\) is a labeled singular graph.
Every point of \(\partial|\mathcal O|\) has cyclic local group.
The singular graph meets \(\partial|\mathcal O|\) transversely.

Conversely, every labeled pair \((M,\Gamma)\) determines an effective orientable smooth \(3\)-orbifold with underlying space \(M\) and singular set \(\Gamma\), uniquely up to orbifold diffeomorphism.
Let \((M',\Gamma')\) be another labeled pair.
Fix orientations on \(M\) and \(M'\).
Suppose that an orientation-preserving diffeomorphism \(F\colon M\to M'\) carries \(\Gamma\) onto \(\Gamma'\) and preserves the edge labels.
Then there is an orientation-preserving orbifold diffeomorphism between the orbifolds determined by \((M,\Gamma)\) and \((M',\Gamma')\).
\end{proposition}

\begin{proof}
By \cite[Theorems~2.3 and~2.4]{cooper2000three}, the local group at each interior point of \(\mathcal O\) is conjugate to a finite subgroup \(K<\SO(3)\).
The quotient space \(S^2/K\) is homeomorphic to \(S^2\).
Hence \(D^3/K\) is homeomorphic to \(D^3\).
The singular set of \(D^3\sslash K\) is either empty, an arc, or three radial arcs meeting at a vertex.
The edge labels and vertex types are those in Definition~\ref{def:labeled-pair}.
These quotient charts give \(|\mathcal O|\) the structure of an orientable topological \(3\)-manifold with a locally tame singular graph.

Let \(x\in\partial|\mathcal O|\).
Choose an invariant metric and an equivariant collar.
Corollary~3.6 of \cite{kankaanrinta2007equivariant} then provides a linear half-space chart at \(x\).
The local group at \(x\) fixes the inward normal vector.
Its action on the tangent plane of the boundary is faithful and orientation-preserving.
Thus the local group is isomorphic to a finite subgroup of \(\SO(2)\) and is cyclic.
If the local group is nontrivial, the image of its fixed set in the quotient half-ball is the inward normal half-line.

Choose a locally finite regular neighborhood \(N\) of \(\Sigma(\mathcal O)\) in \(|\mathcal O|\).
Model \(N\) by \(3\)-balls around the interior vertices, half-balls around the boundary vertices, and tubes along the remaining open edges.
Choose the tube parametrizations so that their end disks are attached smoothly to the boundaries of the vertex neighborhoods.
Use radial segments in the balls and transverse segments in the half-balls.
Circle components have solid torus neighborhoods.
These models define a smooth structure on \(N\).

The singular graph is closed in \(|\mathcal O|\).
By the relative smoothing theorem~\cite[Corollary~4.7]{pardon2021smoothing}, the smooth structure on a neighborhood of \(\Sigma(\mathcal O)\) extends to \(|\mathcal O|\).
With this smooth structure, the singular graph meets \(\partial|\mathcal O|\) transversely.

Conversely, let \((M,\Gamma)\) be a labeled pair.
A neighborhood of an open edge labeled \(q\) is modeled on an interval times \(D^2\sslash\Cyc{q}\).
At each interior vertex, the incident edge labels determine a finite subgroup \(K<\SO(3)\) up to conjugacy.
The quotient \(D^3\sslash K\) has the prescribed singular graph.
At each point of \(\Gamma\cap\partial M\), the local model is a cyclic quotient of a half-ball.

Cooper, Hodgson, and Kerckhoff proved the existence of an orientable \(3\)-orbifold with these local models~\cite[Theorem~2.5]{cooper2000three}.
Kleiner and Lott identified orientable \(3\)-orbifolds by their underlying manifolds and labeled singular graphs~\cite[Section~2.3.4, p.~113]{kleinerlott2014orbifolds}.
Together with their description of smooth orbifold structures~\cite[Section~2.1]{kleinerlott2014orbifolds}, this classification determines the associated smooth orbifold up to orbifold diffeomorphism.

Let \(F\colon M\to M'\) be an orientation-preserving diffeomorphism carrying \(\Gamma\) onto \(\Gamma'\) and preserving the edge labels.
Transport the orbifold structure on \(M\) to \(M'\) by \(F\).
The transported orbifold and the orbifold determined by \((M',\Gamma')\) have the same oriented labeled pair.
The classification yields an orientation-preserving orbifold diffeomorphism between these orbifolds.
Composing it with the map induced by \(F\) gives an orientation-preserving orbifold diffeomorphism between the orbifolds determined by \((M,\Gamma)\) and \((M',\Gamma')\).
\end{proof}

\subsection{Smoothing corners}\label{subsec:corner-smoothing}

Throughout this paper, we smooth the corners of orbifold handles and their attachments as in~\cite[Chapter~18]{fukaya2020kuranishi}.
For any prescribed open set containing the corners, choose the smoothing to agree with the original orbifold outside that set.
For \(H_\rho\), choose a convex rounding of \(D(V^-)\times D(V^+)\) depending only on the norms in the two factors and symmetric under their interchange.
Use the same rounding for all handles.
The rounded product admits an orientation-preserving equivariant diffeomorphism to \(D^4\) whose boundary restriction is radial projection.
Taking quotients identifies \(H_\rho\) with \(D^4\sslash\rho(G)\).
Radial rays identify the original attaching and belt regions with regions of the smoothed boundary.
We retain the notation \(A_\rho,B_\rho,C_\rho\) under this identification.

\subsection{Boundary models}\label{subsec:fixed-region-models}

We fix the oriented labeled models for the attaching and belt regions used in Section~\ref{sec:boundary-regions}.

By Proposition~\ref{prop:3orb-diffeo-criterion}, the singular sets of \(A_\rho\) and \(B_\rho\) meet their common boundary \(C_\rho\) transversely.
Hence no trivalent vertex of either singular graph lies in \(C_\rho\).

\begin{definition}\label{def:solid-torus-boundary-models}
Equip \(S^1\) and \(D^2\) with their standard orientations.
For an integer \(q\geq1\), let \(\Cyc{q}\) act on \(D^2\) by rotations.
Set \(T_q=S^1\times(D^2\sslash\Cyc{q})\).
Equip \(T_q\) with the product orientation.
\end{definition}

\begin{definition}\label{def:pillow-boundary-models}
Regard \(S^1\) as the unit circle and \(D^2\) as the closed unit disk in \(\mathbb C\).
Equip \(D^2\times S^1\) with the product orientation.

Let \(\mathfrak r\colon D^2\times S^1\to D^2\times S^1\) be given by
\[
\mathfrak r(z,v)=(\overline z,\overline v).
\]
Let \(\mathfrak r_\partial\) be the restriction of \(\mathfrak r\) to \(S^1\times S^1\).
Set
\[
\Sigma=(S^1\times S^1)\sslash\langle\mathfrak r_\partial\rangle.
\]
Equip \(\Sigma\) with the quotient orientation.
The orbifold \(\Sigma\) is \(S^2(2,2,2,2)\).
We call \(S^2(2,2,2,2)\) the \emph{pillow}.
Its four cone points are represented by \((\pm1,\pm1)\).

For an integer \(q\geq1\), let
\[
\mathfrak s_q(z,v)=(e^{2\pi i/q}z,v).
\]
Set \(\Gamma_q=\langle\mathfrak s_q,\mathfrak r\rangle\).
The generators satisfy
\[
\mathfrak s_q^q=\mathfrak r^2=1,
\qquad
\mathfrak r\mathfrak s_q\mathfrak r=\mathfrak s_q^{-1}.
\]
Set
\[
\mathcal P_q=(D^2\times S^1)\sslash\Gamma_q.
\]
We call an orbifold diffeomorphic to \(\mathcal P_q\) for some \(q\geq1\) a \emph{solid pillow}.
The map
\[
(z,v)\longmapsto(z^q,v)
\]
induces an orbifold diffeomorphism
\[
\orbpartial\mathcal P_q\longrightarrow\Sigma.
\]
We use this diffeomorphism to identify \(\orbpartial\mathcal P_q\) with \(\Sigma\).
\end{definition}

We use the terminology for oriented Seifert fibered \(3\)-orbifolds from \cite[Section~3]{mecchia2015fibered}.
Mecchia and Seppi described the fiber over a corner reflector as an exceptional interval fiber with a neighborhood orbifold diffeomorphic to a solid pillow~\cite[pp.~818--819, Figure~4]{mecchia2015fibered}.
The product fibration of \(D^2\times S^1\) induces a Seifert fibration on \(\mathcal P_q\).

The underlying space \(|\mathcal P_q|\) is a \(3\)-ball.
For \(q=1\), the singular set consists of two properly embedded arcs of order \(2\).
For \(q\geq2\), the singular set consists of an edge of order \(q\) joining two vertices of type \((2,2,q)\) and four edges of order \(2\) joining these vertices to \(\orbpartial\mathcal P_q\).
First take the quotient by \(\langle\mathfrak s_q\rangle\) using \((z,v)\mapsto(z^q,v)\).
The quotient by the involution induced by \(\mathfrak r\) then gives this singular graph.

\begin{definition}
Orient \(\mathbb R^3\) by its standard ordered basis.
Let \(D^3\subset\mathbb R^3\) be the unit ball.
Equip \(S^2=\orbpartial D^3\) with its boundary orientation.
For a finite subgroup \(K<\SO(3)\), set
\[
N_K=D^3\sslash K,
\qquad
\Sigma_K=S^2\sslash K.
\]

Orient \(D^1=[-1,1]\) by its increasing coordinate.
Equip \(S^0=\orbpartial D^1\) with its boundary orientation.
Equip \(S^0\times N_K\) with the product orientation.
Thus \(S^0\times N_K=N_K\sqcup(-N_K)\), with components ordered by \(+1,-1\).

For an oriented closed \(2\)-orbifold \(\Sigma\), set
\[
I(\Sigma)=-(D^1\times\Sigma),
\]
where the minus sign reverses the product orientation.
The canonical map \(S^0\times\Sigma\to\orbpartial I(\Sigma)\) is orientation-reversing on each component.
\end{definition}

\begin{definition}\label{def:line-region-models}
Let \(K<\SO(3)\) be finite.
Let \(\nu\in\Orth(3)\setminus\SO(3)\) normalize \(K\).
Assume that \(\nu^2\in K\).
For \(k\in K\), define
\[
\widehat k(t,x)=(t,kx),
\qquad
\widehat\nu(t,x)=(-t,\nu x)
\]
on \(D^1\times S^2\).
Set
\[
\widehat K_\nu
=
\langle\widehat k,\widehat\nu\mid k\in K\rangle
\]
and
\[
\mathcal M(K,\nu)
=
(D^1\times S^2)\sslash\widehat K_\nu.
\]
Equip \(\mathcal M(K,\nu)\) with the orientation opposite to the quotient orientation induced by the product orientation on \(D^1\times S^2\).
\end{definition}

The group \(\widehat K_\nu\) is finite.
Its action on \(D^1\times S^2\) is effective and orientation-preserving.
The inclusion \(\{1\}\times S^2\hookrightarrow D^1\times S^2\) induces an orientation-reversing orbifold diffeomorphism
\[
\Sigma_K\longrightarrow\orbpartial\mathcal M(K,\nu).
\]
If \(\nu'=k\nu\) for some \(k\in K\), then \(\widehat K_{\nu'}=\widehat K_\nu\).
If \((K',\nu')=(hKh^{-1},h\nu h^{-1})\) for some \(h\in\SO(3)\), then \((t,x)\mapsto(t,hx)\) induces an orientation-preserving orbifold diffeomorphism
\[
\mathcal M(K,\nu)\longrightarrow\mathcal M(K',\nu').
\]

For the cyclic and dihedral models in Tables~\ref{tab:cyclic-line-regions} and~\ref{tab:dihedral-line-regions}, identify
\[
\mathbb R^3=\mathbb R\oplus\mathbb R^2,
\]
with the first summand equal to the first coordinate axis.
Set
\[
R_\theta
=
\begin{pmatrix}
1&0&0\\
0&\cos\theta&-\sin\theta\\
0&\sin\theta&\cos\theta
\end{pmatrix},
\qquad
J
=
\begin{pmatrix}
-1&0&0\\
0&1&0\\
0&0&-1
\end{pmatrix}.
\]
For an integer \(q\geq1\), set
\[
\Cyc{q}=\langle R_{2\pi/q}\rangle,
\qquad
\Dih{2q}=\langle R_{2\pi/q},J\rangle.
\]
The group \(\Dih{2}\) is cyclic of order \(2\).
The cyclic quotient orbifolds satisfy
\[
\Sigma_{\Cyc{1}}=S^2,
\qquad
\Sigma_{\Cyc{q}}=S^2(q,q)\quad(q\geq2).
\]
The dihedral quotient orbifolds satisfy
\[
\Sigma_{\Dih{2}}=S^2(2,2),
\qquad
\Sigma_{\Dih{2q}}=S^2(2,2,q)\quad(q\geq2).
\]
In particular, \(N_{\Cyc{1}}=D^3\).

Set \(\theta_q=(1-q)\pi/q\) and
\[
\alpha=-I_3,
\qquad
\gamma=
\begin{pmatrix}
1&0&0\\
0&1&0\\
0&0&-1
\end{pmatrix},
\qquad
\delta_q=
\begin{pmatrix}
-1&0&0\\
0&\cos\theta_q&-\sin\theta_q\\
0&\sin\theta_q&\cos\theta_q
\end{pmatrix}.
\]
The matrices \(\alpha\) and \(\delta_q\) normalize both \(\Cyc{q}\) and \(\Dih{2q}\).
The matrix \(\gamma\) normalizes \(\Cyc{q}\).
They satisfy
\[
\gamma R_{2\pi/q}\gamma^{-1}=R_{2\pi/q}^{-1},
\qquad
\delta_qR_{2\pi/q}\delta_q^{-1}=R_{2\pi/q},
\qquad
\delta_qJ\delta_q^{-1}=R_{2\pi/q}J,
\]
and
\[
\alpha^2=\gamma^2=I_3,
\qquad
\delta_q^2=R_{2\pi/q}.
\]

For the tetrahedral model, choose the standard inclusion \(\mathbf T<\mathbf O\) and coordinates in which
\[
\omega_{\mathbf T}
=
\begin{pmatrix}
0&-1&0\\
1&0&0\\
0&0&1
\end{pmatrix}
\in\mathbf O\setminus\mathbf T.
\]
Set \(\varepsilon_{\mathbf T}=-\omega_{\mathbf T}\).
Since \(\mathbf T\) has index \(2\) in \(\mathbf O\),
\[
\omega_{\mathbf T}\mathbf T\omega_{\mathbf T}^{-1}
=
\mathbf T,
\qquad
\omega_{\mathbf T}^2\in\mathbf T.
\]
Hence \(\varepsilon_{\mathbf T}\) normalizes \(\mathbf T\) and satisfies
\(\varepsilon_{\mathbf T}^2\in\mathbf T\).
Let \(\omega'_{\mathbf T}\in\mathbf O\setminus\mathbf T\) and set \(\varepsilon'_{\mathbf T}=-\omega'_{\mathbf T}\).
There is \(k\in\mathbf T\) such that
\[
\omega'_{\mathbf T}=k\omega_{\mathbf T}.
\]
Hence
\[
\varepsilon'_{\mathbf T}=k\varepsilon_{\mathbf T}.
\]
The equality \(\widehat K_{k\nu}=\widehat K_\nu\) therefore gives
\[
\widehat{\mathbf T}_{\varepsilon'_{\mathbf T}}
=
\widehat{\mathbf T}_{\varepsilon_{\mathbf T}}.
\]
Thus \(\mathcal M(\mathbf T,\varepsilon_{\mathbf T})\) is independent of \(\omega_{\mathbf T}\in\mathbf O\setminus\mathbf T\).

The polyhedral quotient orbifolds are
\[
\Sigma_{\mathbf T}=S^2(2,3,3),
\qquad
\Sigma_{\mathbf O}=S^2(2,3,4),
\qquad
\Sigma_{\mathbf I}=S^2(2,3,5).
\]

Let \(\rho\colon G\to\SO(4)\) be a faithful representation of a finite group.
Let \(\mathbb R^4=L\oplus W\) be a \(\rho(G)\)-invariant orthogonal splitting with \(\dim L=1\).
Let \(\chi_L\colon G\to\{\pm1\}\) be the character of the action on \(L\).

We identify the boundary regions for the index-\(1\) splitting \(L\oplus W\).
The index-\(3\) case follows by ordered duality.

Choose \(u\in S(L)\).
Choose an orientation-preserving linear isometry
\[
L\oplus W\longrightarrow\mathbb R\oplus\mathbb R^3
\]
that sends \(u\) to \((1,0,0,0)\), maps \(L\) onto the first summand, and maps \(W\) onto the second summand.
Under this identification, we regard \(\rho(g)|_W\) as an element of \(\Orth(3)\) for every \(g\in G\).

If \(\chi_L\) is trivial, let \(K<\SO(3)\) be the image of the action of \(G\) on \(W\).
If \(\chi_L\) is nontrivial, let \(K<\SO(3)\) be the image of the action of \(\ker\chi_L\) on \(W\).
If \(\chi_L\) is nontrivial, choose \(g_0\in G\) with \(\chi_L(g_0)=-1\) and set \(\nu=\rho(g_0)|_W\).

\begin{lemma}\label{lem:line-region-models}
Let \(t\in S^0\) and \(x\in S^2\).
Write \([x]\) for the class of \(x\) in \(\Sigma_K\).
\begin{enumerate}
\item[\textup{(i)}]
If \(\chi_L\) is trivial, there are orientation-preserving orbifold diffeomorphisms
\[
A_\rho\cong S^0\times N_K,
\qquad
B_\rho\cong I(\Sigma_K).
\]
Under these identifications, the boundary gluing map is
\begin{equation}\label{eq:line-gluing-disconnected}
\begin{aligned}
\Phi_\rho&\colon\orbpartial I(\Sigma_K)\longrightarrow S^0\times\Sigma_K,\\
(t,[x])&\longmapsto(t,[x]).
\end{aligned}
\end{equation}

\item[\textup{(ii)}]
If \(\chi_L\) is nontrivial, then \(\nu=\rho(g_0)|_W\) lies in
\(\Orth(3)\setminus\SO(3)\), normalizes \(K\), and satisfies \(\nu^2\in K\).
There are orientation-preserving orbifold diffeomorphisms
\[
A_\rho\cong N_K,
\qquad
B_\rho\cong\mathcal M(K,\nu).
\]
Under these identifications, the boundary gluing map is
\begin{equation}\label{eq:line-gluing-connected}
\begin{aligned}
\Phi_\rho&\colon\orbpartial\mathcal M(K,\nu)\longrightarrow\Sigma_K,\\
[(1,x)]&\longmapsto[x].
\end{aligned}
\end{equation}
\end{enumerate}
Both maps reverse the boundary orientations.
\end{lemma}

\begin{proof}
Suppose first that \(\chi_L\) is trivial.
The group fixes both points of \(S(L)\).
The chosen coordinates therefore induce an orbifold diffeomorphism \(A_\rho\to S^0\times N_K\).
It preserves orientation by the boundary orientation on \(S^0\).
The action on \(D(L)\) is trivial.
Hence \(B_\rho\cong D^1\times\Sigma_K\) as an orbifold.
The orientation on \(B_\rho\) induced by the handle is opposite to the product orientation.
Hence there is an orientation-preserving orbifold diffeomorphism \(B_\rho\to I(\Sigma_K)\).

Suppose that \(\chi_L\) is nontrivial.
Since \(\rho(g_0)\in\SO(4)\) and \(\chi_L(g_0)=-1\), we have \(\det(\nu)=-1\).
The normality of \(\ker\chi_L\) in \(G\) implies \(\nu K\nu^{-1}=K\).
Since \(g_0^2\in\ker\chi_L\), we have \(\nu^2\in K\).
The action on \(S(L)\) is transitive.
The stabilizer of \(u\) is \(\ker\chi_L\).
The inclusion \(\{u\}\times D(W)\hookrightarrow S(L)\times D(W)\) thus induces an orbifold diffeomorphism \(N_K\to A_\rho\).
It preserves orientation because the chosen coordinates send \(u\) to the positive endpoint of \(D^1\).
The action of \(G\) on \(D(L)\times S(W)\) has image \(\widehat K_\nu\).
By the orientation convention in Definition~\ref{def:line-region-models}, the induced orbifold diffeomorphism \(B_\rho\to\mathcal M(K,\nu)\) is orientation-preserving.

Both model identifications are induced by the quotient of \(S(L)\times S(W)\) on the common boundary.
If \(\chi_L\) is trivial, the point represented by \((t,x)\) has coordinates \((t,[x])\) in both boundary models.
If \(\chi_L\) is nontrivial, the representative \((1,x)\) has coordinates \([(1,x)]\) in \(\orbpartial\mathcal M(K,\nu)\) and \([x]\) in \(\orbpartial N_K\).
These coordinates yield \eqref{eq:line-gluing-disconnected} and \eqref{eq:line-gluing-connected}, respectively.
The orientation conventions for \(I(\Sigma_K)\) and \(\mathcal M(K,\nu)\) make both maps orientation-reversing.
\end{proof}

\subsection{Boundary coordinates and mapping class groups}\label{subsec:rational-tangle-models}

We fix the boundary coordinates and mapping class conventions used to express the boundary gluing maps of Definition~\ref{def:boundary-gluing-models}.
Proposition~\ref{lem:ordered-two-corner-gluing} computes the boundary gluing map for a class of index-\(2\) orbifold handles.

For an oriented \(2\)-orbifold \(S\), write \(\MCG(S)\) for its orbifold mapping class group and \(\MCG^+(S)\) for its orientation-preserving subgroup.

\begin{lemma}\label{lem:gluing-isotopy}
Let \(\mathcal O_1\) and \(\mathcal O_2\) be oriented smooth orbifolds with boundary.
For \(i=1,2\), let \(A_i\) be a union of boundary components of \(\mathcal O_i\).
Let \(f_0,f_1\colon A_1\to A_2\) be orientation-reversing orbifold diffeomorphisms that are isotopic through orbifold diffeomorphisms.
Then there is an orientation-preserving orbifold diffeomorphism
\[
\mathcal O_1\cup_{f_0}\mathcal O_2
\longrightarrow
\mathcal O_1\cup_{f_1}\mathcal O_2
\]
that restricts to the identity on \(\mathcal O_1\) and maps \(\mathcal O_2\) onto \(\mathcal O_2\).
\end{lemma}

\begin{proof}
Let \(f_t\colon A_1\to A_2\) be a smooth orbifold isotopy from \(f_0\) to \(f_1\).
Set \(g_t=f_t\circ f_0^{-1}\).
By the orbifold collar construction in Choi~\cite[Section~4.4.2]{choi2012geometric}, choose an orbifold collar of \(A_2\) in \(\mathcal O_2\) and identify it with \(A_2\times[0,1)\).
Choose a smooth function \(\rho\colon[0,1)\to[0,1]\) such that \(\rho(s)=1\) for \(0\leq s\leq1/4\) and \(\rho(s)=0\) for \(1/2\leq s<1\).
Define \(G(x,s)=(g_{\rho(s)}(x),s)\) on the collar.
Its inverse is \((x,s)\mapsto(g_{\rho(s)}^{-1}(x),s)\).
Since \(g_0=\Id_{A_2}\), the map \(G\) is the identity for \(s\geq1/2\).
Thus \(G\) extends by the identity to an orbifold diffeomorphism of \(\mathcal O_2\).
Each \(g_t\) preserves orientation, so \(G\) preserves orientation.

The restriction of \(G\) to \(A_2\) is \(g_1=f_1\circ f_0^{-1}\).
Hence \(G\circ f_0=f_1\).
Therefore \(\Id_{\mathcal O_1}\) and \(G\) induce an orientation-preserving orbifold diffeomorphism
\[
\mathcal O_1\cup_{f_0}\mathcal O_2
\longrightarrow
\mathcal O_1\cup_{f_1}\mathcal O_2.
\]
For \((x,s)\in A_2\times[0,1/4]\), we have \(G(x,s)=((f_1\circ f_0^{-1})(x),s)\).
Hence the induced map is smooth on the glued orbifold.
\end{proof}

For ordinary \(2\)-handles, the attaching embedding determines the framing of the attaching knot~\cite[Section~2.1]{nicolaescu2007morse}.
Changing the framing changes this embedding.
The boundary identification \(\varphi_\rho\) of the abstract handle is unchanged.

\begin{remark}\label{rem:gluing-not-unique}
The converse of Lemma~\ref{lem:gluing-isotopy} does not hold in general.
The orbifold diffeomorphism type of the glued orbifold need not determine the isotopy class of the gluing map.

Let \(\mathcal O_i\) and \(A_i\) be as in Lemma~\ref{lem:gluing-isotopy}.
Let \(f\colon A_1\to A_2\) be an orientation-reversing orbifold diffeomorphism.
Suppose that \(h\colon A_1\to A_1\) extends to an orientation-preserving orbifold diffeomorphism of \(\mathcal O_1\).
Then there is an orientation-preserving orbifold diffeomorphism
\[
\mathcal O_1\cup_f\mathcal O_2
\longrightarrow
\mathcal O_1\cup_{f\circ h^{-1}}\mathcal O_2.
\]
\end{remark}

For torus boundary components, we use meridian--longitude bases.
For boundary components whose underlying space is \(S^2\) with four cone points of order \(2\), we use the coordinates fixed in Definition~\ref{def:pillow-boundary-coordinates}.

\begin{definition}\label{def:torus-boundary-coordinates}
For an integer \(q\geq1\), let \(T_q\) be the model in Definition~\ref{def:solid-torus-boundary-models}.
Choose \(z_0\in S^1\) and \(w_0\in\orbpartial(D^2\sslash\Cyc{q})\).
Let \(\mu_q\) be the positively oriented boundary of the orbifold disk \(\{z_0\}\times(D^2\sslash\Cyc{q})\).
Let \(\lambda_q\) be the positively oriented curve \(S^1\times\{w_0\}\).
We use the same symbols for their homology classes.
Then \((\mu_q,\lambda_q)\) is an oriented basis of \(H_1(\orbpartial T_q;\mathbb Z)\).

When the attaching and belt models are \(T_{q_A}\) and \(T_{q_B}\), we write \((\mu_A,\lambda_A)\) and \((\mu_B,\lambda_B)\) for the corresponding oriented bases.
\end{definition}

\begin{definition}\label{def:pillow-boundary-coordinates}
Let \(\Sigma\) be the pillow in Definition~\ref{def:pillow-boundary-models}.
Using its quotient presentation, identify the torus double cover of \(\Sigma\) with \(\widetilde\Sigma=\mathbb R^2/\mathbb Z^2\) by \([s,t]\mapsto(e^{2\pi i s},e^{2\pi i t})\).
Let \(a,b\in H_1(\widetilde\Sigma;\mathbb Z)\) be the classes represented by \(t\mapsto[t,0]\) and \(t\mapsto[0,t]\), respectively.
The ordered pair \((a,b)\) is the standard oriented basis of \(H_1(\widetilde\Sigma;\mathbb Z)\).
Assign \emph{slope} \(y/x\) to a primitive class \(xa+yb\), with slope \(\infty\) when \(x=0\).
The image of \(D^2\times\{i\}\) in \(\mathcal P_q\) is an orbifold disk diffeomorphic to \(D^2\sslash\Cyc{q}\).
A lift of its boundary to \(\widetilde\Sigma\) represents \(a\).
Thus this disk has slope \(0\).

The deck involution of \(\widetilde\Sigma\to\Sigma\) is \(-\Id\).
Order its four fixed points by
\[
x_1=[0,1/2],\qquad
x_2=[1/2,1/2],\qquad
x_3=[0,0],\qquad
x_4=[1/2,0].
\]

Let \(\iota\colon\Sigma\to\Sigma\) be the orientation-reversing orbifold diffeomorphism induced by \((u,v)\mapsto(u,\overline v)\).
We use the same symbol for its class in \(\MCG(\Sigma)\).

Let \(\widetilde\tau_a\), \(\widetilde\tau_b\), and \(\widetilde\sigma_a\) be the torus automorphisms represented in the ordered basis \((a,b)\) by
\[
\widetilde\tau_a=
\begin{pmatrix}
1&1\\
0&1
\end{pmatrix},
\qquad
\widetilde\tau_b=
\begin{pmatrix}
1&0\\
-1&1
\end{pmatrix},
\qquad
\widetilde\sigma_a=
\begin{pmatrix}
0&-1\\
1&0
\end{pmatrix}.
\]
These automorphisms commute with \(-\Id\) and induce orientation-preserving orbifold diffeomorphisms \(\tau_a\), \(\tau_b\), and \(\sigma_a\) of \(\Sigma\).
We use the same symbols for their classes in \(\MCG^+(\Sigma)\).
The lift of \(\iota\) is represented by
\[
R=
\begin{pmatrix}
1&0\\
0&-1
\end{pmatrix}.
\]
The map \(\sigma_a\) carries the pairing \(\{\{x_1,x_2\},\{x_3,x_4\}\}\) to \(\{\{x_1,x_3\},\{x_2,x_4\}\}\).
We compose maps and mapping classes from right to left.
\end{definition}

Every class in \(\MCG^+(\Sigma)\) has a representative induced by an affine map of the torus double cover of the form
\[
x\longmapsto Mx+\eta,
\qquad M\in\mathrm{SL}(2,\mathbb Z),\quad
\eta\in\tfrac12\mathbb Z^2/\mathbb Z^2,
\]
where \(M\) is determined up to sign~\cite[Proposition~2.7]{farbmargalit2012primer}.
Mapping classes may permute \(x_1,\ldots,x_4\).

\begin{figure}[!htbp]
\centering

\begin{minipage}[t]{0.49\linewidth}
\centering
\(\tau_a\)
\par\smallskip
\includegraphics[height=25mm]{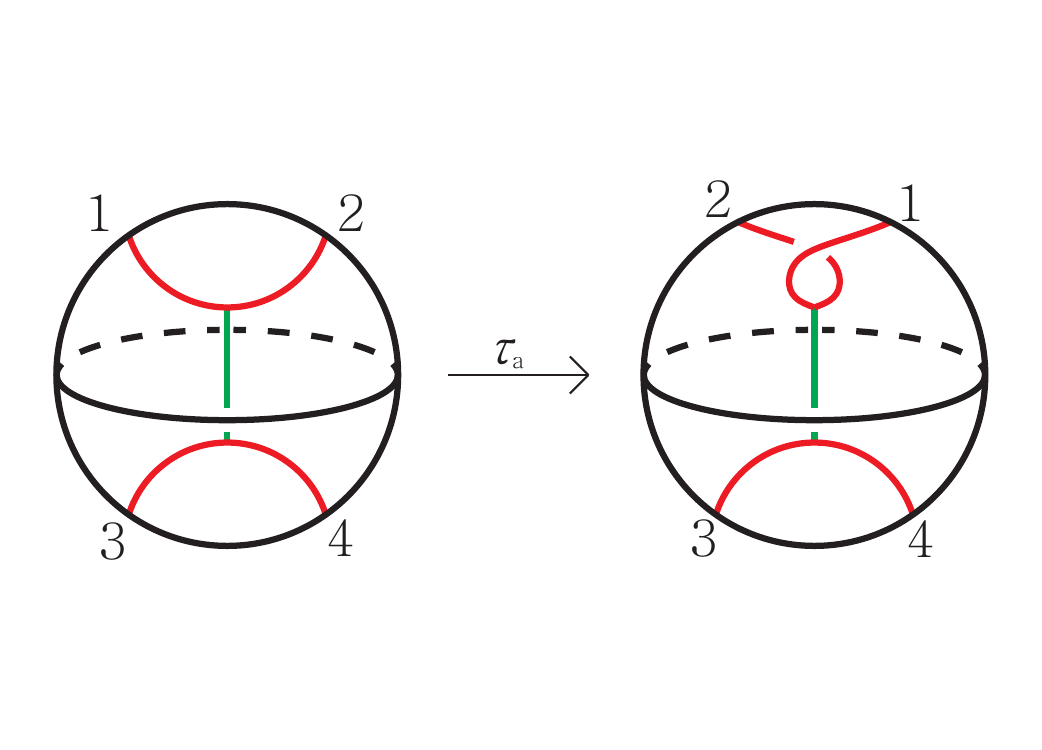}
\end{minipage}
\hfill
\begin{minipage}[t]{0.49\linewidth}
\centering
\(\tau_b\)
\par\smallskip
\includegraphics[height=25mm]{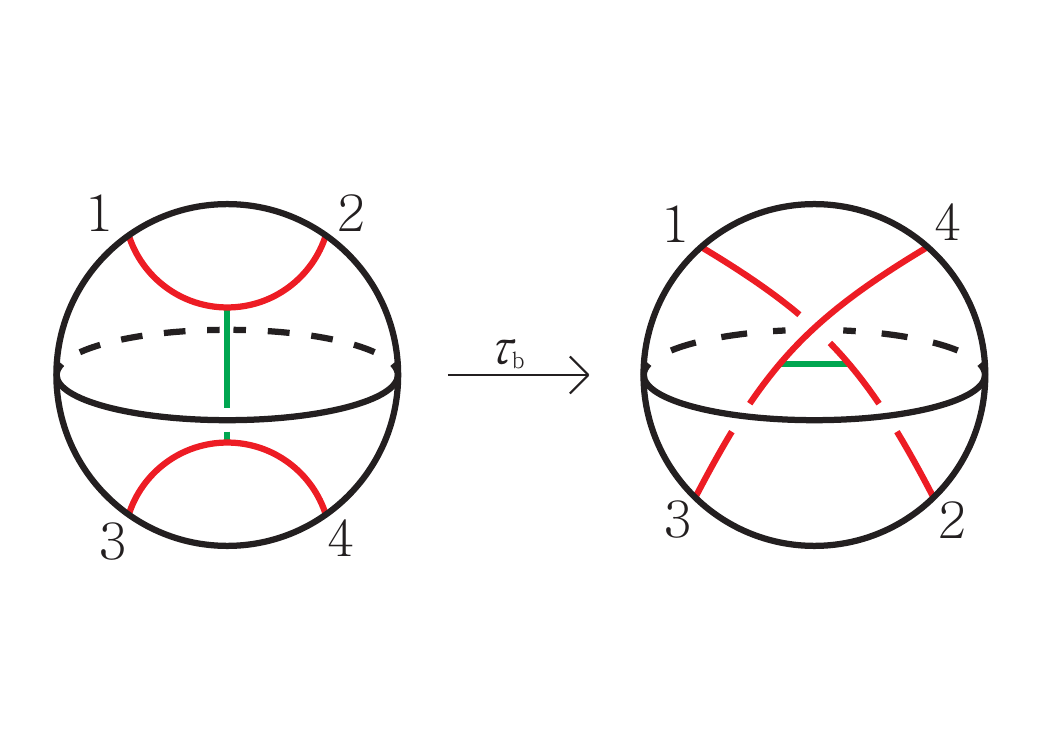}
\end{minipage}

\smallskip

\begin{minipage}[t]{0.49\linewidth}
\centering
\(\sigma_a\)
\par\smallskip
\includegraphics[height=25mm]{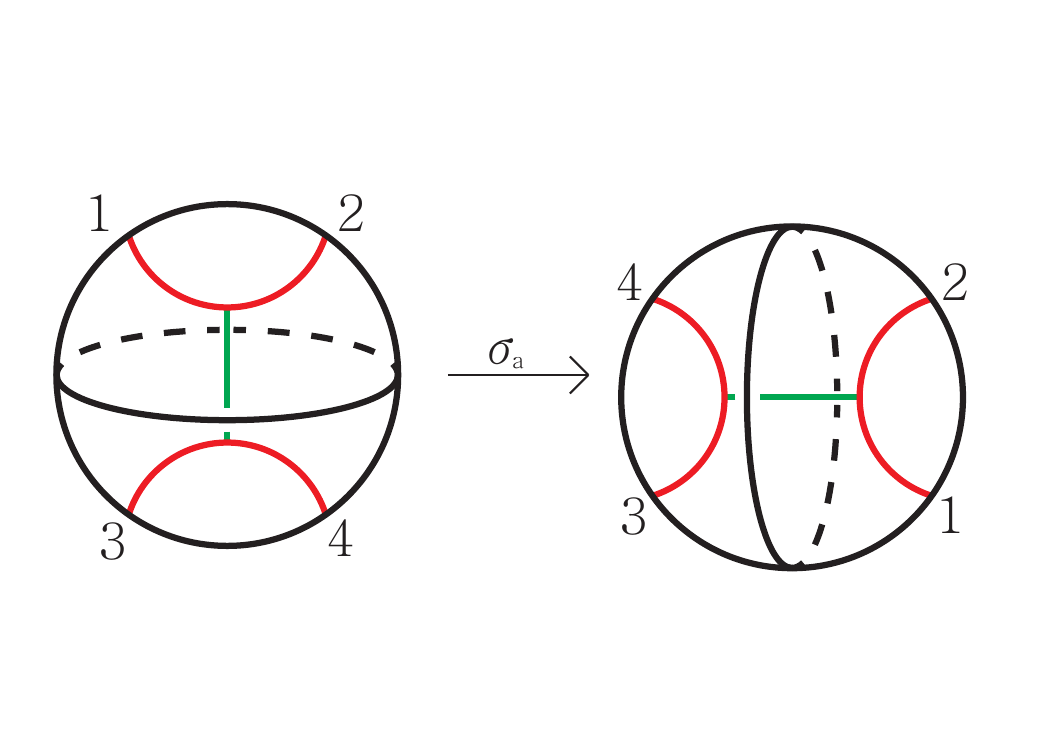}
\end{minipage}

\caption{Representatives of the mapping classes \(\tau_a\), \(\tau_b\), and \(\sigma_a\) on the pillow \(\Sigma\).}\label{fig:mcg-S2222}
\end{figure}

The slope-\(0\) trivial tangle is marked by the identity boundary parametrization.
The slope-\(\infty\) trivial tangle is marked by the boundary parametrization \(\sigma_a\) relative to the slope-\(0\) model.
Figure~\ref{fig:mcg-S2222} shows representatives of \(\tau_a\), \(\tau_b\), and \(\sigma_a\).

Every diagram of a solid pillow uses the ordering \(x_1,\ldots,x_4\) and the boundary coordinates of Definition~\ref{def:pillow-boundary-coordinates}.

\begin{figure}[htbp]
\centering
\begin{tabular}{@{}c@{\hspace{16mm}}c@{}}
\includegraphics[height=35mm]{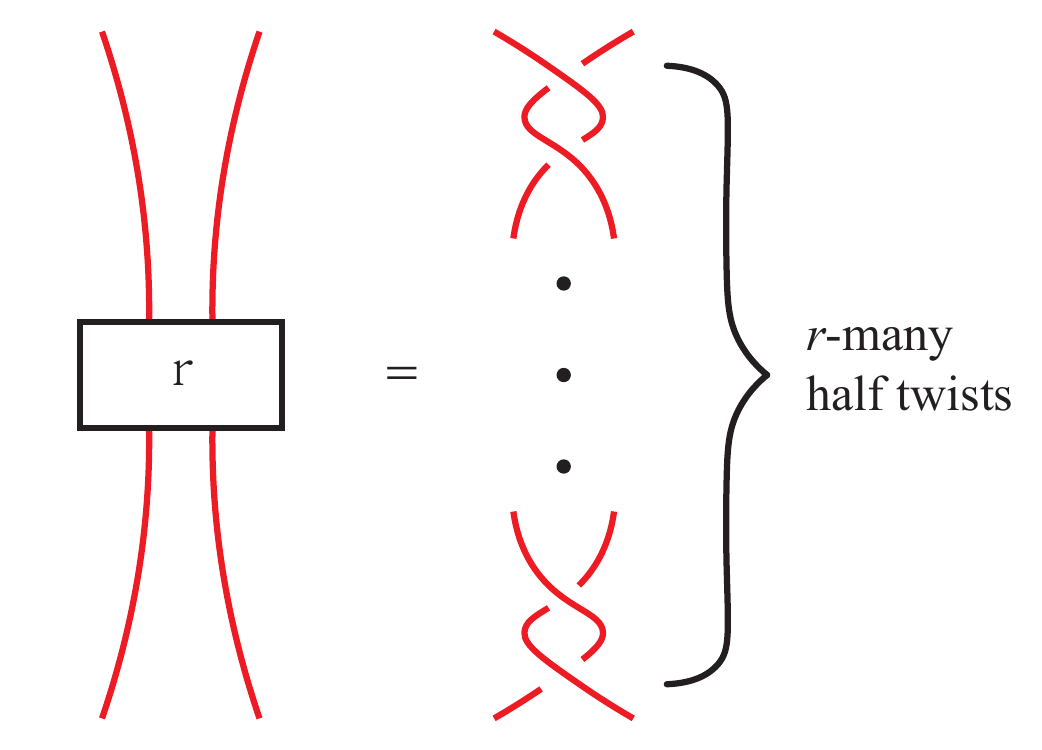}
&
\includegraphics[height=35mm]{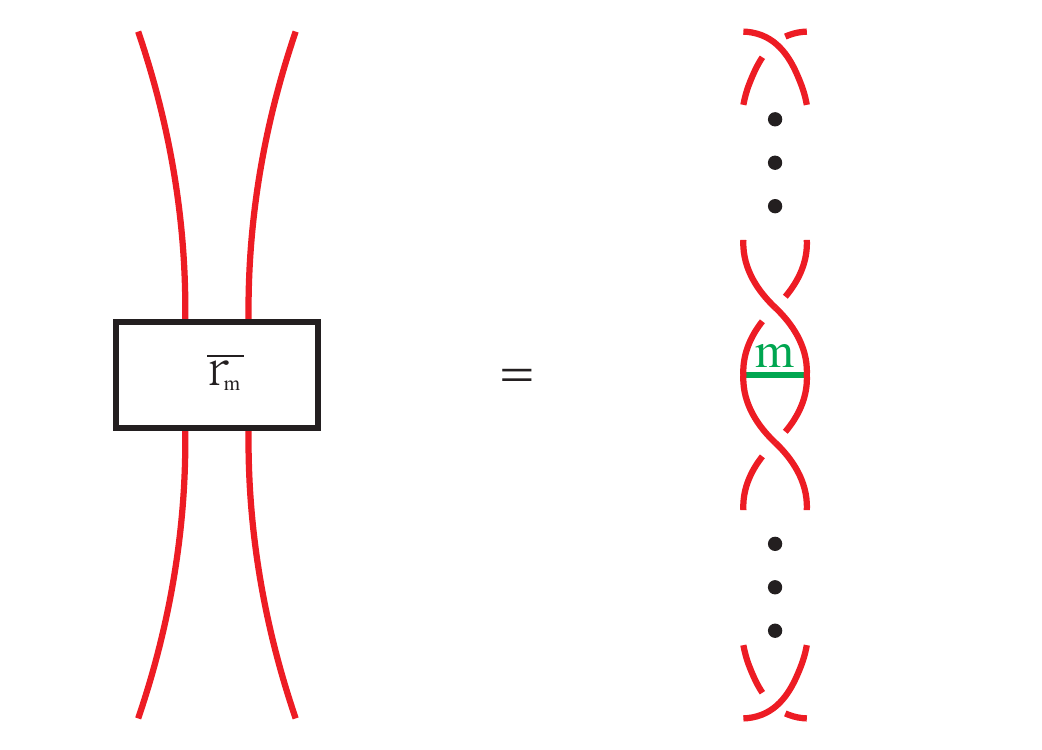}
\\[-1mm]
{\small\textup{(a) \(r\) signed half-twists}} &
{\small\textup{(b) The box \(\overline r_m\)}}
\end{tabular}
\caption{Half-twists and edge labels.}\label{fig:half-twist}
\end{figure}

\begin{definition}\label{def:half-twist-box-convention}
Let \(r\in\mathbb Z\).
A rectangular box labeled \(r\), through which two strands pass, represents \(r\) signed half-twists.
If \(r>0\), the half-twists are right-handed.
If \(r<0\), the half-twists are left-handed.
If \(r=0\), the strands pass through the box without a half-twist.
Figure~\ref{fig:half-twist} shows this convention.

For an integer \(m\geq2\), a box labeled \(\overline r_m\) contains an edge of order \(m\) joining the two strands in addition to the \(r\) signed half-twists.
\end{definition}

We use Dunbar's description of Seifert orbifolds by rational tangles~\cite{dunbar1981fibered}.
Mecchia and Seppi recall this construction in~\cite[Section~5.2, p.~838]{mecchia2015fibered}.
We first express each local Seifert model in the fixed boundary coordinates.
We then compute the gluing for a base disk with two corner reflectors and no cone points.
The Euler number and the two local invariants determine the exponent of \(\tau_b\) in the boundary gluing map.

We use the convention of Mecchia-Seppi for local invariants.
Dunbar used the opposite sign~\cite[p.~818]{mecchia2015fibered}.

\begin{definition}\label{def:dunbar-tangle-slope}
Let \(x\) be a corner reflector in the base of an oriented Seifert fibration.
Write its local invariant in the convention of Mecchia-Seppi as \(\lambda=p_\lambda/q_\lambda\), where \(p_\lambda,q_\lambda\in\mathbb Z\) and \(q_\lambda>0\).
The denominator \(q_\lambda\) is the corner order of \(x\).
The fraction \(p_\lambda/q_\lambda\) need not be reduced.
Set \(c=\gcd(p_\lambda,q_\lambda)\).
The interior of the exceptional interval fiber over \(x\) has isotropy order \(c\).

Let \(\beta\) be the unique rational number such that \(\beta\equiv-\lambda\pmod1\) and \(-1/2<\beta\leq1/2\).
The parameter \(\beta\) uses Dunbar's sign convention.
Set \(s=-\beta\).
Then \(-1/2\leq s<1/2\) and \(s\equiv\lambda\pmod1\).
\end{definition}

For each integer \(c\geq1\), the product fibration on \(\mathcal P_c\) has local invariant \(0/c\) and rational tangle slope \(0\).

\begin{definition}\label{def:continued-fraction-boundary-word}
Set \(t(0)=\Id\).
For \(s\in[-1/2,1/2)\cap\mathbb Q\) with \(s\neq0\), write
\[
\frac{1}{\lvert s\rvert}
=
a_1+\cfrac{1}{a_2+\cfrac{1}{\cdots+\cfrac{1}{a_\ell}}}
\]
as the unique finite continued fraction with
\(a_i\in\mathbb Z_{>0}\), \(a_1\geq2\), and \(a_\ell\geq2\).
Set \(\varepsilon=\operatorname{sgn}(s)\).
Define \(t(s)\in\MCG^+(\Sigma)\) by
\[
t(s)=
\begin{cases}
\tau_a^{\varepsilon a_1}\tau_b^{-\varepsilon a_2}
\cdots\tau_a^{\varepsilon a_\ell}\sigma_a,
& \ell\text{ odd},\\[2mm]
\tau_a^{\varepsilon a_1}\tau_b^{-\varepsilon a_2}
\cdots\tau_b^{-\varepsilon a_\ell},
& \ell\text{ even}.
\end{cases}
\]
The twist factors alternate between powers of \(\tau_a\) and \(\tau_b\),
with one factor for each coefficient \(a_i\) in increasing order of \(i\).
For \(\ell=1\), the formula means
\(t(s)=\tau_a^{\varepsilon a_1}\sigma_a\).
For \(\ell=2\), it means
\(t(s)=\tau_a^{\varepsilon a_1}\tau_b^{-\varepsilon a_2}\).
The rightmost factor acts first.
\end{definition}

For odd \(\ell\), the terminal factor \(\sigma_a\) specifies the boundary parametrization of the initial slope-\(\infty\) tangle.
The pairing of the four cone points alone does not determine this factor.

\begin{lemma}\label{lem:pillow-continued-fraction-word}
Write \(s=p/q\in[-1/2,1/2)\cap\mathbb Q\) in lowest terms with \(q>0\).
The product of the fixed torus lifts of the factors defining \(t(s)\) sends \(a\) to \(\pm(qa+pb)\).
\end{lemma}

\begin{proof}
For \(s=0\), one has \(p=0\), \(q=1\), and \(t(s)=\Id\).
Suppose that \(s\neq0\), and use its continued fraction from Definition~\ref{def:continued-fraction-boundary-word}.
On slopes, the fixed lifts satisfy
\[
\widetilde\tau_a^{\,u}(x)=\frac{1}{u+1/x},
\qquad
\widetilde\tau_b^{-u}(x)=u+x,
\qquad
\widetilde\sigma_a(0)=\infty
\quad(u\in\mathbb Z).
\]
Applying the factors from right to left sends slope \(0\) to
\[
\cfrac{1}{\varepsilon a_1+
\cfrac{1}{\varepsilon a_2+
\cfrac{1}{\cdots+\cfrac{1}{\varepsilon a_\ell}}}}
=\frac{p}{q}.
\]
The image of \(a\) is primitive because the lift acts by an integral matrix of determinant one.
A primitive class of slope \(p/q\) is \(\pm(qa+pb)\).
\end{proof}

\begin{lemma}\label{lem:dunbar-boundary-word}
Let \(c\geq1\) be an integer.
Write \(s=p/q\in[-1/2,1/2)\cap\mathbb Q\) in lowest terms with \(q>0\).
The solid pillow \(\mathcal P_c\) admits a Seifert fibration with local invariant \(cp/cq\) in the convention of Definition~\ref{def:dunbar-tangle-slope}.
There is an orientation-preserving boundary marking
\(\mathfrak m\colon\Sigma\to\orbpartial\mathcal P_c\) such that
\[
[\mathfrak m]=t(s)^{-1}
\]
in the fixed boundary coordinates of Definition~\ref{def:pillow-boundary-models}.
In the torus coordinates determined by \(\mathfrak m\),
the regular boundary fibers have slope \(\infty\).
\end{lemma}

\begin{proof}
Use the products of the fixed linear representatives in Definition~\ref{def:pillow-boundary-coordinates} to represent \(t(s)\).
Let \(\mathfrak m\colon\Sigma\to\orbpartial\mathcal P_c\) be the boundary marking represented by \(t(s)^{-1}\) in the fixed boundary coordinates.
By Lemma~\ref{lem:pillow-continued-fraction-word}, the first column of a torus lift of \(t(s)\) is \((q,p)^T\), after changing its sign if necessary.
Changing this sign does not change the map of \(\Sigma\).
The marked meridian is therefore \(qa+pb\).

Write the regular fiber in the unmarked coordinates as \(\alpha a+qb=(\mathfrak m)_*b\).
Since the boundary lift has determinant one, \(\gcd(\alpha,q)=1\) and \(\alpha p\equiv1\pmod q\).
On the solid torus with singular core of order \(c\), the circle orbits are induced by
\[
(z,v)\longmapsto
\bigl(e^{i\alpha t/c}z,e^{iqt}v\bigr),
\qquad t\in\mathbb R,
\]
where rotation of \(z\) through \(2\pi/c\) is identified by the local cyclic action.
This action has period \(2\pi\) on the quotient.
Complex conjugation in both factors takes each orbit to an orbit and reverses its parameter.
Consequently the fibration induces a Seifert fibration of \(\mathcal P_c\).
The map from \(D^2\times S^1\) to the cyclic solid-torus quotient given by
\[
(w,e^{i\theta})\longmapsto
\bigl[e^{i\alpha\theta/c}w,e^{iq\theta}\bigr]
\]
is unchanged by
\[
(w,e^{i\theta})\longmapsto
\bigl(e^{2\pi i/(cq)}w,e^{i(\theta-2\pi p/q)}\bigr),
\]
because \(\alpha p\equiv1\pmod q\).
It identifies the quotient by this cyclic group of order \(cq\) with the solid torus of core isotropy \(c\).
The generator rotates the base by \(2\pi/(cq)\) and the fiber by \(-2\pi cp/(cq)\).
The quotient by complex conjugation in both factors is the local model of the interval fiber.
Thus its local invariant is \(cp/cq\) in the convention of~\cite[pp.~818--819]{mecchia2015fibered}.
The interior of the interval fiber has isotropy order \(c\).
When \(cq=1\), this fiber lies over a marked regular mirror point.
\end{proof}

The rational tangle construction of Dunbar realizes these meridian slopes by continued fraction diagrams~\cite[Section~4]{dunbar1988geometric}.
We use the full boundary marking in Lemma~\ref{lem:dunbar-boundary-word}, including the terminal factor \(\sigma_a\) when the continued fraction has odd length.
The slope alone does not specify this marking.

Let \(\mathcal O\) be an oriented Seifert \(3\)-orbifold whose base is a disk with two ordered corner reflectors and no cone points.
Choose a smooth embedded arc in the base separating the two corners, with endpoints in the interiors of mirror edges.
In a collar of each mirror edge, require the arc to have constant boundary coordinate.
Let \(A\) and \(B\) be the inverse images of the two closed regions bounded by this arc.
Let \(e\) be the Euler number.
Let \(\beta_A,\beta_B\) be the normalized local invariants in Definition~\ref{def:dunbar-tangle-slope}.
For \(X\in\{A,B\}\), set \(s_X=-\beta_X\).
Let \(c_X\) be the isotropy order of the exceptional interval fiber in \(X\).

\begin{definition}\label{def:family-11-central-box-labels}
Set
\[
\kappa=-2e+\beta_A+\beta_B=-2e-s_A-s_B.
\]
This is an integer by the compatibility relation for Seifert invariants~\cite[Proposition~3.4]{mecchia2015fibered}.
Construct the two marked pillows of Lemma~\ref{lem:dunbar-boundary-word} with parameters \((c_A,s_A)\) and \((c_B,s_B)\).
In their marked boundary coordinates, glue the second pillow to the first by \(\tau_b^\kappa\iota\).
We use the markings on \(A\) and \(B\) obtained from this construction by the fiber-preserving identification in the proof of Proposition~\ref{lem:ordered-two-corner-gluing}.
The integer \(\kappa\) is the \emph{central-box label} in Dunbar's convention and is the negative of the label in Definition~\ref{def:half-twist-box-convention}.
\end{definition}

\begin{proposition}\label{lem:ordered-two-corner-gluing}
Let \(\mathcal O\) be an oriented Seifert \(3\)-orbifold whose base is a disk with two ordered corner reflectors and no cone points.
There is a fiber-preserving, orientation-preserving orbifold diffeomorphism from the orbifold constructed in Definition~\ref{def:family-11-central-box-labels} to \(\mathcal O\).
This diffeomorphism maps the two pillows onto \(A\) and \(B\), respectively.
For the resulting model identifications \(f_X\colon X\to\mathcal P_{c_X}\), the boundary gluing map
\(\Phi\colon\orbpartial\mathcal P_{c_B}\to\orbpartial\mathcal P_{c_A}\)
satisfies
\[
[\Phi]=t(s_A)^{-1}\tau_b^\kappa\iota t(s_B).
\]
\end{proposition}

\begin{proof}
On the torus cover, the gluing map \(\tau_b^\kappa\iota\) is represented by
\[
(x,y)\longmapsto(x,-\kappa x-y).
\]
This map preserves the vertical fibers and commutes with the deck involution.
Hence the two local Seifert fibrations glue.
The base is a disk with two corner reflectors and no cone points.
By Lemma~\ref{lem:dunbar-boundary-word}, the local invariants agree with those of \(\mathcal O\).

Dunbar's Euler-number formula~\cite[pp.~74--75]{dunbar1988geometric} gives
\[
-\frac{\kappa+s_A+s_B}{2}=e.
\]
By Corollary~3.2 of~\cite{mecchia2015fibered}, there is a fiber-preserving, orientation-preserving orbifold diffeomorphism \(F\) from the constructed orbifold to \(\mathcal O\).

If \(F\) interchanges the corner fibers, then \(s_A=s_B\) and \(c_A=c_B\).
The two model pillows are therefore copies of the same fibered orbifold.
Since the fixed representative of \(\tau_b^\kappa\iota\) is an involution, exchanging these copies defines a fiber-preserving, orientation-preserving orbifold diffeomorphism of their union.
Precompose \(F\) with this diffeomorphism if necessary.
Then \(F\) preserves the prescribed corner order.

The image of the common boundary of the model pillows under \(F\) projects to a separating arc in the base of \(\mathcal O\).
This arc is isotopic to the arc defining \(A\) and \(B\), relative to neighborhoods of the corner reflectors.
Lift this isotopy to a fiber-preserving orbifold isotopy of \(\mathcal O\) and compose it with \(F\).
We may therefore assume that \(F\) maps the two model pillows onto \(A\) and \(B\), respectively.

For \(X\in\{A,B\}\), set
\[
f_X=(F|_{\mathcal P_{c_X}})^{-1}
\colon X\longrightarrow\mathcal P_{c_X}.
\]
Let \(\mu_X\colon\Sigma\to\orbpartial\mathcal P_{c_X}\) be the boundary marking used in the construction.
By Lemma~\ref{lem:dunbar-boundary-word},
\[
[\mu_X]=t(s_X)^{-1}.
\]
The boundary gluing map in these model coordinates therefore satisfies
\[
\begin{aligned}
[\Phi]
&=[\mu_A]\,\tau_b^\kappa\iota\,[\mu_B]^{-1}\\
&=t(s_A)^{-1}\tau_b^\kappa\iota t(s_B).
\end{aligned}
\]
\end{proof}

\section{Boundary regions and their gluing}\label{sec:boundary-regions}

In this section, we identify the attaching and belt regions for the ordered handle classes of indices \(1\), \(2\), and \(3\) classified in Part~I~\cite{choilee2026part1}.
We compute their boundary gluing maps.
We also characterize diffeomorphisms of attachments relative to a common exterior.

Section~\ref{subsec:gluing-map} fixes the family parameters and the conventions used for the gluing maps.
Tables~\ref{tab:index-two-models}--\ref{tab:polyhedral-regions} record the attaching and belt regions and their boundary gluing maps.
Section~\ref{subsec:family-boundary-calculations} proves Theorem~\ref{thm:computed-boundary-regions} by computing the regions and gluing maps for each family.
Appendix~\ref{subsec:illustrated-boundary-models} records the corresponding labeled diagrams.
Section~\ref{subsec:relative-handle-attachments} proves the relative diffeomorphism criterion and invariance under isotopy of the attaching embedding.

Let \(\rho\colon G\to\SO(4)\) be a faithful representation of a finite group.
Let \(\mathbb R^4=V^-\oplus V^+\) be a \(\rho(G)\)-invariant orthogonal splitting with \(1\leq\dim V^-\leq3\).
We use the notation \(H_\rho\), \(A_\rho\), \(B_\rho\), \(C_\rho\), and \(\varphi_\rho\) from Definition~\ref{def:quotient-handle-boundary-regions}.
In the family calculations, write \(A=A_\rho\) and \(B=B_\rho\).

\subsection{The computed regions and gluing maps}\label{subsec:gluing-map}

Let \(H_\rho(V^-,V^+)\) be as in Definition~\ref{def:quotient-handle-boundary-regions}.
The \textup{ordered dual} in Part~I~\cite[Definition~1.2]{choilee2026part1} has orbifold handle \(H_\rho(V^+,V^-)\).
Ordered duality interchanges the attaching and belt regions and replaces \(\Phi_\rho\) by \(\Phi_\rho^{-1}\).

We identify \(\mathbb R^4\) with \(\mathbb H\).
Let \(\mu\colon\mathbb H\to\mathbb H\) denote quaternionic conjugation, so that \(\mu(q)=\overline q\).
Set \(\rho^\mu(g)=\mu\rho(g)\mu^{-1}\).
The \textup{factor interchange} of Part~I~\cite[Definition~3.3]{choilee2026part1} replaces \(\rho\) by \(\rho^\mu\).
It sends \(V^-\) and \(V^+\) to \(\mu V^-\) and \(\mu V^+\), respectively.
Quaternionic conjugation induces orientation-reversing orbifold diffeomorphisms between the corresponding attaching regions and between the corresponding belt regions.
For a representative obtained by factor interchange, we use these diffeomorphisms to transport \(A_\rho\), \(B_\rho\), and \(\Phi_\rho\).

We use the completed Du Val family notation, admissibility conditions, and representative parameters from Part~I~\cite[Section~3.1 and Definition~3.2]{choilee2026part1} and Mecchia-Seppi~\cite[Table~1, pp.~814--815]{mecchia2015fibered}.
\textup{Families} \(1\) and \(11\) have parameters \((m,n,r,s)\).
\textup{Families} \(1'\) and \(11'\) have parameters \((m,n,r',s)\).
For the primed families, the parameter \(r\) of Mecchia-Seppi equals \(2r'\).
The parameters satisfy \(\gcd(s,r)=1\) in the unprimed families and \(\gcd(s,2r')=1\) in the primed families.
In \textup{Families} \(1\) and \(11\), we choose an odd representative of \(s\) modulo \(r\).
If \(s\) is even, then \(r\) is odd.
Replacing \(s\) by \(s+r\) does not change the isomorphism class of the quotient.

The following table lists the notation and parameter conditions for the cyclic subfamilies in Part~I~\cite[Definition~3.10 and Propositions~A.5--A.6]{choilee2026part1}.
\[
\begin{array}{c|c}
\textnormal{subfamily}&\textnormal{parameter conditions}\\
\hline
1_{\mathrm{(iii)}}&m=n=1,\quad (1-s)/r\in\mathbb Z\\
1_{\mathrm{(iv)}}&m=n=1,\quad (1+s)/r\in\mathbb Z\\
1'_{\mathrm{(iii)}}&m=n=1,\quad (1-s)/r'\in\mathbb Z\\
1'_{\mathrm{(iv)}}&m=n=1,\quad (1+s)/r'\in\mathbb Z.
\end{array}
\]
For \textup{Families} \(11\) and \(11'\), these conditions apply to the abelian subgroups \(G_1\subset G_{11}\) and \(G_{1'}\subset G_{11'}\), respectively.
Part~I proves that \(G_{11}\) has an invariant line precisely when \(G_1\) belongs to \textup{Family} \(1_{\mathrm{(iii)}}\) or \(1_{\mathrm{(iv)}}\)~\cite[Lemma~3.14]{choilee2026part1}.
It also proves that \(G_{11'}\) has an invariant line precisely when \(G_{1'}\) belongs to \textup{Family} \(1'_{\mathrm{(iii)}}\) or \(1'_{\mathrm{(iv)}}\)~\cite[Lemma~3.14]{choilee2026part1}.

Assume that \(r'>2\).
Then the two primed subfamilies are disjoint.
Set \(\ell'=(1-s)/r'\) in \textup{Family} \(1'_{\mathrm{(iii)}}\) and \(\ell'=(1+s)/r'\) in \textup{Family} \(1'_{\mathrm{(iv)}}\).
The notation \(\ell'\) is defined only for \(r'>2\).

Suppose that \(\ell'\) is odd.
Since \(\gcd(s,2r')=1\), the integer \(s\) is odd.
The definition of \(\ell'\) then implies that \(r'\) is even.
Write \(r'=2r''\).
Since \(r'>2\), we have \(r''\geq2\).

\noindent\textbf{Gluing conventions.}
We use the boundary identification of Section~\ref{subsec:corner-smoothing}.

The calculations specify orientation-preserving orbifold diffeomorphisms from \(A_\rho\) and \(B_\rho\) to the models in Tables~\ref{tab:index-two-models}--\ref{tab:polyhedral-regions}.
Under these identifications, \(\Phi_\rho\) is the boundary gluing map defined in Definition~\ref{def:boundary-gluing-models}.

For torus boundary components, we express \((\Phi_\rho)_*\) in the meridian--longitude bases of Definition~\ref{def:torus-boundary-coordinates}.
For pillow boundary components, we use the identification with \(\Sigma\) from Definition~\ref{def:pillow-boundary-models} and the coordinates of Definition~\ref{def:pillow-boundary-coordinates}.
We express \([\Phi_\rho]\) using elements of \(\MCG^+(\Sigma)\) and \(\iota\).

For index-\(1\) handles, we use the boundary maps in Lemma~\ref{lem:line-region-models}.
Equation~\eqref{eq:line-gluing-connected} applies when the common boundary is connected.
Equation~\eqref{eq:line-gluing-disconnected} applies when it has two components.
The map in Equation~\eqref{eq:line-gluing-disconnected} preserves the \(S^0\)-coordinate.
These maps specify the correspondence between the cone points.

We write \(\mathbf i,\mathbf j,\mathbf k\) for the quaternion units.
All spans are taken over \(\mathbb R\).
We identify
\[
\mathbb R^4
=
\mathbb H
=
\mathbb C\oplus\mathbb C\mathbf j
=
\operatorname{span}\{1,\mathbf i\}
\oplus
\operatorname{span}\{\mathbf j,\mathbf k\}.
\]
We identify the two summands with \(\mathbb R^2\) using the ordered bases \((1,\mathbf i)\) and \((\mathbf j,\mathbf k)\), respectively.
Part~I defines, for \(\theta\in\mathbb R\),
\[
M_\theta
=
\begin{bmatrix}
\cos\theta&\sin\theta\\
\sin\theta&-\cos\theta
\end{bmatrix},
\qquad
W_\theta
=
\{(v,M_\theta v)\mid v\in\mathbb R^2\}
\]
\cite[Definition~3.17]{choilee2026part1}.

For \textup{Family} \(1'\) with \(m=n=1\) and \(r'=2\), write \(\mathbb R^4=U_{\mathrm{triv}}\oplus U_{\mathrm{nontriv}}\).
Part~I shows that \(G_{1'}\) acts trivially on \(U_{\mathrm{triv}}\) and by the nontrivial real character on \(U_{\mathrm{nontriv}}\)~\cite[Section~4.2]{choilee2026part1}.

For \textup{Family} \(11'\) with \(r'>2\) and even \(\ell'\), we use the representative chosen in Part~I~\cite[proof of Proposition~4.4]{choilee2026part1}.
In these coordinates, set \(L_0=\operatorname{span}\{1\}\) and \(L_1=\operatorname{span}\{\mathbf i\}\).
The group \(G_{11'}\) fixes \(L_0\) pointwise and acts nontrivially on \(L_1\).

For \textup{Families} \(1\), \(1'\), \(11\), and \(11'\), set \(h=\gcd(m,n)\), \(m_0=m/h\), and \(n_0=n/h\).
We adopt the definition of \(a\) from Mecchia-Seppi~\cite[pp.~826 and~830]{mecchia2015fibered}.
\[
a=
\begin{cases}
\gcd(n_0-sm_0,n_0+sm_0,2m_0n_0r),
&\textup{for Families }1,11,\\
\gcd(n_0-sm_0,n_0+sm_0,2m_0n_0r'),
&\textup{for Families }1',11'.
\end{cases}
\]

For \textup{Families} \(1\) and \(11\), we use \(b_1,b_2,f_1,f_2,f,d,g,\overline g,\overline c\) from Mecchia-Seppi's Table~3~\cite[Table~3, p.~832]{mecchia2015fibered}.
For \textup{Families} \(1'\) and \(11'\), we use \(b_1,b_2,f,d,g,\overline g,\overline c\) from Mecchia-Seppi's Table~2~\cite[Table~2, p.~828]{mecchia2015fibered}.
In both cases, \(a\) is defined as above.
The symbols \(\overline g\) and \(\overline c\) denote the modular inverses specified in these tables.
In particular, \(g\overline g\equiv1\pmod f\).

For \(G\in\{G_{11},G_{11'}\}\), let \(\lambda_A^{\mathrm{loc}}\) and \(\lambda_B^{\mathrm{loc}}\) denote the local invariants of the Seifert fibration of \(S^3\sslash G\) corresponding to \(A\) and \(B\), respectively.
Let \(e_G\) denote the Euler number of this fibration.
Mecchia and Seppi record the following values in Tables~2 and~3~\cite[Tables~2 and~3, pp.~828 and~832]{mecchia2015fibered}.
\[
(\lambda_A^{\mathrm{loc}},\lambda_B^{\mathrm{loc}},e_G)
=
\begin{cases}
\left(
-\dfrac{g\overline c f_1b_1h}{nr},
\dfrac{d\overline c f_2b_2h}{nr},
-\dfrac{m}{nr}
\right),
& G=G_{11},\\[3mm]
\left(
-\dfrac{g\overline c b_1h}{nr'},
\dfrac{d\overline c b_2h}{nr'},
-\dfrac{m}{2nr'}
\right),
& G=G_{11'}.
\end{cases}
\]

For \(G\in\{G_{11},G_{11'}\}\), let \(s_A\) and \(s_B\) be the slopes associated with \(\lambda_A^{\mathrm{loc}}\) and \(\lambda_B^{\mathrm{loc}}\) by Definition~\ref{def:dunbar-tangle-slope}.
Let \(\kappa_G\) denote the integer \(\kappa\) in Definition~\ref{def:family-11-central-box-labels} with \(e=e_G\).

For \(G\in\{G_1,G_{1'}\}\), let \(q_G\in\{0,\ldots,f-1\}\) be the unique integer such that \(q_G\equiv d\overline g\pmod f\).
Then \(\gcd(f,q_G)=1\).
Let \(v_G\in\{0,\ldots,f-1\}\) be the unique integer such that \(q_Gv_G\equiv1\pmod f\).
Set \(u_G=(1-q_Gv_G)/f\).

\begin{theorem}\label{thm:computed-boundary-regions}
Let \(\rho\colon G\to\SO(4)\) be a faithful representation of a nontrivial finite group.
Let \(\mathbb R^4=V^-\oplus V^+\) be a \(\rho(G)\)-invariant orthogonal splitting.
Assume that the splitting is one of the index-\(1\) or index-\(2\) representatives in Part~I, the ordered dual of an index-\(1\) representative, or a splitting obtained from one of these by factor interchange~\cite[Theorem~B and Table~1]{choilee2026part1}.

Tables~\ref{tab:index-two-models}--\ref{tab:polyhedral-regions} identify oriented labeled models for the attaching region \(A_\rho\) and the belt region \(B_\rho\).
With the fixed model identifications and boundary coordinates, these tables determine the orbifold isotopy class \([\Phi_\rho]\) of the boundary gluing map in Definition~\ref{def:boundary-gluing-models}.
\end{theorem}

The splittings in Table~\ref{tab:index-two-models} and their ordered duals represent all index-\(2\) classes.

\newcommand{\IndexTwoTableContent}{%
\begingroup
\setlength{\tabcolsep}{1.5pt}
\begin{tabular}{@{}L{0.16\linewidth}L{0.22\linewidth}L{0.11\linewidth}L{0.11\linewidth}L{0.36\linewidth}@{}}
\hline
Group and parameter range
&
Representative splitting
&
Attaching region \(A_\rho\)
&
Belt region \(B_\rho\)
&
Boundary \(C_\rho\) and gluing map
\\
\hline
\(G_1\)
&
\(\operatorname{span}\{1,\mathbf i\}\oplus\operatorname{span}\{\mathbf j,\mathbf k\}\)
&
\(T_{f_1b_1h}\)
&
\(T_{f_2b_2h}\)
&
\(T^2\), with
\(\displaystyle
(\Phi_\rho)_*
=
\begin{pmatrix}
-q_G&u_G\\
f&v_G
\end{pmatrix}\)
\\
\(G_1\), \(m=n=1\), \(r=2\)
&
\(\operatorname{span}\{1,\mathbf j\}\oplus\operatorname{span}\{\mathbf i,\mathbf k\}\)
&
\(\mathcal P_1\)
&
\(\mathcal P_1\)
&
\(S^2(2,2,2,2)\), with
\([\Phi_\rho]=\tau_b^2\iota\)
\\
\(G_{1'}\)
&
\(\operatorname{span}\{1,\mathbf i\}\oplus\operatorname{span}\{\mathbf j,\mathbf k\}\)
&
\(T_{b_1h}\)
&
\(T_{b_2h}\)
&
\(T^2\), with
\(\displaystyle
(\Phi_\rho)_*
=
\begin{pmatrix}
-q_G&u_G\\
f&v_G
\end{pmatrix}\)
\\
\(G_{1'}\), \(m=n=1\), \(r'=2\)
&
\(\operatorname{span}\{1,\mathbf j\}\oplus\operatorname{span}\{\mathbf i,\mathbf k\}\)
&
\(\mathcal P_1\)
&
\(\mathcal P_1\)
&
\(S^2(2,2,2,2)\), with
\([\Phi_\rho]=\sigma_a\iota\)
\\
\(G_{11}\)
&
\(\operatorname{span}\{1,\mathbf i\}\oplus\operatorname{span}\{\mathbf j,\mathbf k\}\)
&
\(\mathcal P_{f_1b_1h}\)
&
\(\mathcal P_{f_2b_2h}\)
&
\(S^2(2,2,2,2)\), with
\([\Phi_\rho]=t(s_A)^{-1}\tau_b^{\kappa_{G_{11}}}\iota t(s_B)\)
\\
\(G_{11'}\)
&
\(\operatorname{span}\{1,\mathbf i\}\oplus\operatorname{span}\{\mathbf j,\mathbf k\}\)
&
\(\mathcal P_{b_1h}\)
&
\(\mathcal P_{b_2h}\)
&
\(S^2(2,2,2,2)\), with
\([\Phi_\rho]=t(s_A)^{-1}\tau_b^{\kappa_{G_{11'}}}\iota t(s_B)\)
\\
\(G_3\), \(m=1\), \(n\geq2\)
&
\(W_\theta\oplus W_\theta^\perp\)
&
\(\mathcal P_1\)
&
\(\mathcal P_1\)
&
\(S^2(2,2,2,2)\), with
\([\Phi_\rho]=\tau_b^{-2n}\iota\)
\\
\(G_{34}\), \(m=1\), \(n\geq3\) odd
&
\(W_\theta\oplus W_\theta^\perp\)
&
\(\mathcal P_1\)
&
\(\mathcal P_1\)
&
\(S^2(2,2,2,2)\), with
\([\Phi_\rho]=\tau_b^{-n}\iota\)
\\
\hline
\end{tabular}
\endgroup
}

For the two solid-torus rows of Table~\ref{tab:index-two-models}, a Dehn twist about a meridional disk extends over \(T_q\).
It changes the boundary longitude by an integral multiple of the meridian.
We choose the model identifications so that \((\Phi_\rho)_*\) is represented by the matrix specified in the table.

Tables~\ref{tab:cyclic-line-regions} and~\ref{tab:dihedral-line-regions} describe the subfamilies with an invariant line.
In these tables, \(m=n=1\).
For \(G_1\) and \(G_{11}\), the parameters satisfy \((1-s)/r\in\mathbb Z\) or \((1+s)/r\in\mathbb Z\).
For \(G_{1'}\) and \(G_{11'}\), they satisfy \((1-s)/r'\in\mathbb Z\) or \((1+s)/r'\in\mathbb Z\).

\newcommand{\CyclicLineTableContent}{%
\begin{tabular}{@{}L{0.16\linewidth}L{0.14\linewidth}L{0.18\linewidth}L{0.25\linewidth}L{0.20\linewidth}@{}}
\hline
Group and parameter range
&
Orbit of invariant lines
&
Attaching region \(A_\rho\)
&
Belt region \(B_\rho\)
&
Boundary \(C_\rho\) and gluing map
\\
\hline
\(G_1\), \(r\geq1\)
&
unique orbit
&
\(N_{\Cyc{r}}\)
&
\(\mathcal M(\Cyc{r},\alpha)\)
&
\(\Sigma_{\Cyc{r}}\), \eqref{eq:line-gluing-connected}
\\
\(G_{1'}\), \(r'>2\), \(\ell'\) even
&
unique orbit
&
\(S^0\times N_{\Cyc{r'}}\)
&
\(I(\Sigma_{\Cyc{r'}})\)
&
\(S^0\times\Sigma_{\Cyc{r'}}\), \eqref{eq:line-gluing-disconnected}
\\
\(G_{1'}\), \(r'>2\), \(\ell'\) odd
&
unique orbit
&
\(N_{\Cyc{r''}}\)
&
\(\mathcal M(\Cyc{r''},\delta_{r''})\)
&
\(\Sigma_{\Cyc{r''}}\), \eqref{eq:line-gluing-connected}
\\
\(G_{1'}\), \(r'=2\)
&
lines in \(U_{\mathrm{triv}}\)
&
\(S^0\times N_{\Cyc{2}}\)
&
\(I(\Sigma_{\Cyc{2}})\)
&
\(S^0\times\Sigma_{\Cyc{2}}\), \eqref{eq:line-gluing-disconnected}
\\
\(G_{1'}\), \(r'=2\)
&
lines in \(U_{\mathrm{nontriv}}\)
&
\(N_{\Cyc{1}}\)
&
\(\mathcal M(\Cyc{1},\delta_1)\)
&
\(\Sigma_{\Cyc{1}}\), \eqref{eq:line-gluing-connected}
\\
\hline
\end{tabular}
}

\newcommand{\DihedralLineTableContent}{%
\begin{tabular}{@{}L{0.16\linewidth}L{0.14\linewidth}L{0.18\linewidth}L{0.25\linewidth}L{0.20\linewidth}@{}}
\hline
Group and parameter range
&
Orbit of invariant lines
&
Attaching region \(A_\rho\)
&
Belt region \(B_\rho\)
&
Boundary \(C_\rho\) and gluing map
\\
\hline
\(G_{11}\)
&
unique orbit
&
\(N_{\Dih{2r}}\)
&
\(\mathcal M(\Dih{2r},\alpha)\)
&
\(\Sigma_{\Dih{2r}}\), \eqref{eq:line-gluing-connected}
\\
\(G_{11'}\), \(r'>2\), \(\ell'\) even
&
\(L_0\)
&
\(S^0\times N_{\Dih{2r'}}\)
&
\(I(\Sigma_{\Dih{2r'}})\)
&
\(S^0\times\Sigma_{\Dih{2r'}}\), \eqref{eq:line-gluing-disconnected}
\\
\(G_{11'}\), \(r'>2\), \(\ell'\) even
&
\(L_1\)
&
\(N_{\Cyc{r'}}\)
&
\(\mathcal M(\Cyc{r'},\gamma)\)
&
\(\Sigma_{\Cyc{r'}}\), \eqref{eq:line-gluing-connected}
\\
\(G_{11'}\), \(r'>2\), \(\ell'\) odd
&
unique orbit
&
\(N_{\Dih{2r''}}\)
&
\(\mathcal M(\Dih{2r''},\delta_{r''})\)
&
\(\Sigma_{\Dih{2r''}}\), \eqref{eq:line-gluing-connected}
\\
\(G_{11'}\), \(r'=2\)
&
fixed lines
&
\(S^0\times N_{\Dih{4}}\)
&
\(I(\Sigma_{\Dih{4}})\)
&
\(S^0\times\Sigma_{\Dih{4}}\), \eqref{eq:line-gluing-disconnected}
\\
\(G_{11'}\), \(r'=2\)
&
lines on which \(G_{11'}\) acts nontrivially
&
\(N_{\Dih{2}}\)
&
\(\mathcal M(\Dih{2},\delta_1)\)
&
\(\Sigma_{\Dih{2}}\), \eqref{eq:line-gluing-connected}
\\
\hline
\end{tabular}
}
In the second column of Table~\ref{tab:polyhedral-regions}, we write \(\mathbf J\) for \(\mathbf I\), following Dunbar~\cite[Figures~3 and~10, pp.~129 and~139]{dunbar1994nonfibering}.

\newcommand{\PolyhedralTableContent}{%
\begin{tabular}{@{}L{0.12\linewidth}L{0.19\linewidth}L{0.17\linewidth}L{0.25\linewidth}L{0.20\linewidth}@{}}
\hline
Group
&
Ambient quotient in Dunbar's notation
&
Attaching region \(A_\rho\)
&
Belt region \(B_\rho\)
&
Boundary \(C_\rho\) and gluing map
\\
\hline
\(G_{21}\)
&
\(\mathbf T\times_{\mathbf T}\mathbf T\), Figure~3
&
\(N_{\mathbf T}\)
&
\(\mathcal M(\mathbf T,\alpha)\)
&
\(\Sigma_{\mathbf T}\), \eqref{eq:line-gluing-connected}
\\
\(G_{21'}\)
&
\(\mathbf T\times_{\mathbf T}\mathbf T^1\), Figure~10
&
\(S^0\times N_{\mathbf T}\)
&
\(I(\Sigma_{\mathbf T})\)
&
\(S^0\times\Sigma_{\mathbf T}\), \eqref{eq:line-gluing-disconnected}
\\
\(G_{26}\)
&
\(\mathbf O\times_{\mathbf O}\mathbf O\), Figure~3
&
\(N_{\mathbf O}\)
&
\(\mathcal M(\mathbf O,\alpha)\)
&
\(\Sigma_{\mathbf O}\), \eqref{eq:line-gluing-connected}
\\
\(G_{26'}\)
&
\(\mathbf O\times_{\mathbf O}\mathbf O^1\), Figure~10
&
\(S^0\times N_{\mathbf O}\)
&
\(I(\Sigma_{\mathbf O})\)
&
\(S^0\times\Sigma_{\mathbf O}\), \eqref{eq:line-gluing-disconnected}
\\
\(G_{26''}\)
&
\(\mathbf O\times_{\mathbf O}\mathbf O^2\), Figure~11
&
\(N_{\mathbf T}\)
&
\(\mathcal M(\mathbf T,\varepsilon_{\mathbf T})\)
&
\(\Sigma_{\mathbf T}\), \eqref{eq:line-gluing-connected}
\\
\(G_{31}\)
&
\(\mathbf J\times_{\mathbf J}\mathbf J\), Figure~3
&
\(N_{\mathbf I}\)
&
\(\mathcal M(\mathbf I,\alpha)\)
&
\(\Sigma_{\mathbf I}\), \eqref{eq:line-gluing-connected}
\\
\(G_{31'}\)
&
\(\mathbf J\times_{\mathbf J}\mathbf J^1\), Figure~10
&
\(S^0\times N_{\mathbf I}\)
&
\(I(\Sigma_{\mathbf I})\)
&
\(S^0\times\Sigma_{\mathbf I}\), \eqref{eq:line-gluing-disconnected}
\\
\hline
\end{tabular}
}

\begingroup
\setlength{\rotFPtop}{0pt plus 1fil}
\setlength{\rotFPbot}{0pt plus 1fil}
\begin{sidewaystable}[p]
\centering
\footnotesize
\setlength{\tabcolsep}{1.5pt}
\renewcommand{\arraystretch}{1.05}
\setlength{\abovecaptionskip}{2pt}
\setlength{\belowcaptionskip}{1pt}
\makeatletter
\def\@captionfont{\scriptsize\normalfont}
\def\@captionheadfont{\scriptsize\scshape}
\makeatother
\IndexTwoTableContent
\caption{Attaching regions, belt regions, common boundaries, and gluing maps for index-\(2\) handles. See Figures~\ref{fig:G1-2h}, \ref{fig:G3-2h}\textup{(a)}, and~\ref{fig:G11-2h}.}\label{tab:index-two-models}
\CyclicLineTableContent
\caption{Attaching regions, belt regions, and gluing maps for index-\(1\) handles in \textup{Families} \(1\) and \(1'\). See Figures~\ref{fig:G1-1h-rneq2}--\ref{fig:G1p-1h-odd}.}\label{tab:cyclic-line-regions}
\DihedralLineTableContent
\caption{Attaching regions, belt regions, and gluing maps for index-\(1\) handles in \textup{Families} \(11\) and \(11'\). See Figures~\ref{fig:G1-1h-rneq2}\textup{(a)} and~\ref{fig:G11-1h-att}--\ref{fig:G11p-even-belt}.}\label{tab:dihedral-line-regions}
\PolyhedralTableContent
\caption{Attaching regions, belt regions, and gluing maps for index-\(1\) handles in the polyhedral families. For the diagrams in Appendix~\ref{subsec:illustrated-boundary-models}, see Figures~\ref{fig:G21-1h}--\ref{fig:G31p-1h}.}\label{tab:polyhedral-regions}
\end{sidewaystable}
\endgroup

\FloatBarrier

For the index-\(1\) \(G_1\) models in Table~\ref{tab:cyclic-line-regions}, edges of order \(1\) are omitted from the diagrams.
For \(r=1\), \(S^3\sslash G_1\cong\mathbb{RP}^3\) and \(B_\rho\cong\mathbb{RP}^3\setminus\operatorname{int}D^3\).
For even \(r\geq2\), the underlying space of \(S^3\sslash G_1\) is \(S^3\), and its singular set is a Hopf link with labels \(r\) and \(2\).
For odd \(r\geq3\), the underlying space is \(\mathbb{RP}^3\), and the singular set is a core circle, labeled \(r\), of a genus-one Heegaard splitting~\cite[Table~3, p.~832]{mecchia2015fibered}.

Dunbar represented the complete labeled singular graphs of the ambient quotients in Figures~3, 10, and~11~\cite[pp.~129 and~139]{dunbar1994nonfibering}.
These figures are cited in the second column of Table~\ref{tab:polyhedral-regions}.

\newcommand{\HandleDiagram}[1]{\includegraphics[height=35mm,keepaspectratio]{#1}}

\subsection{The family calculations}\label{subsec:family-boundary-calculations}

\begin{lemma}\label{lem:cyclic-solid-torus-calculation}
Consider the splitting \(\operatorname{span}\{1,\mathbf i\}\oplus\operatorname{span}\{\mathbf j,\mathbf k\}\) for \textup{Families} \(1\) and \(1'\).

For \textup{Family} \(1\), there are orientation-preserving orbifold diffeomorphisms \(A\cong T_{f_1b_1h}\) and \(B\cong T_{f_2b_2h}\).
For \textup{Family} \(1'\), there are orientation-preserving orbifold diffeomorphisms \(A\cong T_{b_1h}\) and \(B\cong T_{b_2h}\).

Use the ordered meridian--longitude bases \((\mu_B,\lambda_B)\) and \((\mu_A,\lambda_A)\).
Under the fixed model identifications, the boundary gluing map in either family is represented by
\[
(\Phi_\rho)_*=
\begin{pmatrix}
-q_G&u_G\\
f&v_G
\end{pmatrix}.
\]
\end{lemma}

\begin{proof}
With the identification in Section~\ref{subsec:corner-smoothing}, \(A\) and \(B\) are the quotients of the two invariant solid tori in the standard genus-one decomposition of \(S^3\).
They correspond to \(T_2\) and \(T_1\), respectively, in the notation of Mecchia-Seppi~\cite[p.~821]{mecchia2015fibered}.

Mecchia-Seppi record the isotropy orders along the cores of \(A\) and \(B\) in \textup{Family} \(1\) as \(f_1b_1h\) and \(f_2b_2h\), respectively~\cite[Table~3, p.~832]{mecchia2015fibered}.
They record the isotropy orders along the cores of \(A\) and \(B\) in \textup{Family} \(1'\) as \(b_1h\) and \(b_2h\), respectively~\cite[Table~2, p.~828]{mecchia2015fibered}.
These identifications give the required oriented models for \(A\) and \(B\).
First take the quotient by the subgroup generated by elements fixing one of the two core circles.
This quotient preserves the two solid tori.
In the resulting sphere, the remaining action is free cyclic of order \(f\), with generator
\[
(z_1,z_2)\longmapsto
\bigl(e^{2\pi i g/f}z_1,e^{2\pi i d/f}z_2\bigr)
\]
in the notation of Mecchia-Seppi~\cite[pp.~827 and~830--832]{mecchia2015fibered}.
Apply their quotient calculation~\cite[Lemma~4.1, pp.~821--822]{mecchia2015fibered} to \(T_2\), interchanging \(d\) and \(g\).
Since the meridian of \(T_1\) is a longitude of \(T_2\), this gives oriented meridian--longitude bases satisfying
\[
(\Phi_\rho)_*(\mu_B)
=
f\lambda_A-d\overline g\,\mu_A.
\]
For \(f=1\), take \(\overline g=0\).
Then \((\Phi_\rho)_*(\mu_B)=\lambda_A\).
Write \(d\overline g=q_G+kf\) for some \(k\in\mathbb Z\).
A \(k\)-fold meridional twist of the model for \(A\) fixes \(\mu_A\) and sends \(\lambda_A\) to \(\lambda_A+k\mu_A\).
Postcompose the diffeomorphism from \(A\) to its model with this twist.
Then
\[
(\Phi_\rho)_*(\mu_B)
=
f\lambda_A-q_G\mu_A.
\]

Write \((\Phi_\rho)_*(\lambda_B)=u\mu_A+v\lambda_A\) for integers \(u\) and \(v\).
Since \(\Phi_\rho\) reverses boundary orientations, the induced matrix has determinant \(-1\).
Hence \(fu+q_Gv=1\).
Therefore \(v\equiv v_G\pmod f\).

A meridional twist of the model for \(B\) changes \(v\) by an integral multiple of \(f\).
Choose the diffeomorphism from \(B\) to its model so that \(v=v_G\).
The equality \(fu+q_Gv=1\) then gives \(u=u_G\).
Therefore
\[
(\Phi_\rho)_*
=
\begin{pmatrix}
-q_G&u_G\\
f&v_G
\end{pmatrix}.
\]
\end{proof}

\begin{lemma}\label{lem:cyclic-pillow-calculation}
Consider the splitting \(\operatorname{span}\{1,\mathbf j\}\oplus\operatorname{span}\{\mathbf i,\mathbf k\}\) in \textup{Family} \(1\) with \(m=n=1\) and \(r=2\), or in \textup{Family} \(1'\) with \(m=n=1\) and \(r'=2\).

In both cases, there are orientation-preserving orbifold diffeomorphisms \(A\cong\mathcal P_1\) and \(B\cong\mathcal P_1\).
Under these identifications, the boundary gluing map satisfies
\[
[\Phi_\rho]
=
\begin{cases}
\tau_b^2\iota, & \textup{in Family }1,\\
\sigma_a\iota, & \textup{in Family }1'.
\end{cases}
\]
\end{lemma}

\begin{proof}
Set \(V^-=\operatorname{span}\{1,\mathbf j\}\) and \(V^+=\operatorname{span}\{\mathbf i,\mathbf k\}\).
Identify these subspaces with \(\mathbb C\) using the ordered bases \((1,\mathbf j)\) and \((\mathbf i,-\mathbf k)\), respectively.
The concatenated basis has the orientation of \((1,\mathbf i,\mathbf j,\mathbf k)\).
Write \(x,y\) for the resulting coordinates.
For \textup{Family} \(1'\) with \(s\equiv3\pmod4\), replace \((x,y)\) by \((ix,iy)\).

The generators in Part~I~\cite[Section~3.3]{choilee2026part1} then give
\[
G_1=\langle-I_4,\eta\rangle,
\qquad
G_{1'}=\langle\eta\rangle,
\]
where
\[
-I_4(x,y)=(-x,-y),
\qquad
\eta(x,y)=(\overline x,\overline y).
\]
Set \(\widetilde A=S^1_x\times D^2_y\) and \(\widetilde B=D^2_x\times S^1_y\).
On their common boundary, write
\[
x=e^{2\pi iu},
\qquad
y=e^{2\pi iv}.
\]
The involution \(\eta\) acts by \((u,v)\mapsto(-u,-v)\).
We use \((u,v)=(0,0)\) as the origin of the boundary coordinates.

\smallskip

\noindent\textup{Family \(1\).}
The action of \(-I_4\) on the common boundary is
\[
(u,v)\longmapsto(u+\tfrac12,v+\tfrac12).
\]
Thus the torus double cover of \(C_\rho\) is \(\mathbb R^2/\Lambda\), where
\[
\Lambda
=
\mathbb Z^2+
\mathbb Z(\tfrac12,\tfrac12).
\]
Regard the vectors of \(\Lambda\) as homology classes of this torus.

The quotient solid tori have meridian--longitude bases
\[
\begin{aligned}
\mu_A&=(0,1),
&
\lambda_A&=(\tfrac12,\tfrac12),\\
\mu_B&=(1,0),
&
\lambda_B&=(\tfrac12,\tfrac12).
\end{aligned}
\]
These bases agree with the boundary orientations of the two regions.

The disk and circle angular coordinates determined by these bases are \((v-u,2u)\) for \(\widetilde A/\langle-I_4\rangle\) and \((u-v,2v)\) for \(\widetilde B/\langle-I_4\rangle\).
Together with the disk radii \(|y|\) and \(|x|\), respectively, they define orientation-preserving diffeomorphisms of the quotient solid tori with \(D^2\times S^1\).
In these coordinates, the involution induced by \(\eta\) is \((z,w)\mapsto(\overline z,\overline w)\).
Hence these diffeomorphisms induce orientation-preserving orbifold diffeomorphisms
\[
A\cong\mathcal P_1,
\qquad
B\cong\mathcal P_1.
\]

On the common boundary,
\[
\mu_B=-\mu_A+2\lambda_A,
\qquad
\lambda_B=\lambda_A.
\]
The boundary gluing map therefore lifts to the torus automorphism represented by
\[
M_1=
\begin{pmatrix}
-1&0\\
2&1
\end{pmatrix}.
\]
The chosen coordinates have the same origin, so the lift has no translation part.
By Definition~\ref{def:pillow-boundary-coordinates},
\[
M_1R
=
\begin{pmatrix}
-1&0\\
2&-1
\end{pmatrix}
=
-\widetilde\tau_b^{\,2}.
\]
Since \(-\Id\) is the deck transformation of the torus double cover, \([\Phi_\rho]=\tau_b^2\iota\).

\smallskip

\noindent\textup{Family \(1'\).}
The torus double cover of \(C_\rho\) is \(\mathbb R^2/\mathbb Z^2\).
The meridian--longitude bases are
\[
\begin{aligned}
\mu_A&=(0,1),
&
\lambda_A&=(1,0),\\
\mu_B&=(1,0),
&
\lambda_B&=(0,1).
\end{aligned}
\]
Interchanging the factors of \(\widetilde A\) and retaining the product coordinates of \(\widetilde B\) identify both solid tori with \(D^2\times S^1\).
These identifications preserve orientation and intertwine \(\eta\) with \(\mathfrak r\).
They therefore induce orientation-preserving orbifold diffeomorphisms
\[
A\cong\mathcal P_1,
\qquad
B\cong\mathcal P_1.
\]

On the common boundary,
\[
\mu_B=\lambda_A,
\qquad
\lambda_B=\mu_A.
\]
Thus the boundary gluing map lifts to the torus automorphism represented by
\[
M_{1'}=
\begin{pmatrix}
0&1\\
1&0
\end{pmatrix}.
\]
The lift has no translation part.
By Definition~\ref{def:pillow-boundary-coordinates},
\[
M_{1'}R
=
\begin{pmatrix}
0&-1\\
1&0
\end{pmatrix}
=
\widetilde\sigma_a.
\]
Therefore \([\Phi_\rho]=\sigma_a\iota\).
\end{proof}

\begin{lemma}\label{lem:pillow-shear}
Let \(G=G_3\) with \(m=1\) and \(n\geq2\), or \(G=G_{34}\) with \(m=1\) and odd \(n\geq3\).
Consider the splitting \(V^-=W_\theta\) and \(V^+=(V^-)^\perp\).
Set
\[
d_G=
\begin{cases}
2n, & G=G_3,\\
n, & G=G_{34}.
\end{cases}
\]
There are orientation-preserving orbifold diffeomorphisms \(A\cong\mathcal P_1\) and \(B\cong\mathcal P_1\).
Under these identifications, the boundary gluing map satisfies
\[
[\Phi_\rho]=\tau_b^{-d_G}\iota.
\]
\end{lemma}

\begin{proof}
By Part~I~\cite[Proposition~4.6]{choilee2026part1}, it suffices to take \(\theta=0\).
Use the ordered bases
\[
\left(\frac{1+\mathbf j}{\sqrt2},\frac{\mathbf i-\mathbf k}{\sqrt2}\right),
\qquad
\left(\frac{1-\mathbf j}{\sqrt2},\frac{-\mathbf i-\mathbf k}{\sqrt2}\right)
\]
for \(V^-\) and \(V^+\), respectively.
Their concatenation has the orientation of \((1,\mathbf i,\mathbf j,\mathbf k)\).
Identify each summand with \(\mathbb C\) using these bases.
Multiply both complex coordinates by \(e^{-\pi i/4}\), and write \(x,y\) for the resulting coordinates.

Set \(\omega=e^{-2\pi i/d_G}\).
The generator formulas in Part~I~\cite[proofs of Propositions~A.9 and~A.16]{choilee2026part1} give \(G=\langle\upsilon,\eta\rangle\), where
\[
\upsilon(x,y)=(\omega x,\omega^{-1}y),
\qquad
\eta(x,y)=(\overline x,\overline y).
\]
Set \(\widetilde A=S^1_x\times D^2_y\) and \(\widetilde B=D^2_x\times S^1_y\).
On their common boundary, write
\[
x=e^{2\pi iu},
\qquad
y=e^{2\pi iv}.
\]
The action of \(\upsilon\) is translation by \((-1/d_G,1/d_G)\).
Thus the quotient by \(\langle\upsilon\rangle\) is the torus \(\mathbb R^2/\Lambda\), where
\[
\Lambda
=
\mathbb Z^2+
\mathbb Z(-1/d_G,1/d_G).
\]
The involution \(\eta\) acts on this torus by \((u,v)\mapsto(-u,-v)\).
Hence \(\mathbb R^2/\Lambda\) is the torus double cover of \(C_\rho\).

Regard vectors of \(\Lambda\) as homology classes of this torus.
The quotients of \(\widetilde A\) and \(\widetilde B\) by \(\langle\upsilon\rangle\) have meridian--longitude bases
\[
\begin{aligned}
\mu_A&=(0,1),
&
\lambda_A&=(1/d_G,-1/d_G),\\
\mu_B&=(1,0),
&
\lambda_B&=(-1/d_G,1/d_G).
\end{aligned}
\]
The corresponding disk and circle coordinates are \((xy,x^{d_G})\) on \(\widetilde A\) and \((xy,y^{d_G})\) on \(\widetilde B\).
These maps are invariant under \(\upsilon\), and their fibers are the \(\langle\upsilon\rangle\)-orbits.
Local choices of a \(d_G\)-th root of the circle coordinate give smooth inverses on the quotients.

In the oriented polar coordinates \((|y|,\arg y,\arg x)\) and \((|x|,\arg x,\arg y)\), respectively, both maps have Jacobian determinant \(d_G>0\).
They therefore induce orientation-preserving diffeomorphisms of the quotient solid tori with \(D^2\times S^1\).
Under these identifications, \(\eta\) induces \((z,w)\mapsto(\overline z,\overline w)\).
Thus there are orientation-preserving orbifold diffeomorphisms
\[
A\cong\mathcal P_1,
\qquad
B\cong\mathcal P_1.
\]

In \(\Lambda\), the meridian--longitude classes satisfy
\[
\mu_B=\mu_A+d_G\lambda_A,
\qquad
\lambda_B=-\lambda_A.
\]
The boundary gluing map therefore lifts to the torus automorphism represented by
\[
M=
\begin{pmatrix}
1&0\\
d_G&-1
\end{pmatrix}.
\]
Both coordinate systems have the same origin, so the lift has no translation part.
By Definition~\ref{def:pillow-boundary-coordinates},
\[
M
=
\widetilde\tau_b^{-d_G}R.
\]
Therefore \([\Phi_\rho]=\tau_b^{-d_G}\iota\).
\end{proof}

We apply Dunbar's construction to \textup{Families} \(11\) and \(11'\).
Mecchia and Seppi obtain their quotient topology by applying Proposition~2.10 of Dunbar's thesis~\cite[p.~812]{mecchia2015fibered}.
We use the boundary markings of Section~\ref{subsec:rational-tangle-models} to express the gluing map.

For \(G\in\{G_{11},G_{11'}\}\), Mecchia and Seppi describe the base of the Seifert fibration of \(S^3\sslash G\) as a disk with two corner reflectors and no cone points~\cite[Tables~2 and~3, pp.~828 and~832]{mecchia2015fibered}.
For the splitting \(\operatorname{span}\{1,\mathbf i\}\oplus\operatorname{span}\{\mathbf j,\mathbf k\}\), the regions \(A\) and \(B\) are the inverse images of neighborhoods of the two corner reflectors.

The exceptional interval fibers have isotropy orders
\[
(c_A,c_B)=
\begin{cases}
(f_1b_1h,f_2b_2h), & G=G_{11},\\
(b_1h,b_2h), & G=G_{11'}.
\end{cases}
\]
The corresponding rational tangle slopes are \(s_A\) and \(s_B\).
Apply the construction preceding Proposition~\ref{lem:ordered-two-corner-gluing} to \(A\) and \(B\).
With \(e=e_G\), Definition~\ref{def:family-11-central-box-labels} gives \(\kappa=\kappa_G\).
Proposition~\ref{lem:ordered-two-corner-gluing} therefore gives
\[
[\Phi_\rho]
=
t(s_A)^{-1}\tau_b^{\kappa_G}\iota t(s_B).
\]

\begin{proof}[\textup{Proof of Theorem~\ref{thm:computed-boundary-regions}}]
Part~I classifies the ordered handle classes of indices \(1\), \(2\), and \(3\) and specifies representatives for these classes~\cite[Theorem~B and Table~1]{choilee2026part1}.
Lemmas~\ref{lem:cyclic-solid-torus-calculation}, \ref{lem:cyclic-pillow-calculation}, and~\ref{lem:pillow-shear} establish the \(G_1\), \(G_{1'}\), \(G_3\), and \(G_{34}\) rows of Table~\ref{tab:index-two-models}.
The calculation for \textup{Families} \(11\) and \(11'\) establishes the remaining two rows.

\smallskip

\noindent\textup{Invariant lines.}
Let \(L\subset\mathbb R^4\) be an invariant line such that \(L\oplus L^\perp\) represents one of the classes classified in Part~I.
Write \(\chi_L\colon G\to\{\pm1\}\) for the character on \(L\).
If \(\chi_L\) is trivial, let \(K\) be the image of \(G\) on \(L^\perp\).
If \(\chi_L\) is nontrivial, let \(K\) be the image of \(\ker\chi_L\) on \(L^\perp\).
In this case, choose \(g_0\in G\setminus\ker\chi_L\) and set \(\nu=\rho(g_0)|_{L^\perp}\).
Lemma~\ref{lem:line-region-models} determines the attaching and belt regions from \(K\) and, when \(\chi_L\) is nontrivial, from \(\nu\).

For \(G_1\) and \(G_{11}\), the character \(\chi_L\) is nontrivial.
The corresponding groups \(K\) are \(\Cyc{r}\) and \(\Dih{2r}\), respectively.
In both families, \(\nu\) is represented by \(\alpha\).

For \(G_{1'}\) with \(r'>2\) and even \(\ell'\), the character \(\chi_L\) is trivial and \(K=\Cyc{r'}\).
For \(r'>2\) and odd \(\ell'\), the character \(\chi_L\) is nontrivial and \((K,\nu)=(\Cyc{r''},\delta_{r''})\).
For \(r'=2\), the lines in \(U_{\mathrm{triv}}\) have trivial character and \(K=\Cyc{2}\).
The lines in \(U_{\mathrm{nontriv}}\) have nontrivial character and satisfy \((K,\nu)=(\Cyc{1},\delta_1)\).

For \(G_{11'}\) with \(r'>2\) and even \(\ell'\), the line \(L_0\) has trivial character and \(K=\Dih{2r'}\).
The line \(L_1\) has nontrivial character and satisfies \((K,\nu)=(\Cyc{r'},\gamma)\).
For \(r'>2\) and odd \(\ell'\), the character \(\chi_L\) is nontrivial and \((K,\nu)=(\Dih{2r''},\delta_{r''})\).

For \(G_{11'}\) with \(r'=2\), a line \(L\) fixed pointwise by \(G_{11'}\) has \(\chi_L=1\) and \(K=\Dih{4}\).
For the orbit of lines with nontrivial \(\chi_L\), use the diagonal representative in Part~I~\cite[proof of Proposition~4.4]{choilee2026part1} and choose \(L=\operatorname{span}\{\mathbf i\}\).
In the ordered basis \((1,\mathbf j,\mathbf k)\) of \(L^\perp\), the image of \(\ker\chi_L\) is generated by \(\operatorname{diag}(1,-1,-1)\).
The element \(\operatorname{diag}(1,-1,1,-1)\) restricts to \(\operatorname{diag}(1,1,-1)\) on \(L^\perp\).
The permutation \(P(a,b,c)=(c,a,b)\) preserves orientation and conjugates these transformations to \(J\) and \(\delta_1\), respectively.
Thus \((K,\nu)=(\Dih{2},\delta_1)\) in these coordinates.

Lemma~\ref{lem:line-region-models} identifies the attaching and belt regions in Tables~\ref{tab:cyclic-line-regions} and~\ref{tab:dihedral-line-regions}.
The same lemma identifies the boundary gluing maps as \eqref{eq:line-gluing-connected} and \eqref{eq:line-gluing-disconnected}.
Tables~\ref{tab:cyclic-line-regions} and~\ref{tab:dihedral-line-regions} record the applicable formula in each case.

\smallskip

\noindent\textup{Polyhedral families.}
Let \(L\oplus W\) be a representative invariant splitting with \(\dim L=1\).
Write \(\chi_L\colon G\to\{\pm1\}\) for the character on \(L\).
For \(G_{21}\), \(G_{26}\), \(G_{26''}\), and \(G_{31}\), the character \(\chi_L\) is nontrivial.
The corresponding groups \(K\) are \(\mathbf T\), \(\mathbf O\), \(\mathbf T\), and \(\mathbf I\), respectively.
The corresponding values of \(\nu\) are \(\alpha\), \(\alpha\), \(\varepsilon_{\mathbf T}\), and \(\alpha\).

For \(G_{21'}\), \(G_{26'}\), and \(G_{31'}\), the character \(\chi_L\) is trivial.
The corresponding groups \(K\) are \(\mathbf T\), \(\mathbf O\), and \(\mathbf I\), respectively.
Lemma~\ref{lem:line-region-models} identifies the attaching and belt regions recorded in Table~\ref{tab:polyhedral-regions}.

For \(G_{21}\), \(G_{26}\), \(G_{26''}\), and \(G_{31}\), the gluing map is \eqref{eq:line-gluing-connected}.
For \(G_{21'}\), \(G_{26'}\), and \(G_{31'}\), it is \eqref{eq:line-gluing-disconnected}.
These maps are recorded in Table~\ref{tab:polyhedral-regions}.

Ordered duality interchanges the attaching and belt regions and replaces \(\Phi_\rho\) by \(\Phi_\rho^{-1}\).
For a representative obtained by factor interchange, the diffeomorphisms induced by \(\mu\) transport \(A_\rho\), \(B_\rho\), and \(\Phi_\rho\) to the corresponding models.
Ordered duality and factor interchange account for the remaining representatives in the statement of the theorem.
\end{proof}

\subsection{Diffeomorphisms relative to the exterior}\label{subsec:relative-handle-attachments}

We compare attachments of a fixed handle along embeddings with the same image.
The diffeomorphisms of the resulting cobordisms are required to fix a common exterior pointwise.
Theorem~\ref{thm:relative-handle-attachments} establishes a criterion in terms of extension over the belt region.

\begin{definition}\label{def:relative-handle-attachments}
Fix an oriented orbifold handle \(H\) as in Definition~\ref{def:quotient-handle-boundary-regions}.
Write \(A\) and \(B\) for its attaching and belt regions.
Write \(\varphi\colon\orbpartial B\to\orbpartial A\) for its boundary identification.
Let \(\mathcal O\) be a closed oriented \(3\)-orbifold.
Let \(\psi_0,\psi_1\colon A\to\mathcal O\) be orientation-reversing full orbifold embeddings with \(\psi_0(A)=\psi_1(A)\).
For \(i=0,1\), set
\[
W_i=(\mathcal O\times[0,1])\cup_{\psi_i}H.
\]

Use the smooth handles and attachments of Section~\ref{subsec:corner-smoothing}.
Choose a common closed subset \(E\subset\mathcal O\times[0,1]\) and compact suborbifolds \(N_i\subset W_i\) with smooth boundary such that
\[
N_i=\overline{W_i\setminus E}.
\]
Choose orientation-preserving orbifold diffeomorphisms \(\eta_i\colon H\to N_i\) satisfying
\begin{equation}\label{eq:common-exterior-pieces}
W_i=E\cup N_i,
\qquad
\eta_i(A)=E\cap N_i,
\qquad
\eta_i(B)=N_i\cap\orbpartial W_i.
\end{equation}
Require
\begin{equation}\label{eq:common-exterior-parametrizations}
\eta_0|_A=\eta_1|_A\circ\psi_1^{-1}\circ\psi_0.
\end{equation}
The suborbifolds \(N_0\) and \(N_1\) agree on an open subset of \(\mathcal O\times[0,1]\) containing \(E\cap N_i\).

The attachments \(W_0\) and \(W_1\) are \emph{equivalent relative to \(E\)} if there is an orientation-preserving orbifold diffeomorphism \(F\colon W_0\to W_1\) satisfying \(F|_E=\Id_E\).
\end{definition}

The suborbifolds \(E\) and \(N_i\) and the diffeomorphisms \(\eta_i\) are obtained using compatible orbifold collars of the common attaching image~\cite[Section~4.4.2]{choi2012geometric}\cite[Theorems~3.5 and~7.2]{kankaanrinta2007equivariant}.
We use the smoothing of Section~\ref{subsec:corner-smoothing} for both attachments.

\begin{lemma}\label{lem:handle-boundary-extension}
Let \(H\) be an oriented orbifold handle.
Every orientation-preserving orbifold diffeomorphism \(f\colon\orbpartial H\to\orbpartial H\) extends to an orientation-preserving orbifold diffeomorphism \(H\to H\).
\end{lemma}

\begin{proof}
By Section~\ref{subsec:corner-smoothing}, identify \(H\) with \(D^4\sslash G\), where \(G<\SO(4)\) is finite.
The main theorem of~\cite{mecchiaseppi2019isometry} implies that \(f\) is isotopic to an isometry of \(S^3\sslash G\).
Hence there is a smooth orbifold isotopy \(f_t\), \(0\leq t\leq1\), from an isometry \(f_0=a\) to \(f_1=f\).
The isometry \(a\) preserves orientation and lifts to a linear isometry \(L\in\SO(4)\) with \(LGL^{-1}=G\).

Lift the isotopy to \(S^3\), starting with \(L|_{S^3}\), and denote its lift by \(\widetilde f_t\).
For each \(g\in G\), the two lifts \(\widetilde f_t\circ g\) and \(\widetilde f_t\) differ by an element of the finite group \(G\).
Since \(G\) is finite, this element is constant in \(t\) and equals \(LgL^{-1}\) at \(t=0\).
Thus
\[
\widetilde f_t(gp)=(LgL^{-1})\widetilde f_t(p).
\]
Choose a smooth function \(\chi\colon[0,1]\to[0,1]\) with \(\chi=0\) on \([0,1/3]\) and \(\chi=1\) on \([2/3,1]\).
Define
\[
\widetilde h(0)=0,
\qquad
\widetilde h(rp)=r\widetilde f_{\chi(r)}(p)
\quad (p\in S^3,\ 0<r\leq1).
\]
For \(\lVert v\rVert\leq1/3\), we have \(\widetilde h(v)=L(v)\).
Its inverse on the punctured ball sends \(rq\) to \(r\widetilde f_{\chi(r)}^{-1}(q)\).
It is therefore a smooth diffeomorphism of the ball, equivariant for the conjugation \(g\mapsto LgL^{-1}\).
Since \(d\widetilde h_0=L\in\SO(4)\), the diffeomorphism \(\widetilde h\) preserves orientation.
The induced orbifold diffeomorphism has boundary restriction \(f\).
\end{proof}

\begin{theorem}\label{thm:relative-handle-attachments}
Assume the hypotheses and use the notation of
Definition~\ref{def:relative-handle-attachments}.
The attachments \(W_0\) and \(W_1\) are equivalent relative to \(E\)
if and only if there is an orientation-preserving orbifold
diffeomorphism \(b\colon B\to B\) satisfying
\begin{equation}\label{eq:relative-belt-boundary}
\psi_1|_{\orbpartial A}\circ\varphi\circ b|_{\orbpartial B}
=\psi_0|_{\orbpartial A}\circ\varphi.
\end{equation}
\end{theorem}

\begin{proof}
Use the suborbifolds \(N_i=\overline{W_i\setminus E}\) and diffeomorphisms \(\eta_i\colon H\to N_i\) from Definition~\ref{def:relative-handle-attachments}.
Suppose that \(F\colon W_0\to W_1\) is an equivalence relative to \(E\).
Since \(N_i=\overline{W_i\setminus E}\), one has \(F(N_0)=N_1\).
By \eqref{eq:common-exterior-pieces} and \eqref{eq:common-exterior-parametrizations}, the map
\[
h=\eta_1^{-1}\circ F|_{N_0}\circ\eta_0
\colon H\longrightarrow H
\]
satisfies \(h|_{A}=\psi_1^{-1}\circ\psi_0\) and \(h(B)=B\).
Its restrictions to the two boundary regions agree on their common boundary.
Therefore
\[
(\psi_1^{-1}\circ\psi_0)|_{\orbpartial A}\circ\varphi
=\varphi\circ h|_{\orbpartial B}.
\]
The map \(b=h|_{B}\) satisfies \eqref{eq:relative-belt-boundary}.

Conversely, suppose that \(b\) satisfies \eqref{eq:relative-belt-boundary}.
Transport \(b\) by \(\eta_0\) and \(\eta_1\) to a diffeomorphism from \(N_0\cap\orbpartial W_0\) to \(N_1\cap\orbpartial W_1\).
Equations~\eqref{eq:common-exterior-parametrizations} and \eqref{eq:relative-belt-boundary} imply that this map is the identity on the common boundary of these regions.
The regions coincide on a collar of their common boundary.
Their common boundary collar has a finite manifold cover inherited from the product handle.
Since the map fixes the boundary, it lifts on a smaller collar.
Apply equivariant collar straightening~\cite[Theorems~9.2 and~9.3]{kankaanrinta2007equivariant} to the lifted collars and pass to the quotient.
The adjusted map is the identity on a smaller collar.
It consequently joins the identity on \(E\cap N_i\) to form a smooth orientation-preserving orbifold diffeomorphism
\[
f\colon\orbpartial N_0\longrightarrow\orbpartial N_1.
\]
The map \(f\) is the identity on an open neighborhood of \(E\cap N_i\) in \(\orbpartial N_i\).

By Lemma~\ref{lem:handle-boundary-extension} and the maps \(\eta_i\), extend \(f\) to an orientation-preserving orbifold diffeomorphism \(h\colon N_0\to N_1\).
Use collars of \(\orbpartial N_i\) whose restrictions to an open subset containing \(E\cap N_i\) coincide.
Apply equivariant collar straightening to \(h\), keeping \(h|_{\orbpartial N_0}=f\) fixed.
The adjusted map equals the identity on an open subset of \(N_0\) containing \(E\cap N_0\).
The maps \(h\) and \(\Id_E\) define a smooth orbifold diffeomorphism \(F\colon W_0\to W_1\).
Its inverse is \(h^{-1}\) on \(N_1\) and \(\Id_E\) on \(E\), so \(F\) is an equivalence relative to \(E\).
\end{proof}

For a compact oriented orbifold \(B\), let \(\operatorname{Diff}_{\mathrm{orb}}^+(B)\) be its group of orientation-preserving orbifold diffeomorphisms.
The group \(\pi_0\operatorname{Diff}_{\mathrm{orb}}^+(B)\) consists of their smooth orbifold isotopy classes.

Theorem~\ref{thm:relative-handle-attachments} has the following equivalent formulation in terms of mapping classes.
The attachments are equivalent relative to \(E\) if and only if
\[
\left[
\varphi^{-1}\circ(\psi_1^{-1}\circ\psi_0)|_{\orbpartial A}\circ\varphi
\right]
\in
\operatorname{im}\!\left(
\pi_0\operatorname{Diff}_{\mathrm{orb}}^+(B)
\longrightarrow
\MCG^+(\orbpartial B)
\right).
\]
Indeed, the collar construction in the proof of Lemma~\ref{lem:gluing-isotopy} extends an isotopy of the boundary to a diffeomorphism of \(B\).
It therefore adjusts a representative of this image to satisfy \eqref{eq:relative-belt-boundary} exactly.
In particular, an isotopy between \(\psi_0\) and \(\psi_1\) through embeddings with the fixed image gives equivalent attachments relative to \(E\).

The following proposition determines which boundary diffeomorphisms extend over \(T_q\) and \(\mathcal P_q\).
Appendix~\ref{app:spherical-boundary-extension} computes the images of the restriction maps for the remaining models.

\begin{proposition}\label{prop:solid-boundary-extension}
Let \(q\geq1\), and let \(f\) be an orientation-preserving orbifold diffeomorphism of \(\orbpartial T_q\) or \(\orbpartial\mathcal P_q\).
In the torus case, let \(M\in\mathrm{SL}(2,\mathbb Z)\) be its matrix in the meridian--longitude basis.
In the pillow case, let \(M\) be the homology matrix of a lift to the torus double cover, determined up to sign.
Then \(f\) extends to an orientation-preserving orbifold diffeomorphism of the solid model if and only if it preserves the meridian slope, equivalently,
\[
M=\pm\begin{pmatrix}1&n\\0&1\end{pmatrix}
\quad\text{for some }n\in\mathbb Z.
\]
\end{proposition}

\begin{proof}
The orbifold fundamental groups of the torus and pillow boundaries are \(\mathbb Z^2\) and \(\mathbb Z^2\rtimes\Cyc{2}\), respectively, where the involution acts by inversion.
The corresponding solid models have groups
\[
\pi_1^{\mathrm{orb}}(T_q)=\Cyc{q}\times\mathbb Z,
\qquad
\pi_1^{\mathrm{orb}}(\mathcal P_q)
=(\Cyc{q}\times\mathbb Z)\rtimes\Cyc{2}.
\]
In the fixed coordinates, each boundary inclusion is surjective and has kernel generated by the meridional translation \((q,0)\).
An extending diffeomorphism preserves this kernel, so \(M(q,0)^T=\pm(q,0)^T\).
Since \(\det M=1\), we have \(M=\pm\begin{pmatrix}1&n\\0&1\end{pmatrix}\) for some \(n\in\mathbb Z\).

To prove sufficiency, use \(D^2\times\mathbb R\) with the transformations
\[
s(z,t)=(e^{2\pi i/q}z,t),\qquad
\ell(z,t)=(z,t+1),\qquad
r(z,t)=(\overline z,-t).
\]
The quotients by \(\langle s,\ell\rangle\) and \(\langle s,\ell,r\rangle\) are the cyclic solid torus and \(\mathcal P_q\), respectively.
For every \(n\in\mathbb Z\), the orientation-preserving diffeomorphism
\[
F_n(z,t)=(e^{2\pi int/q}z,t)
\]
commutes with \(s,r\) and satisfies \(F_n\ell F_n^{-1}=s^n\ell\).
It therefore induces a diffeomorphism of each solid model.
In the boundary coordinates \((q\arg z/(2\pi),t)\), it induces the shear \(\begin{pmatrix}1&n\\0&1\end{pmatrix}\).
The map \(r\) induces a diffeomorphism of the cyclic solid torus with boundary matrix \(-I\).

For the pillow, the maps
\[
(z,t)\longmapsto(e^{\pi i/q}z,t),\qquad
(z,t)\longmapsto(z,t+\tfrac12)
\]
also normalize \(\langle s,\ell,r\rangle\).
They induce the two half-translations of the boundary torus.
By the affine description following Definition~\ref{def:pillow-boundary-coordinates}, these maps realize every boundary mapping class satisfying the matrix condition.
A collar isotopy adjusts an extending representative to have boundary restriction exactly \(f\).
\end{proof}

\begin{corollary}\label{cor:index-two-framing-relative-exterior}
Assume the hypotheses of Definition~\ref{def:relative-handle-attachments}, and suppose that \(H\) has index \(2\).
Set \(h=(\psi_1^{-1}\circ\psi_0)|_{\orbpartial A}\).
In the boundary coordinates of the solid torus or solid pillow model for \(A\), the matrix of \(h\), determined up to sign in the pillow case, has the form
\[
\pm\begin{pmatrix}1&n\\0&1\end{pmatrix}
\]
for a unique integer \(n\).
The attachments are equivalent relative to \(E\) if and only if \(n=0\).
\end{corollary}

\begin{proof}
The diffeomorphism \(\psi_1^{-1}\circ\psi_0\) extends \(h\) over \(A\).
By Proposition~\ref{prop:solid-boundary-extension}, the matrix of \(h\) has the form
\[
\pm\begin{pmatrix}1&n\\0&1\end{pmatrix}
\]
in the fixed boundary coordinates.

Represent the meridian slope of \(B\), transported to \(\orbpartial A\) by \(\varphi\), by a primitive vector \((x,y)\).
Then \(y\neq0\).
Indeed, if \(y=0\), the meridian slopes of \(A\) and \(B\) agree.
By the kernel calculation in Proposition~\ref{prop:solid-boundary-extension} and the orbifold van Kampen theorem, \(\pi_1^{\mathrm{orb}}(A\cup_\varphi B)\) has an infinite cyclic quotient in the torus case and an infinite dihedral quotient in the pillow case.
This contradicts \(A\cup_\varphi B\cong S^3\sslash G\), whose orbifold fundamental group is the finite group \(G\).

By Theorem~\ref{thm:relative-handle-attachments} and Proposition~\ref{prop:solid-boundary-extension}, the attachments are equivalent relative to \(E\) if and only if \(h\) preserves the meridian slope of \(B\) transported to \(\orbpartial A\).
Up to sign, the matrix of \(h\) sends \((x,y)\) to \((x+ny,y)\).
Since \(y\neq0\), this slope is preserved if and only if \(n=0\).
\end{proof}

For an ordinary handle, write \(A=S^1\times D^2\).
Fix an attaching embedding \(\psi_*\).
Set \(d_n(z,w)=(z,z^n w)\) and \(\psi_n=\psi_*\circ d_n\).
The integer in Corollary~\ref{cor:index-two-framing-relative-exterior} for \(\psi_n^{-1}\circ\psi_m\) is \(m-n\).
Thus the attachments determined by \(\psi_m\) and \(\psi_n\) are equivalent relative to a compatible common exterior if and only if \(m=n\).

The following corollary allows the attaching image to vary during the isotopy.

\begin{corollary}\label{cor:attaching-isotopy-cobordisms}
Let \(\mathcal X\) be a compact oriented \(4\)-orbifold with boundary.
Let \(H\) be an oriented orbifold handle with attaching region \(A\).
Let \(\psi_0,\psi_1\colon A\to\orbpartial\mathcal X\) be smoothly isotopic orientation-reversing full orbifold embeddings.
Then the orbifolds
\[
\mathcal X\cup_{\psi_0}H
\qquad\text{and}\qquad
\mathcal X\cup_{\psi_1}H
\]
are orientation-preservingly orbifold diffeomorphic.
\end{corollary}

\begin{proof}
Let \(\psi_t\colon A\to\orbpartial\mathcal X\), \(0\leq t\leq1\), be a smooth isotopy of orientation-reversing full orbifold embeddings.
We use the vector-field proof of isotopy extension~\cite[Theorem~2.4.2]{wall2016differential}.
In compatible orbifold charts, fullness implies that the derivatives of the local lifts of \(\psi_t\) are equivariant under the target local groups.
Extend these derivatives to local vector fields and average over the finite local groups.
An orbifold partition of unity~\cite[Appendix, Lemma~4.2.1]{chenruan2002orbifoldgw} gives a smooth time-dependent orbifold vector field on \(\orbpartial\mathcal X\) that agrees with the derivatives of \(\psi_t\) along its image.
Since \(\orbpartial\mathcal X\) is compact, this vector field has a flow \(g_t\) for \(0\leq t\leq1\).
Uniqueness of integral curves gives
\[
g_0=\Id_{\orbpartial\mathcal X},
\qquad
g_t\circ\psi_0=\psi_t.
\]

Choose an orbifold collar of \(\orbpartial\mathcal X\) in \(\mathcal X\).
Let \(s\) denote its collar coordinate.
Choose a smooth function \(\chi\) that equals \(1\) for \(s\) sufficiently small and vanishes outside a smaller collar.
The map
\[
(x,s)\longmapsto(g_{\chi(s)}(x),s)
\]
extends by the identity to an orientation-preserving orbifold diffeomorphism \(G\colon\mathcal X\to\mathcal X\).
Its boundary restriction is \(g_1\).
Since \(g_1\circ\psi_0=\psi_1\), the maps \(G\) and \(\Id_H\) induce an orientation-preserving orbifold diffeomorphism
\[
\mathcal X\cup_{\psi_0}H
\longrightarrow
\mathcal X\cup_{\psi_1}H.
\]
The equivariant collar comparison of~\cite[Theorems~9.2 and~9.3]{kankaanrinta2007equivariant} identifies the resulting smoothings.
\end{proof}

\section{Surgery cobordisms}\label{sec:surgery-cobordisms}

In this section, we construct a cobordism from a closed connected oriented \(3\)-orbifold to \(S^3\).
We determine the minimum total number of handles of indices \(1\) and \(3\).

Section~\ref{subsec:local-computations-for-cobordisms} analyzes handle attachments that remove polyhedral vertices, dihedral vertices, and singular link components.
Section~\ref{subsec:cobordism-and-reverse} constructs the cobordism to \(S^3\) and describes the handles in its oriented reverse.
Section~\ref{subsec:minimum-number-indices-one-three} characterizes when index-\(2\) handles suffice and determines the minimum number of index-\(1\) and index-\(3\) handles.

\subsection{Handle attachments and the singular set}\label{subsec:local-computations-for-cobordisms}

We first compare the index-\(1\) and index-\(2\) attachments that remove a pair of dihedral vertices.
The index-\(1\) attachment corresponds to Druschel's pairwise removal~\cite[proof of Theorem~3.2, pp.~51--53]{druschel2000cobordism}.
We use the index-\(2\) attachment in the cobordism construction.
It preserves the diffeomorphism type of the underlying \(3\)-manifold.

We retain the notation \(\mathbf T\), \(\mathbf O\), \(\mathbf I\), and \(v_P(\mathcal O)\) from the Introduction and the notation for orbifold handles in Definition~\ref{def:quotient-handle-boundary-regions}.

For each integer \(q\geq2\), let \(\rho_{11',q}\) and \(\rho_{1',q}\) denote the representations in \textup{Families} \(11'\) and \(1'\), respectively, with \(m=n=s=1\) and \(r'=q\).
For these parameters, \(G_{1'}\) belongs to \textup{Family} \(1'_{\mathrm{(iii)}}\), with \((1-s)/r'=0\).

For \(\rho_{11',q}\), the splittings
\[
\operatorname{span}\{1\}\oplus\operatorname{span}\{\mathbf i,\mathbf j,\mathbf k\}
\qquad\text{and}\qquad
\mathbb C\oplus\mathbb C\mathbf j
\]
determine handles of indices \(1\) and \(2\), respectively.
Let \((A_i,B_i)\) denote the attaching and belt regions of the index-\(i\) handle.
By Tables~\ref{tab:dihedral-line-regions} and~\ref{tab:index-two-models},
\[
\begin{aligned}
(A_1,B_1)
&\cong
\left(
S^0\times(D^3\sslash\Dih{2q}),
I\bigl(S^2(2,2,q)\bigr)
\right),\\
(A_2,B_2)
&\cong
(\mathcal P_q,\mathcal P_1).
\end{aligned}
\]

\begin{proposition}\label{prop:two-G11p-handles}
Let \(q\geq2\).
For \(i=1,2\), let \(H_i\) be the index-\(i\) handle associated with \(\rho_{11',q}\) and the invariant splitting specified above.
Let \(A_i\) be its attaching region.
Let \(\mathcal O\) be a connected oriented \(3\)-orbifold.
Let \(f_i\colon A_i\to\operatorname{int}\mathcal O\) be orientation-reversing orbifold embeddings.
Write \(\mathcal O'_i\) for the outgoing orbifold obtained by attaching \(H_i\) along \(f_i\).
Then
\[
|\mathcal O'_1|\cong|\mathcal O|\mathbin{\#}(S^1\times S^2),
\qquad
|\mathcal O'_2|\cong|\mathcal O|.
\]
\end{proposition}

\begin{proof}
The attaching and belt regions of \(H_1\) satisfy
\[
|A_1|\cong D^3\sqcup D^3,
\qquad
|B_1|\cong S^2\times I.
\]
The first attachment removes the interiors of two disjoint \(3\)-balls and joins their boundary spheres by \(S^2\times I\).
Since \(|\mathcal O|\) is connected,
\[
|\mathcal O'_1|
\cong
|\mathcal O|\mathbin{\#}(S^1\times S^2).
\]

Similarly,
\[
|A_2|\cong D^3,
\qquad
|B_2|\cong D^3.
\]
Every homeomorphism of \(S^2\) extends over \(D^3\).
Thus the second attachment preserves the homeomorphism type of the underlying \(3\)-manifold.
Uniqueness of the smooth structure on a \(3\)-manifold implies that it also preserves the diffeomorphism type~\cite[Section~1]{pardon2021smoothing}.
Therefore
\[
|\mathcal O'_2|\cong|\mathcal O|.
\]
\end{proof}

\begin{proposition}\label{prop:change-of-polyhedral-vertices}
Let \(\mathcal O'\) be obtained from a closed oriented \(3\)-orbifold \(\mathcal O\) by a single orbifold handle attachment.
An index-\(2\) attachment preserves each count \(v_P\), where \(P\in\{\mathbf T,\mathbf O,\mathbf I\}\).
For an index-\(1\) handle associated with a representative group \(G\), the numbers of polyhedral vertices removed are
\[
\begin{array}{c|ccc}
G&\mathbf T&\mathbf O&\mathbf I\\\hline
G_{21},\ G_{26''}&1&0&0\\
G_{21'}&2&0&0\\
G_{26}&0&1&0\\
G_{26'}&0&2&0\\
G_{31}&0&0&1\\
G_{31'}&0&0&2
\end{array}.
\]
These attachments introduce no polyhedral vertices.
Every other index-\(1\) attachment preserves all three counts.
An index-\(3\) handle obtained by ordered duality changes each count by the negative of the change for the corresponding index-\(1\) handle.
\end{proposition}

\begin{proof}
The common boundary \(C_\rho\) has only cyclic local groups and therefore contains no polyhedral vertices.
Hence, for \(P\in\{\mathbf T,\mathbf O,\mathbf I\}\),
\[
v_P(\mathcal O')-v_P(\mathcal O)
=
v_P(B_\rho)-v_P(A_\rho).
\]
Let \(p=(x,y)\) be a point of \(S(V^-)\times D(V^+)\) or \(D(V^-)\times S(V^+)\), where \(x\in V^-\) and \(y\in V^+\).
Let \(G_p<\rho(G)\) be the stabilizer of \(p\).
For an index-\(2\) handle, suppose first that \(x\neq0\).
Restriction to \(V^+\) embeds \(G_p\) in \(\Orth(2)\).
Indeed, an element in the kernel fixes \(V^+\) pointwise and fixes \(x\).
Its determinant is one, so it also fixes the remaining line in \(V^-\).
The same argument with the factors interchanged applies when \(y\neq0\).
Every finite subgroup of \(\Orth(2)\) is cyclic or dihedral.
Thus neither boundary region contains a polyhedral vertex, and each \(v_P\) is preserved.

For an index-\(1\) handle, write \(\mathbb R^4=L\oplus W\) with \(\dim L=1\).
For \(p=(x,y)\in D(L)\times S(W)\), we have \(y\neq0\).
Restriction to \(y^\perp\cap W\) embeds \(G_p\) in \(\Orth(2)\) by the same determinant argument.
Hence the belt region contains no polyhedral vertices.
The same argument applies to points \(p=(x,y)\in S(L)\times D(W)\) with \(y\neq0\).
Thus every polyhedral vertex in the attaching region comes from \(S(L)\times\{0\}\).
This set has two orbits with the same local group when \(G\) fixes \(L\) pointwise.
It has one orbit otherwise.
Tables~\ref{tab:cyclic-line-regions}--\ref{tab:polyhedral-regions} identify the local groups of these orbits.
Ordered duality interchanges \(A_\rho\) and \(B_\rho\), so it reverses the sign of \(v_P(B_\rho)-v_P(A_\rho)\).
\end{proof}

\begin{lemma}\label{lem:remove-polyhedral-vertices}
Let \(\mathcal O\) be a closed oriented \(3\)-orbifold.
There is a finite sequence of index-\(1\) handle attachments from \(\mathcal O\) to a closed oriented \(3\)-orbifold with no polyhedral vertices.
\end{lemma}

\begin{proof}
If \(\mathcal O\) has no polyhedral vertices, take the empty sequence.
Otherwise, let \(x\) be a polyhedral vertex with local group \(P\in\{\mathbf T,\mathbf O,\mathbf I\}\).
Choose a neighborhood \(U\) of \(x\) orbifold diffeomorphic to \(D^3\sslash P\).

For \(P=\mathbf T,\mathbf O,\mathbf I\), choose the index-\(1\) handle associated with \(G_{26''}\), \(G_{26}\), or \(G_{31}\), respectively.
By Table~\ref{tab:polyhedral-regions}, its attaching region is orbifold diffeomorphic to \(D^3\sslash P\).
The map \(-\Id_{D^3}\) commutes with \(P\) and induces an orientation-reversing orbifold self-diffeomorphism of \(D^3\sslash P\).
Hence an orientation-reversing attaching embedding into \(U\) exists.

By Table~\ref{tab:polyhedral-regions}, the belt region contains no polyhedral vertices.
The common boundary also contains no polyhedral vertices.
Thus the attachment removes \(x\) without introducing any polyhedral vertices.

Since \(\mathcal O\) has finitely many polyhedral vertices, repeating this construction gives a finite sequence of index-\(1\) handle attachments whose outgoing orbifold has no polyhedral vertices.
\end{proof}

\begin{lemma}\label{lem:remove-dihedral-vertices}
Let \(\mathcal O\) be a closed oriented \(3\)-orbifold whose trivalent vertices are all dihedral.
There is a finite sequence of index-\(2\) handle attachments from \(\mathcal O\) to a closed oriented \(3\)-orbifold whose singular set is a labeled link.
Every handle in the sequence is associated with \(\rho_{11',q}\) for some integer \(q\geq2\).
\end{lemma}

\begin{proof}
If \(\mathcal O\) has no trivalent vertices, take the empty sequence.
Suppose that \(\mathcal O\) has a trivalent vertex.

Suppose first that \(\mathcal O\) has a vertex of type \((2,2,q)\) with \(q>2\).
Its unique incident edge of order \(q\) cannot be a loop, since a loop would contribute two incidences of order \(q\) at the vertex.
The other endpoint is therefore a distinct vertex of type \((2,2,q)\).

Otherwise, every trivalent vertex has type \((2,2,2)\).
Since a loop contributes two incidences at a vertex, some edge must have distinct endpoints.
Thus, in either case, there is an edge \(e\) of order \(q\) with distinct endpoints of type \((2,2,q)\) for some integer \(q\geq2\).

Let \(U\) be a sufficiently small regular neighborhood of \(e\).
Its underlying space is a \(3\)-ball.
Its singular graph consists of \(e\), its two endpoints, and four edges of order \(2\) joining the endpoints to \(\partial|U|\).
This labeled pair is diffeomorphic to the pair defining \(\mathcal P_q\).
Proposition~\ref{prop:3orb-diffeo-criterion} therefore identifies \(U\) with \(\mathcal P_q\) up to orbifold diffeomorphism.
The map
\[
(z,v)\longmapsto(\overline z,v)
\]
on \(D^2\times S^1\) normalizes \(\Gamma_q\) and induces an orientation-reversing orbifold self-diffeomorphism of \(\mathcal P_q\).
Hence an orientation-reversing attaching embedding into \(U\) exists.

Attaching the index-\(2\) handle associated with \(\rho_{11',q}\) replaces \(\mathcal P_q\) by \(\mathcal P_1\).
Since \(\mathcal P_q\) contains exactly two trivalent vertices and \(\mathcal P_1\) contains none, the attachment removes two trivalent vertices without introducing any new ones.

Since \(\mathcal O\) has finitely many trivalent vertices, repeating this construction removes all of them after finitely many attachments.
The resulting orbifold is closed.
Its singular graph has no vertices.
Hence its singular set is a labeled link.
\end{proof}

\begin{lemma}\label{lem:remove-cyclic-components}
Let \(\mathcal O\) be a closed oriented \(3\)-orbifold whose singular set is a labeled link.
There is a finite sequence of index-\(2\) handle attachments from \(\mathcal O\) to a closed oriented \(3\)-manifold.
Every handle in the sequence is associated with \(\rho_{1',q}\) for some integer \(q\geq2\).
\end{lemma}

\begin{proof}
If \(\Sigma(\mathcal O)=\varnothing\), take the empty sequence.
Otherwise, let \(K\) be a singular component of order \(q\geq2\).
Choose an orientation of \(K\).
The orientations of \(\mathcal O\) and \(K\) orient the normal orbifold disk bundle of \(K\).
Its structure group is \(\SO(2)/\Cyc{q}\cong\SO(2)\).
Since \(K\cong S^1\) and this group is connected, the bundle is trivial.
Thus \(K\) has a tubular neighborhood \(N(K)\) orbifold diffeomorphic to \(T_q\).

By Table~\ref{tab:index-two-models}, the index-\(2\) handle associated with \(\rho_{1',q}\) has attaching and belt regions orbifold diffeomorphic to \(T_q\) and \(T_1\), respectively.
Complex conjugation on the \(S^1\)-factor induces an orientation-reversing orbifold self-diffeomorphism of \(T_q\).
Hence an orientation-reversing attaching embedding into \(N(K)\) exists.

Attaching this handle replaces \(N(K)\) by \(T_1\), which is a smooth solid torus.
This attachment removes \(K\) without introducing any new singular points.
Since \(\Sigma(\mathcal O)\) has finitely many components, repeating the construction removes the entire singular set after finitely many attachments.
The resulting orbifold is a closed oriented \(3\)-manifold.
\end{proof}

\subsection{The cobordism and its reverse}\label{subsec:cobordism-and-reverse}

\begin{theorem}\label{thm:cobordism-to-sphere}
Let \(\mathcal O\) be a closed connected oriented \(3\)-orbifold.
There is a compact connected oriented \(4\)-orbifold \(W\) with
\[
\orbpartial W=(-\mathcal O)\sqcup S^3.
\]
The orbifold \(W\) admits a handle decomposition relative to \(\mathcal O\) with the following four stages.
\begin{enumerate}[label=\textup{(\roman*)}]
\item Attach index-\(1\) handles associated with \(G_{26''}\), \(G_{26}\), and \(G_{31}\).
The outgoing orbifold has no polyhedral vertices.
\item Attach index-\(2\) handles associated with \(\rho_{11',q}\), for \(q\geq2\).
The outgoing orbifold has no trivalent vertices.
\item Attach index-\(2\) handles associated with \(\rho_{1',q}\), for \(q\geq2\).
The outgoing orbifold is a smooth \(3\)-manifold.
\item Attach smooth index-\(2\) handles to obtain \(S^3\).
\end{enumerate}
\end{theorem}

\begin{proof}
By Lemma~\ref{lem:remove-polyhedral-vertices}, the first stage removes all polyhedral vertices.
By Lemma~\ref{lem:remove-dihedral-vertices}, the second stage removes all remaining trivalent vertices.
By Lemma~\ref{lem:remove-cyclic-components}, the third stage removes the remaining singular set.

In the first two stages, each attaching region has underlying space \(D^3\).
Each belt region is connected.
Removing an embedded \(3\)-ball from a connected closed \(3\)-manifold leaves a connected complement.
Hence the outgoing boundaries remain connected.

In the third stage, each attaching region is a tubular neighborhood of a singular link component.
Its complement in the underlying \(3\)-manifold is connected.
The belt region is a smooth solid torus.
Thus the outgoing boundaries remain connected.
After the third stage, the outgoing boundary is a closed connected oriented \(3\)-manifold \(M\).

By the Lickorish-Wallace theorem, \(M\) is obtained from \(S^3\) by surgery on a framed link~\cite{lickorish1962representation,wallace1960modifications}.
Reversing this surgery gives a finite sequence of smooth index-\(2\) handle attachments from \(M\) to \(S^3\).

Concatenating the four stages gives a compact connected \(4\)-orbifold \(W\).
The attaching maps reverse orientation, so the handles can be oriented consistently with the cobordism.
Thus \(W\) is oriented and satisfies
\[
\orbpartial W=(-\mathcal O)\sqcup S^3.
\]
\end{proof}

\begin{lemma}\label{lem:oriented-reverse-attachment}
Let \(\rho\colon G\to\SO(4)\) be a faithful representation of a finite group.
Let \(\mathbb R^4=V^-\oplus V^+\) be a \(\rho(G)\)-invariant orthogonal splitting with \(1\leq\dim V^-\leq3\).
Let \(\mu\colon\mathbb H\to\mathbb H\) be quaternionic conjugation, and set \(\rho^\mu(g)=\mu\circ\rho(g)\circ\mu^{-1}\).
Suppose that attaching \(H_\rho(V^-,V^+)\) to an oriented \(3\)-orbifold \(\mathcal O\) produces \(\mathcal O'\).
Then the oriented reverse is obtained by attaching
\[
H_{\rho^\mu}(\mu V^+,\mu V^-)
\]
to \(\mathcal O'\).
This handle has index \(4-\dim V^-\).
\end{lemma}

\begin{proof}
The ordered dual interchanges \(V^-\) and \(V^+\), and hence the attaching and belt regions~\cite[Definition~1.2]{choilee2026part1}.
Quaternionic conjugation reverses the orientation of \(\mathbb H\) and conjugates \(\rho\) to \(\rho^\mu\).
It therefore induces an orientation-reversing orbifold diffeomorphism
\[
H_\rho(V^-,V^+)\longrightarrow H_{\rho^\mu}(\mu V^+,\mu V^-)
\]
that carries the original belt region to the new attaching region and the original attaching region to the new belt region.
The handle \(H_{\rho^\mu}(\mu V^+,\mu V^-)\) realizes the oriented reverse of the attachment.
Its index is \(\dim V^+=4-\dim V^-\).
\end{proof}

By Theorem~\ref{thm:cobordism-to-sphere} and Lemma~\ref{lem:oriented-reverse-attachment}, every closed connected oriented \(3\)-orbifold \(\mathcal O\) can be obtained from \(S^3\) by a finite sequence of orbifold handle attachments.
The sequence begins with smooth index-\(2\) handles.
These are followed by index-\(2\) handles obtained by reversing the attachments associated with \(\rho_{1',q}\) and then \(\rho_{11',q}\), for \(q\geq2\).
The sequence ends with index-\(3\) handles obtained by reversing the index-\(1\) attachments associated with \(G_{26''}\), \(G_{26}\), and \(G_{31}\).

\subsection{The minimum number of handles of indices \texorpdfstring{\(1\) and \(3\)}{1 and 3}}\label{subsec:minimum-number-indices-one-three}

Let \(\mathcal O\) be a closed connected oriented \(3\)-orbifold.
Consider all finite sequences of orbifold handle attachments of indices \(1\), \(2\), and \(3\) from \(\mathcal O\) to \(S^3\).
The number \(c_{\mathrm{poly}}(\mathcal O)\) is the minimum total number of index-\(1\) and index-\(3\) handles among these sequences.

\begin{theorem}\label{thm:index-two-characterization}
Let \(\mathcal O\) be a closed connected oriented \(3\)-orbifold.
There is a finite sequence of index-\(2\) orbifold handle attachments from \(\mathcal O\) to \(S^3\) if and only if \(\mathcal O\) has no polyhedral vertices.
\end{theorem}

\begin{proof}
If \(\mathcal O\) has no polyhedral vertices, the first stage of Theorem~\ref{thm:cobordism-to-sphere} is empty.
Every remaining handle has index \(2\).

Conversely, by Proposition~\ref{prop:change-of-polyhedral-vertices}, index-\(2\) attachments preserve each \(v_P\), where \(P\in\{\mathbf T,\mathbf O,\mathbf I\}\).
A sequence from \(\mathcal O\) to \(S^3\) therefore gives \(v_P(\mathcal O)=v_P(S^3)=0\) for every such \(P\).
\end{proof}

\begin{theorem}\label{thm:minimum-polyhedral-handles}
For every closed connected oriented \(3\)-orbifold \(\mathcal O\),
\[
c_{\mathrm{poly}}(\mathcal O)
=
\sum_{P\in\{\mathbf T,\mathbf O,\mathbf I\}}
\left\lceil
\frac{v_P(\mathcal O)}{2}
\right\rceil.
\]
Moreover, this minimum is attained by a sequence containing no index-\(3\) handles.
\end{theorem}

\begin{proof}
Consider any finite sequence of handle attachments of indices \(1\), \(2\), and \(3\) from \(\mathcal O\) to \(S^3\).
For each \(P\in\{\mathbf T,\mathbf O,\mathbf I\}\), let \(n_P\) be the number of attachments in this sequence that change \(v_P\).
By Proposition~\ref{prop:change-of-polyhedral-vertices}, each counted handle has index \(1\) or \(3\).
Its change in \(v_P\) has absolute value at most \(2\).
Since \(v_P(S^3)=0\), the total change is \(-v_P(\mathcal O)\).
The triangle inequality gives \(v_P(\mathcal O)\leq2n_P\).
Each handle changes at most one of the three counts, so the total number of handles of indices \(1\) and \(3\) is at least
\[
\sum_{P\in\{\mathbf T,\mathbf O,\mathbf I\}}n_P
\geq
\sum_{P\in\{\mathbf T,\mathbf O,\mathbf I\}}
\left\lceil\frac{v_P(\mathcal O)}{2}\right\rceil.
\]
Taking the minimum gives the required lower bound for \(c_{\mathrm{poly}}(\mathcal O)\).

For each \(P\in\{\mathbf T,\mathbf O,\mathbf I\}\), partition the set of vertices with local group \(P\) into pairs and at most one singleton.
Choose pairwise disjoint small neighborhoods of these vertices, each orbifold diffeomorphic to \(N_P=D^3\sslash P\).

For \(P=\mathbf T,\mathbf O,\mathbf I\), attach an index-\(1\) handle associated with \(G_{21'}\), \(G_{26'}\), or \(G_{31'}\), respectively, to each pair.
If one \(\mathbf T\)-vertex remains, attach the handle associated with \(G_{26''}\).
If one \(\mathbf O\)-vertex remains, attach the handle associated with \(G_{26}\).
If one \(\mathbf I\)-vertex remains, attach the handle associated with \(G_{31}\).

The map \(-\Id_{D^3}\) commutes with \(P\) and induces an orbifold diffeomorphism of \(N_P\) that reverses orientation.
Choose the attaching maps to reverse orientation.
By Table~\ref{tab:polyhedral-regions}, these attachments remove the chosen polyhedral vertices without introducing any new polyhedral vertices.
The number of index-\(1\) handles used is
\[
\sum_{P\in\{\mathbf T,\mathbf O,\mathbf I\}}
\left\lceil
\frac{v_P(\mathcal O)}{2}
\right\rceil.
\]

For each attachment, the underlying attaching region consists of one or two disjoint \(3\)-balls.
The belt region is connected.
Removing the interiors of these balls from a connected closed \(3\)-manifold leaves a connected complement.
Hence the outgoing orbifold remains closed and connected.

By Theorem~\ref{thm:index-two-characterization}, there is a finite sequence of index-\(2\) handle attachments from the outgoing orbifold to \(S^3\).
Appending this sequence gives a sequence from \(\mathcal O\) to \(S^3\) that attains the lower bound.
It contains no index-\(3\) handles.
\end{proof}

\begin{corollary}\label{cor:minimum-polyhedral-fillings}
Let \(\mathcal O\) be a closed connected oriented \(3\)-orbifold.
Consider compact connected oriented \(4\)-orbifolds \(W\) with \(\orbpartial W=\mathcal O\) that admit a handle decomposition consisting of one smooth \(4\)-ball and handles of indices \(1\), \(2\), and \(3\).
Among these decompositions, the minimum total number of index-\(1\) and index-\(3\) handles is \(c_{\mathrm{poly}}(\mathcal O)\).
This minimum is attained by a decomposition with no index-\(1\) handles, in which every index-\(3\) handle is associated with a faithful representation of a nontrivial finite group.

In particular, \(\mathcal O\) bounds a compact connected oriented \(4\)-orbifold obtained from a smooth \(4\)-ball by attaching only index-\(2\) handles if and only if \(\mathcal O\) has no polyhedral vertices.
\end{corollary}

\begin{proof}
Consider the sequence from \(\mathcal O\) to \(S^3\) constructed in the proof of Theorem~\ref{thm:minimum-polyhedral-handles}.
It contains exactly \(c_{\mathrm{poly}}(\mathcal O)\) handles of index \(1\) and no handles of index \(3\).

By Lemma~\ref{lem:oriented-reverse-attachment}, reversing this sequence gives a sequence from \(S^3\) to \(\mathcal O\) containing no index-\(1\) handles and exactly \(c_{\mathrm{poly}}(\mathcal O)\) index-\(3\) handles.
The remaining handles have index \(2\).
Every index-\(3\) handle is associated with a representation obtained by conjugating a faithful representation of a nontrivial finite group by \(\mu\).
These representations remain faithful and have nontrivial images.

Capping the incoming \(S^3\) with a smooth \(4\)-ball gives a compact connected oriented \(4\)-orbifold \(W\) with \(\orbpartial W=\mathcal O\).

Conversely, let \(W\) admit a handle decomposition as in the statement.
Removing the initial smooth \(4\)-ball and reversing the remaining cobordism gives a sequence of handle attachments from \(\mathcal O\) to \(S^3\).
The total number of handles of indices \(1\) and \(3\) is unchanged.
By Theorem~\ref{thm:minimum-polyhedral-handles}, this number is at least \(c_{\mathrm{poly}}(\mathcal O)\).

Finally, a decomposition using only index-\(2\) handles exists if and only if \(c_{\mathrm{poly}}(\mathcal O)=0\).
By Theorem~\ref{thm:minimum-polyhedral-handles}, this is equivalent to \(\mathcal O\) having no polyhedral vertices.
\end{proof}

\appendix
\section{Extension over the spherical boundary models}\label{app:spherical-boundary-extension}

In this appendix, we determine which orientation-preserving orbifold diffeomorphisms of the boundaries extend over the index-\(1\) and index-\(3\) models of Section~\ref{subsec:fixed-region-models}.
Every \(\Sigma_K\) is a sphere with at most three cone points.
Its orientation-preserving mapping classes are therefore determined by the permutations of cone points of equal order.

\begin{proposition}\label{prop:spherical-boundary-extension}
Let \(X\) be one of the orbifolds \(N_K\), \(S^0\times N_K\), \(I(\Sigma_K)\), or \(\mathcal M(K,\nu)\) in Tables~\ref{tab:cyclic-line-regions}--\ref{tab:polyhedral-regions}.
Let \(f\colon\orbpartial X\to\orbpartial X\) be an orientation-preserving orbifold diffeomorphism.

Then \(f\) extends to an orientation-preserving orbifold diffeomorphism of \(X\), subject only to the following conditions.
\begin{enumerate}[label=\textup{(\roman*)}]
\item If \(X=I(\Sigma_K)\), the maps induced by \(f\) on the two ends must determine the same permutation of the cone points in the canonical \(\Sigma_K\) coordinates.
\item If \(X=\mathcal M(\Dih{4},\delta_2)\), the map \(f\) must fix the cone point represented by \((1,0,0)\in S^2\).
\end{enumerate}
Interchange of boundary components is allowed.
For \(X=\mathcal M(\Dih{4},\delta_2)\), the image of the boundary restriction homomorphism has order \(2\) and index \(3\) in \(\MCG^+(\Sigma_{\Dih{4}})\cong S_3\).
\end{proposition}

\begin{proof}
Extension depends only on the boundary mapping class, since boundary isotopies extend over collars.

Every permutation of cone points of equal order on \(\Sigma_K\) is induced by an element of \(N_{\SO(3)}(K)\).
For \(K=\Cyc{q}\), \(q\geq2\), the element \(J\) interchanges the two cone points.
For \(K=\Dih{2q}\), \(q\geq3\), the element \(R_{\pi/q}\) interchanges the two cone points of order \(2\).
For \(\Dih{2}=\langle J\rangle\), use \(R_\pi\).
For \(\Dih{4}\), the orientation-preserving signed coordinate permutations induce all permutations of the three cone points.
For \(\mathbf T\), the element \(\omega_{\mathbf T}\) interchanges the two cone points of order \(3\).
The groups \(\MCG^+(\Sigma_K)\) are trivial for \(K=1,\mathbf O,\mathbf I\).
Each of these elements induces an orbifold diffeomorphism of \(N_K\).

The map induced by \(\alpha=-I_3\) on \(N_K\) reverses orientation.
Compose it with a map induced by an element of \(N_{\SO(3)}(K)\) that reverses its permutation of cone points.
The composition reverses orientation and fixes every cone point of \(\Sigma_K\).
Consequently, every orientation-preserving boundary mapping class of \(S^0\times N_K\) extends over its two components.
An interchange of the components is realized by orientation-reversing diffeomorphisms between the two copies of \(N_K\).

In \(I(\Sigma_K)\), each singular arc joins corresponding cone points on the two boundary components.
An extending diffeomorphism must therefore induce the same permutation at both ends.
Conversely, the product of the identity on \(D^1\) with a diffeomorphism of \(\Sigma_K\) realizes any such permutation.
Reversing the interval coordinate and applying an orientation-reversing diffeomorphism of \(\Sigma_K\) fixing the cone points realizes an interchange of the boundary components.

For \(\mathcal M(K,\nu)\), let \(h\in N_{\SO(3)}(K)\).
If \(h\nu h^{-1}\nu^{-1}\in K\), then the map
\[
(t,x)\longmapsto(t,hx)
\]
normalizes \(\widehat K_\nu\) and induces an orientation-preserving orbifold diffeomorphism of \(\mathcal M(K,\nu)\).

For \(\nu=\alpha\), every \(h\in N_{\SO(3)}(K)\) satisfies this condition.
For \((K,\nu)=(\Cyc{q},\gamma)\) or \((\Cyc{q},\delta_q)\), take \(h=J\).
For \((K,\nu)=(\Dih{2q},\delta_q)\), take \(h=R_{\pi/q}\) when \(q\geq3\), and \(h=R_\pi\) when \(q=1\).
For \((K,\nu)=(\mathbf T,\varepsilon_{\mathbf T})\), take \(h=\omega_{\mathbf T}\).
These maps realize every orientation-preserving boundary mapping class of the models in the statement, except for \(\mathcal M(\Dih{4},\delta_2)\).

For \(K=\Dih{4}=\langle R_\pi,J\rangle\), we have
\[
\delta_2^2=R_\pi,
\qquad
\delta_2J\delta_2^{-1}=R_\pi J.
\]
Since \(D^1\times S^2\) is simply connected,
\[
\pi_1^{\mathrm{orb}}(\mathcal M(\Dih{4},\delta_2))
=
\langle K,\delta_2\rangle
\cong\Dih{8}.
\]
The center of this group is \(\{1,R_\pi\}\).
The meridians of the three boundary cone points have images represented by \(R_\pi,J,R_\pi J\), respectively.
The cone point represented by \((1,0,0)\in S^2\) corresponds to \(R_\pi\).
Every extending diffeomorphism preserves the unique nonidentity central element and therefore fixes this cone point.
Conversely, \(R_{\pi/2}\) normalizes \(K\), commutes with \(\delta_2\), and interchanges the other two cone points.
Thus the boundary mapping classes that extend are precisely those fixing the cone point represented by \((1,0,0)\).
\end{proof}

\clearpage
\section{Diagrams of the boundary regions}\label{subsec:illustrated-boundary-models}

\begin{figure}[!htbp]
\centering
\includegraphics[height=35mm,width=0.85\linewidth,keepaspectratio]{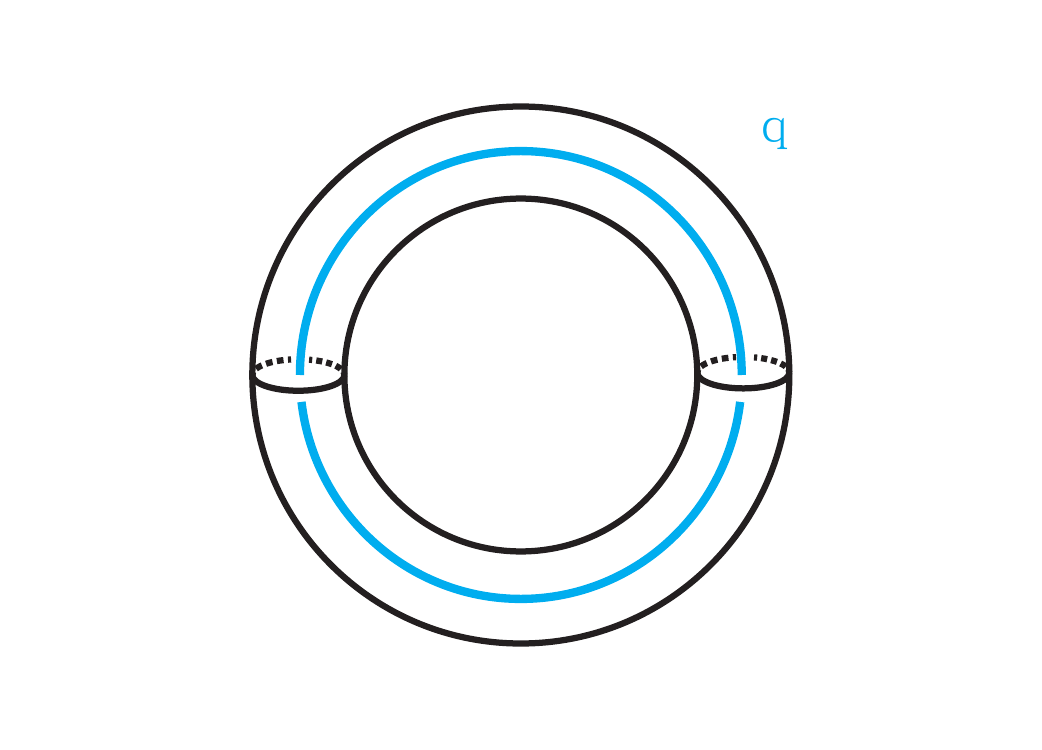}
\caption{The solid torus \(T_q\) for \(q\geq1\). The core is regular for \(q=1\).}\label{fig:G1-2h}
\end{figure}

\begin{figure}[!htbp]
\centering
\begin{tabular}{@{}c@{\hspace{6mm}}c@{}}
\includegraphics[height=35mm,width=0.45\linewidth,keepaspectratio]{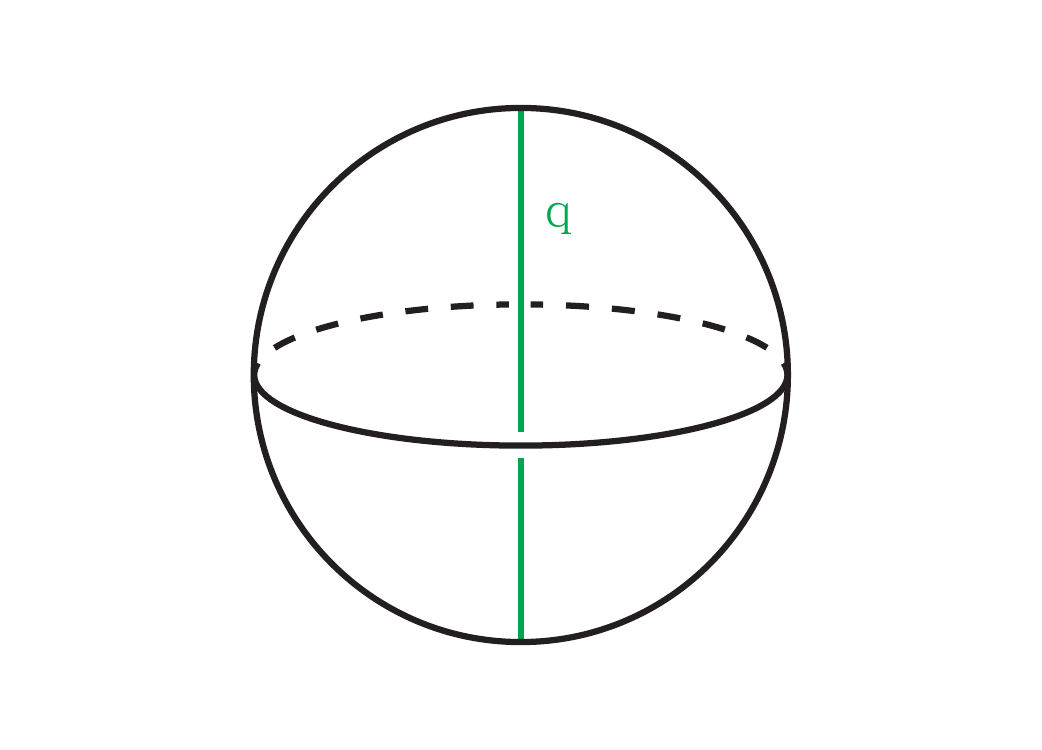} &
\HandleDiagram{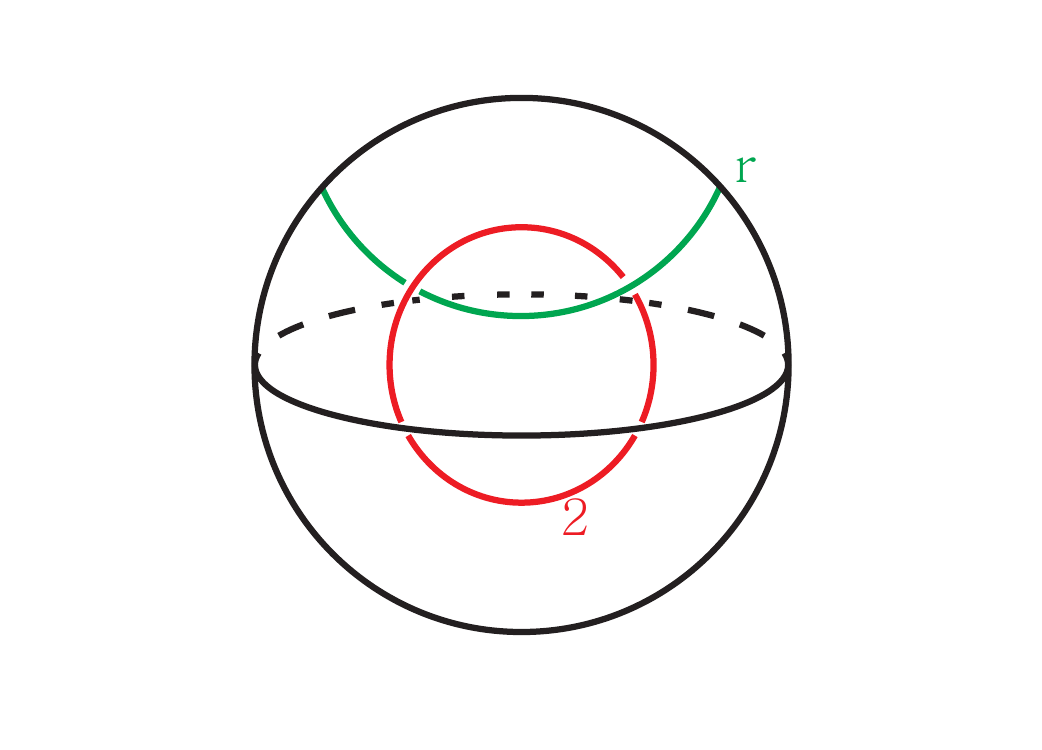} \\[-1mm]
{\small\textup{(a)}} & {\small\textup{(b)}}
\end{tabular}
\caption{\textup{(a)} \(N_{\Cyc{q}}\) for \(q\geq2\). \textup{(b)} \(\mathcal M(\Cyc{r},\alpha)\) for even \(r\geq2\).}\label{fig:G1-1h-rneq2}
\end{figure}

\begin{figure}[!htbp]
\centering
\includegraphics[height=35mm,width=0.85\linewidth,keepaspectratio]{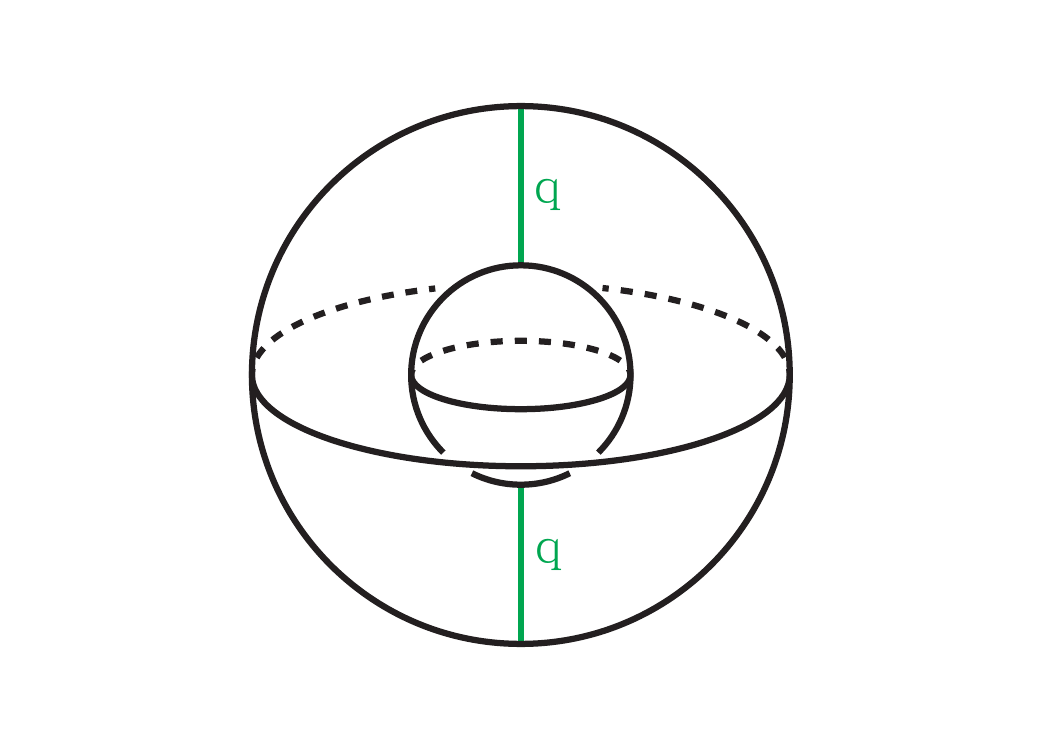}
\caption{The model \(I(\Sigma_{\Cyc{q}})\) for \(q\geq2\).}\label{fig:G1p-1h-even}
\end{figure}

\begin{figure}[!htbp]
\centering
\HandleDiagram{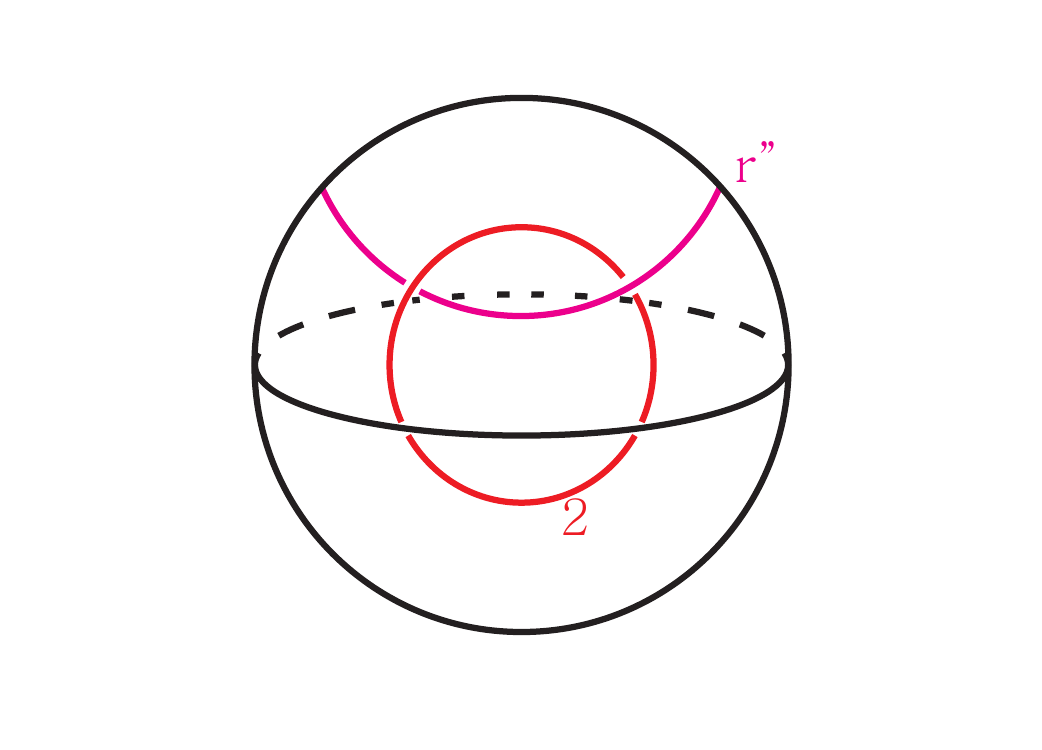}
\caption{The model \(\mathcal M(\Cyc{r''},\delta_{r''})\) for odd \(r''\geq3\).}\label{fig:G1p-1h-odd}
\end{figure}

\begin{figure}[!htbp]
\centering
\begin{tabular}{@{}c@{\hspace{6mm}}c@{}}
\HandleDiagram{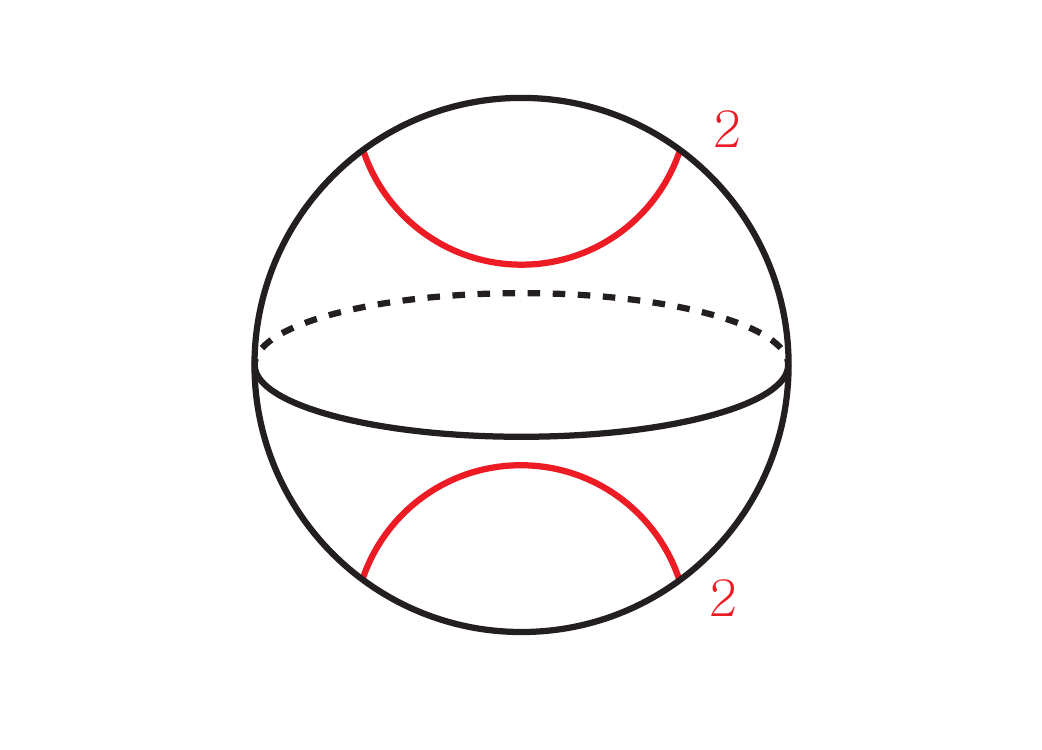} &
\HandleDiagram{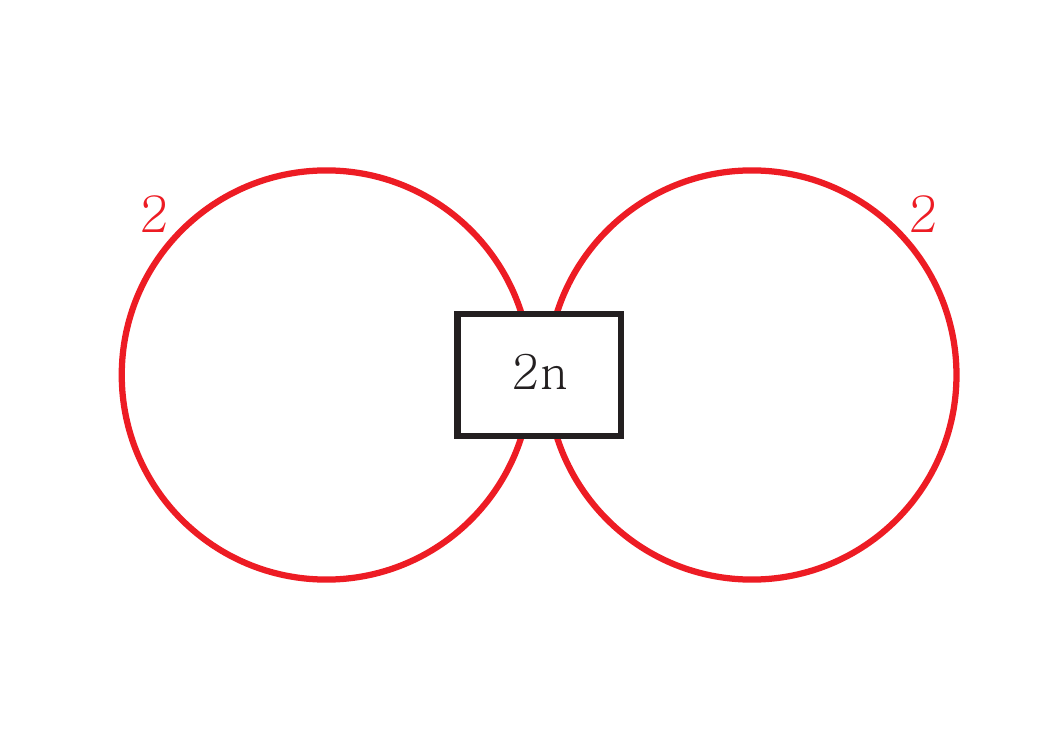} \\[-1mm]
{\small\textup{(a)}} & {\small\textup{(b)}}
\end{tabular}
\caption{\textup{(a)} The solid pillow \(\mathcal P_1\). \textup{(b)} The singular link of \(S^3\sslash G_3\).}\label{fig:G3-2h}
\end{figure}

\begin{figure}[!htbp]
\centering
\includegraphics[height=35mm,width=0.85\linewidth,keepaspectratio]{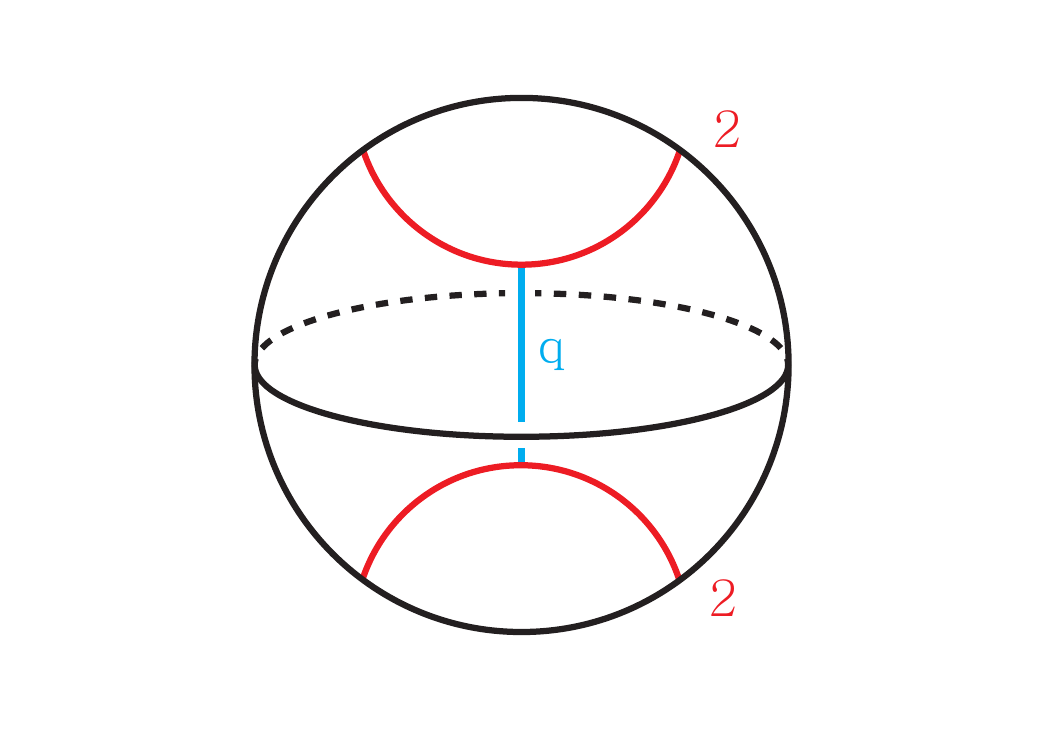}
\caption{The solid pillow \(\mathcal P_q\) for \(q\geq2\).}\label{fig:G11-2h}
\end{figure}

\begin{figure}[!htbp]
\centering
\HandleDiagram{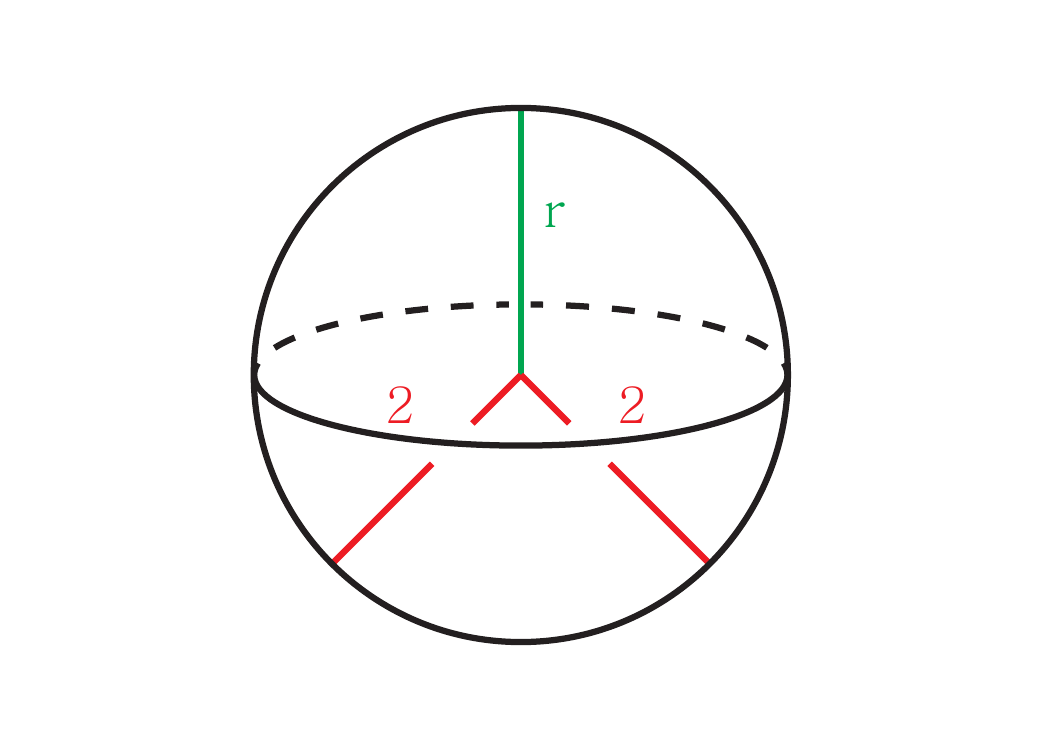}
\caption{The model \(N_{\Dih{2r}}\) for \(r\geq2\).}\label{fig:G11-1h-att}
\end{figure}

\begin{figure}[!htbp]
\centering
\begin{tabular}{@{}c@{\hspace{6mm}}c@{}}
\HandleDiagram{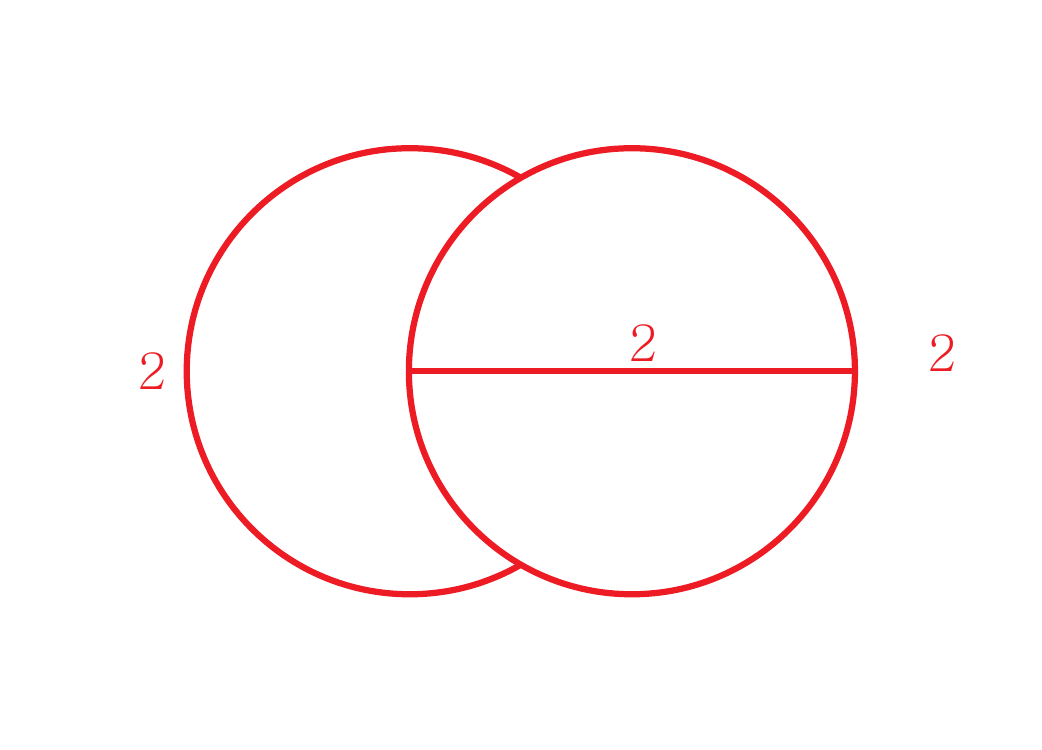} &
\HandleDiagram{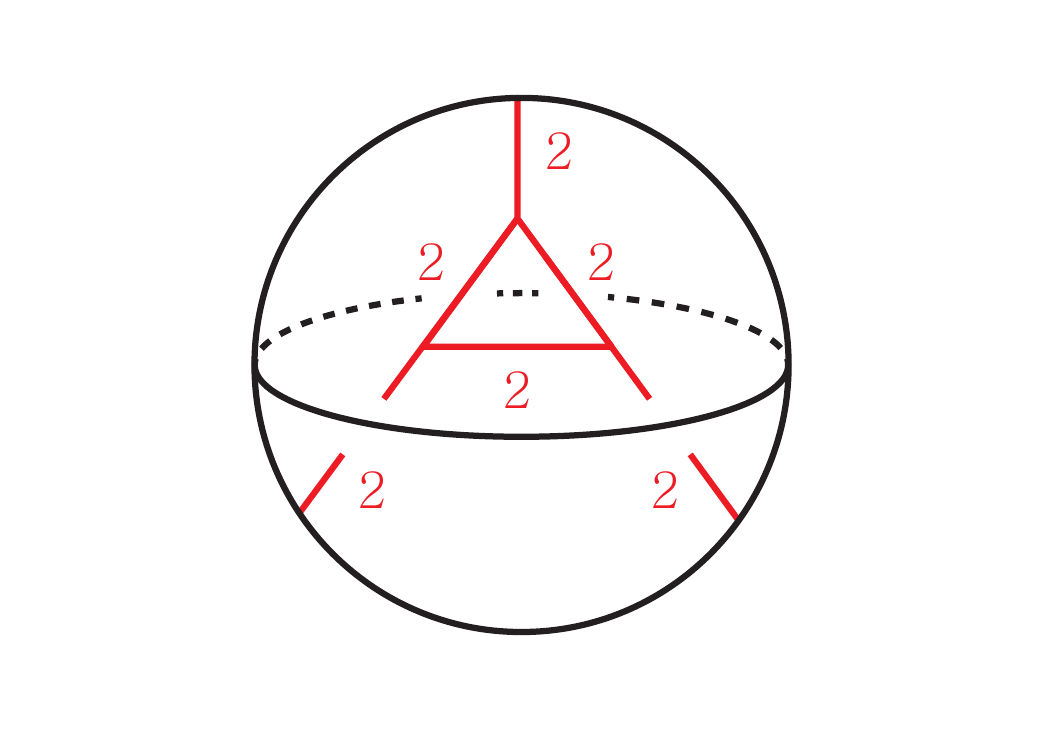} \\[-1mm]
{\small\textup{(a)}} & {\small\textup{(b)}}
\end{tabular}
\caption{\textup{Family} \(11\) for \(r=2\). \textup{(a)} The singular set of \(S^3\sslash G_{11}\). \textup{(b)} \(\mathcal M(\Dih{4},\alpha)\).}\label{fig:G11-r2}
\end{figure}

\begin{figure}[!htbp]
\centering
\begin{tabular}{@{}c@{\hspace{6mm}}c@{}}
\HandleDiagram{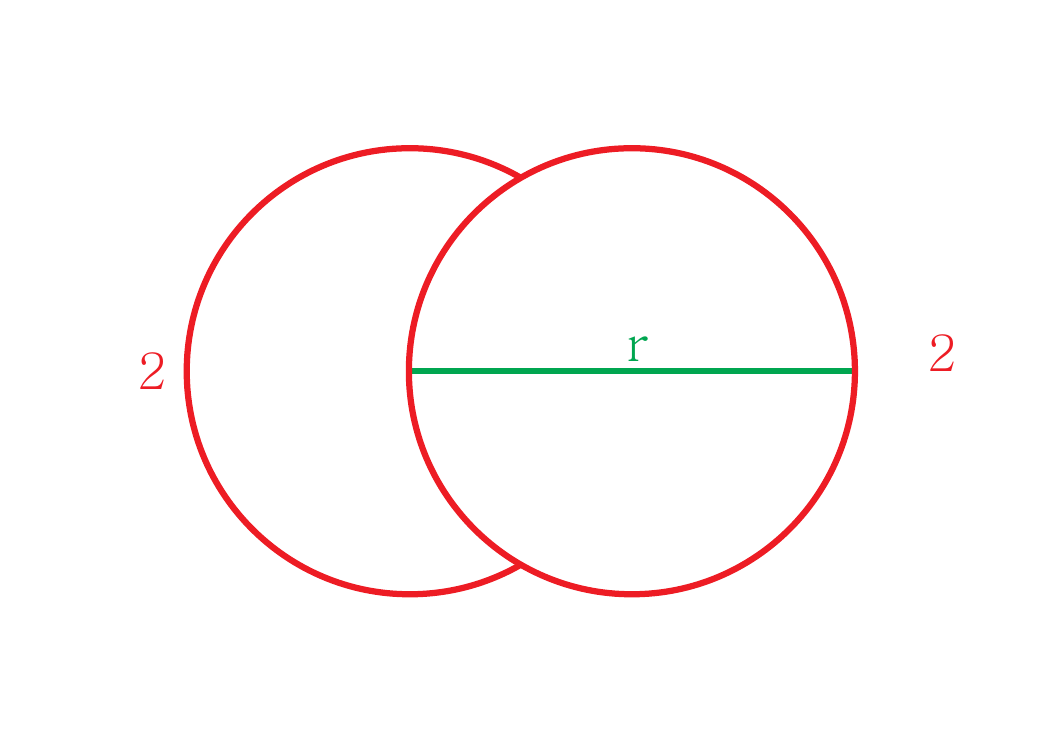} &
\HandleDiagram{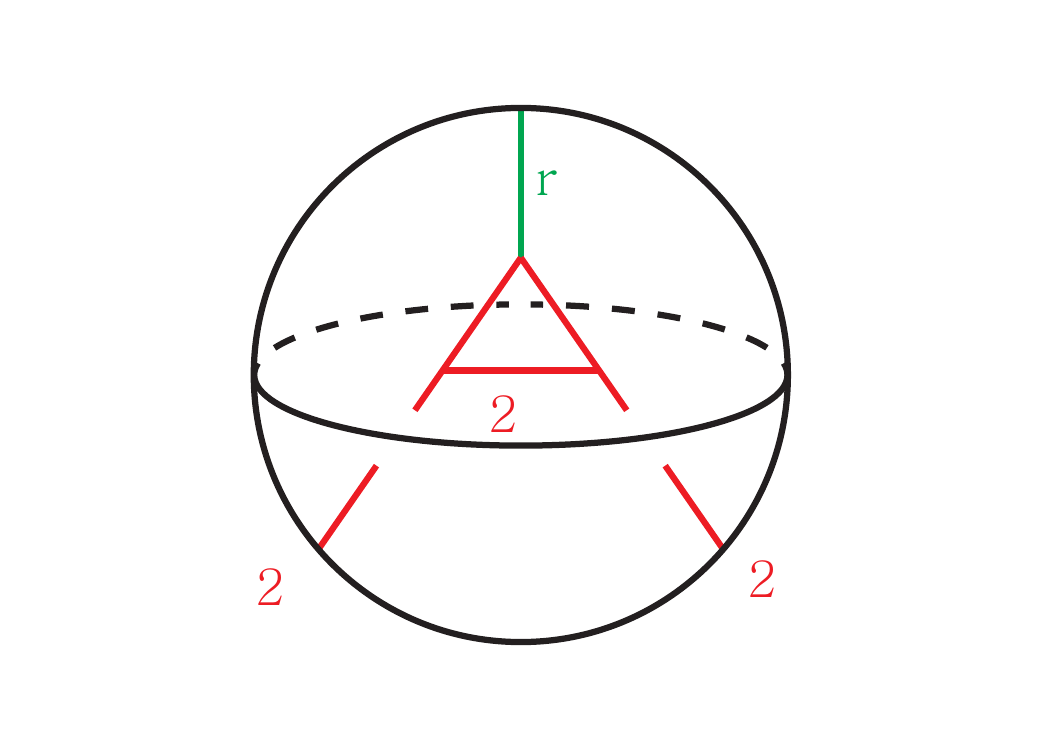} \\[-1mm]
{\small\textup{(a)}} & {\small\textup{(b)}}
\end{tabular}
\caption{\textup{Family} \(11\) for even \(r\geq4\). \textup{(a)} The singular set of \(S^3\sslash G_{11}\). \textup{(b)} \(\mathcal M(\Dih{2r},\alpha)\).}\label{fig:G11-r-even}
\end{figure}

\begin{figure}[!htbp]
\centering
\begin{tabular}{@{}c@{\hspace{6mm}}c@{}}
\HandleDiagram{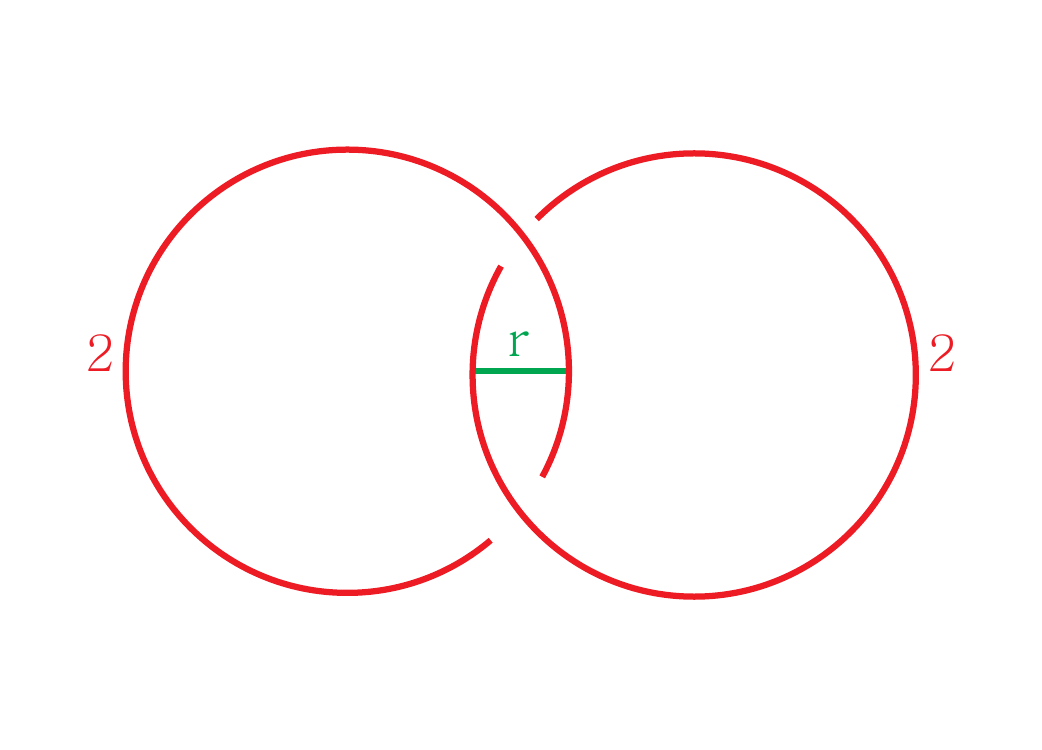} &
\HandleDiagram{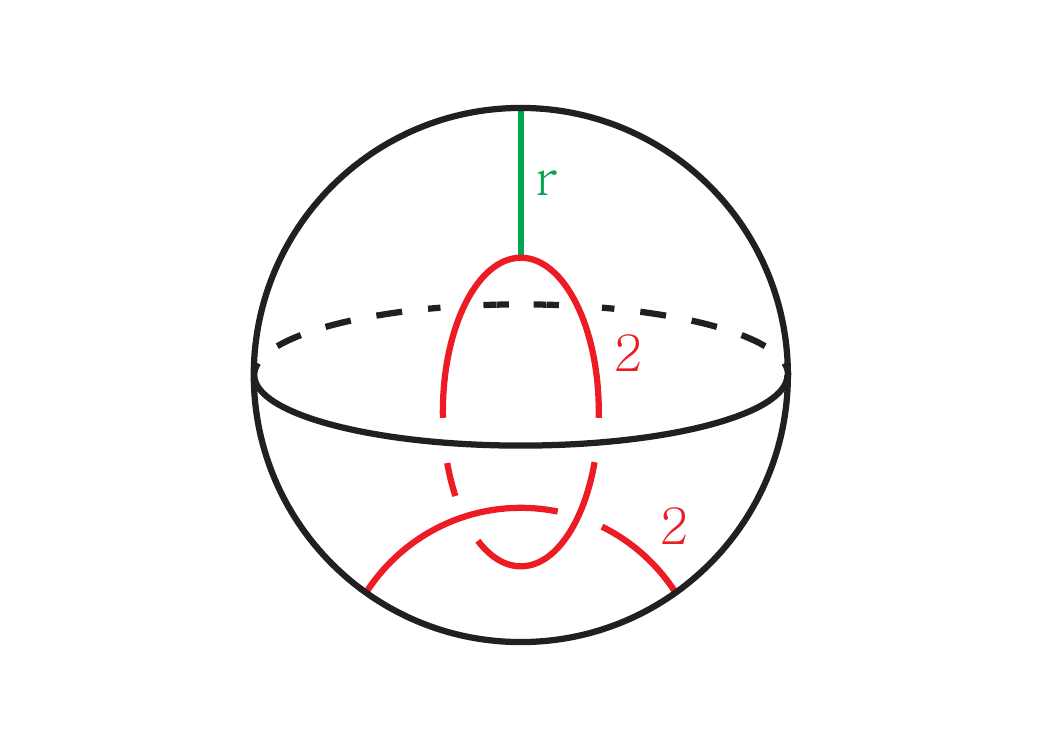} \\[-1mm]
{\small\textup{(a)}} & {\small\textup{(b)}}
\end{tabular}
\caption{\textup{Family} \(11\) for odd \(r\geq3\). \textup{(a)} The singular set of \(S^3\sslash G_{11}\). \textup{(b)} \(\mathcal M(\Dih{2r},\alpha)\).}\label{fig:G11-r-odd}
\end{figure}

\begin{figure}[!htbp]
\centering
\HandleDiagram{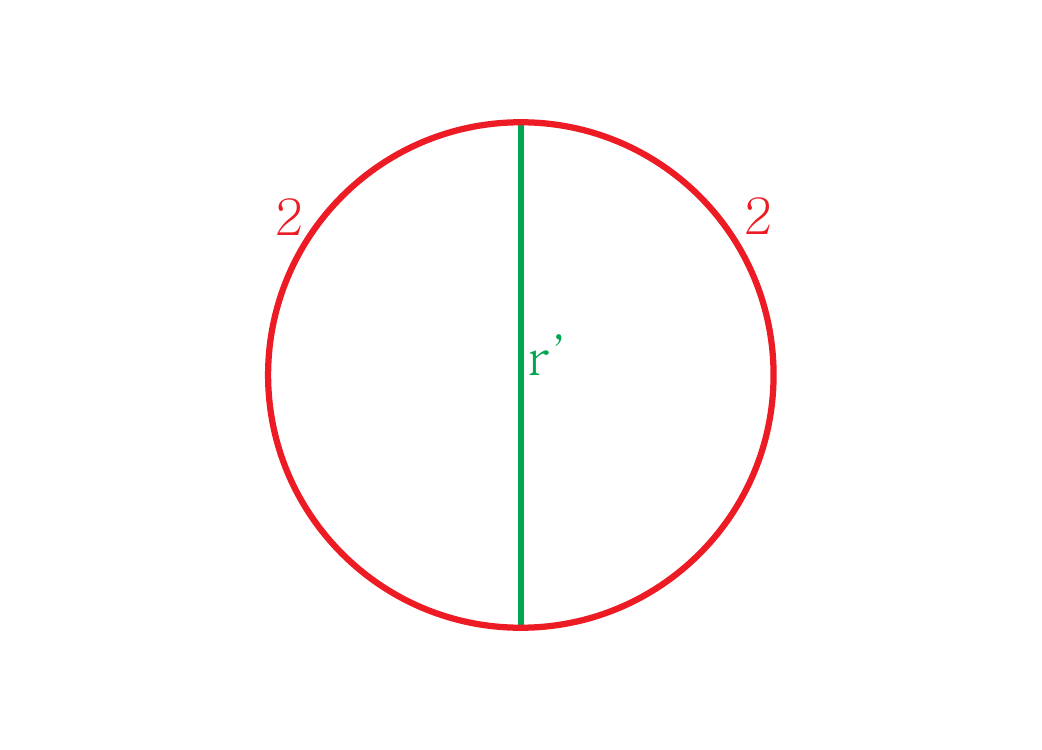}
\caption{The singular set of \(S^3\sslash G_{11'}\) for \(r'>2\) and even \(\ell'\).}\label{fig:G11p-even-ambient}
\end{figure}

\begin{figure}[!htbp]
\centering
\begin{tabular}{@{}c@{\hspace{6mm}}c@{}}
\HandleDiagram{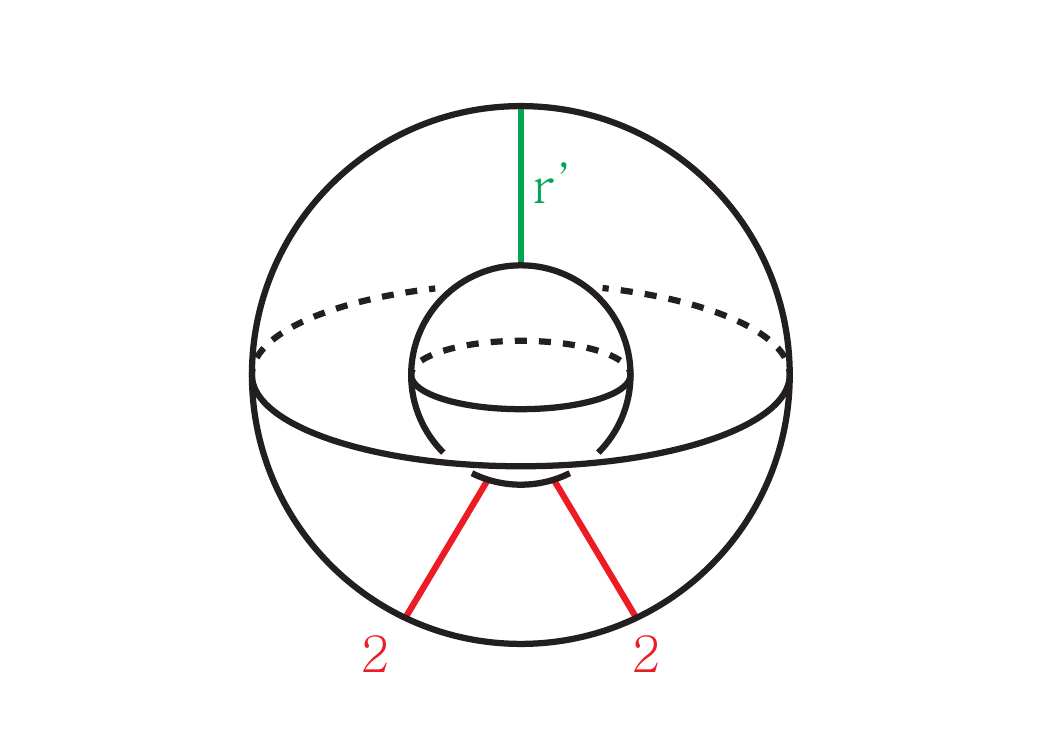} &
\HandleDiagram{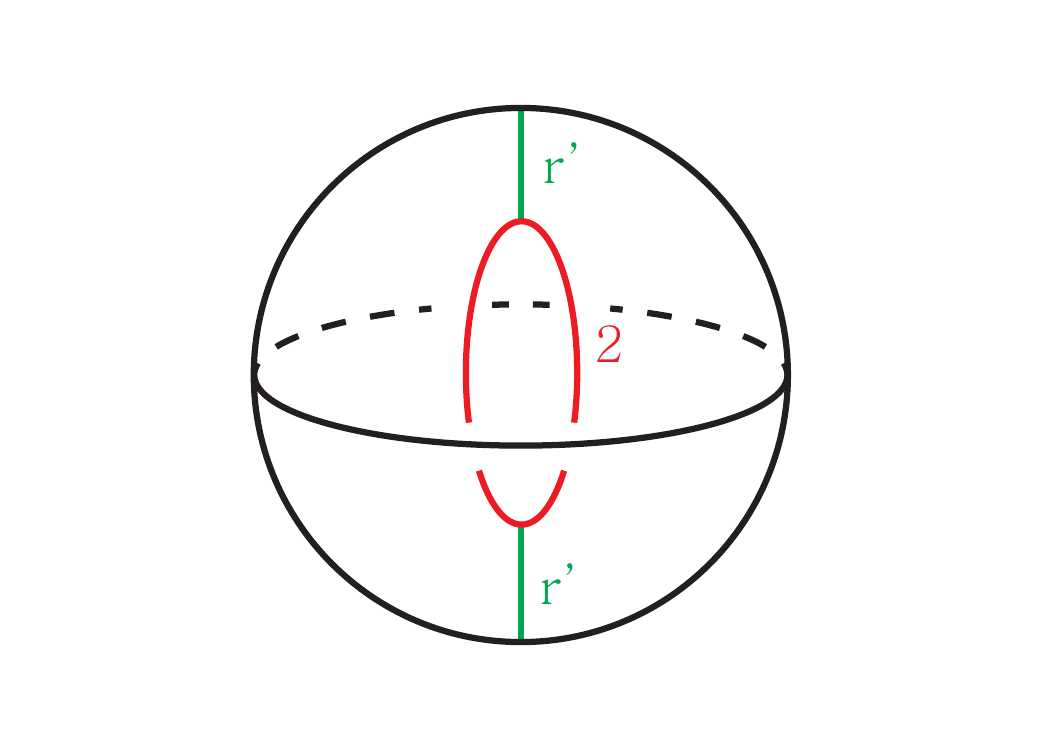} \\[-1mm]
{\small\textup{(a)}} & {\small\textup{(b)}}
\end{tabular}
\caption{Belt regions for \(G_{11'}\) with \(r'>2\) and even \(\ell'\). \textup{(a)} \(I(\Sigma_{\Dih{2r'}})\) for \(L_0\). \textup{(b)} \(\mathcal M(\Cyc{r'},\gamma)\) for \(L_1\).}\label{fig:G11p-even-belt}
\end{figure}

\begin{figure}[!htbp]
\centering
\begin{tabular}{@{}c@{\hspace{6mm}}c@{}}
\HandleDiagram{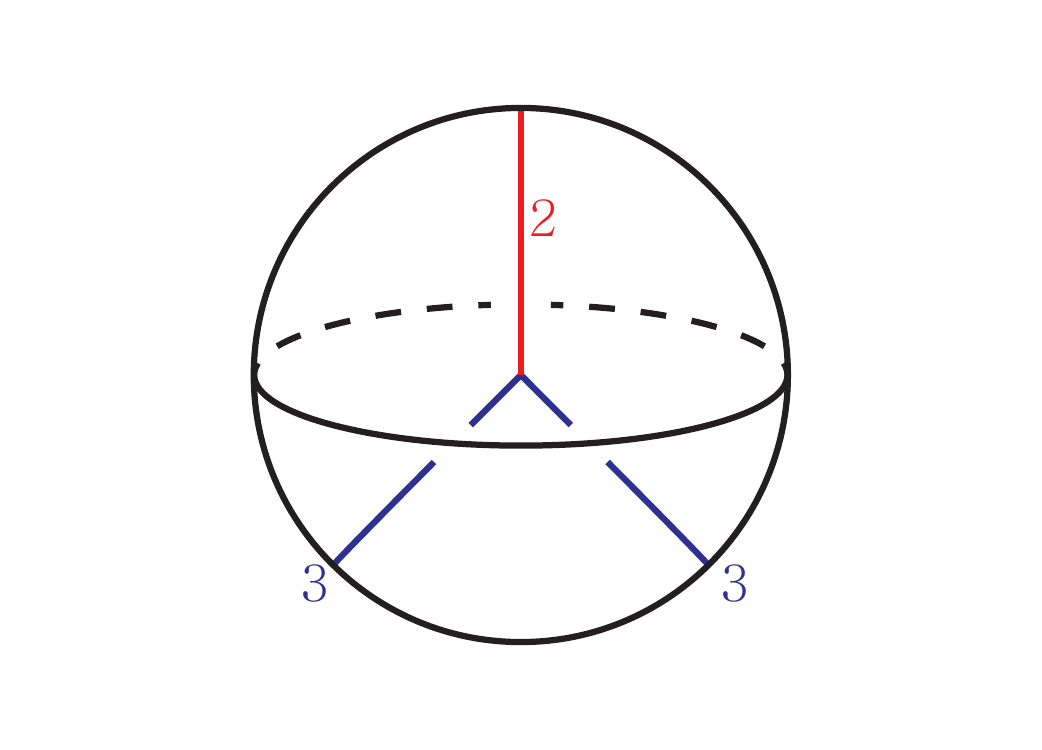} &
\HandleDiagram{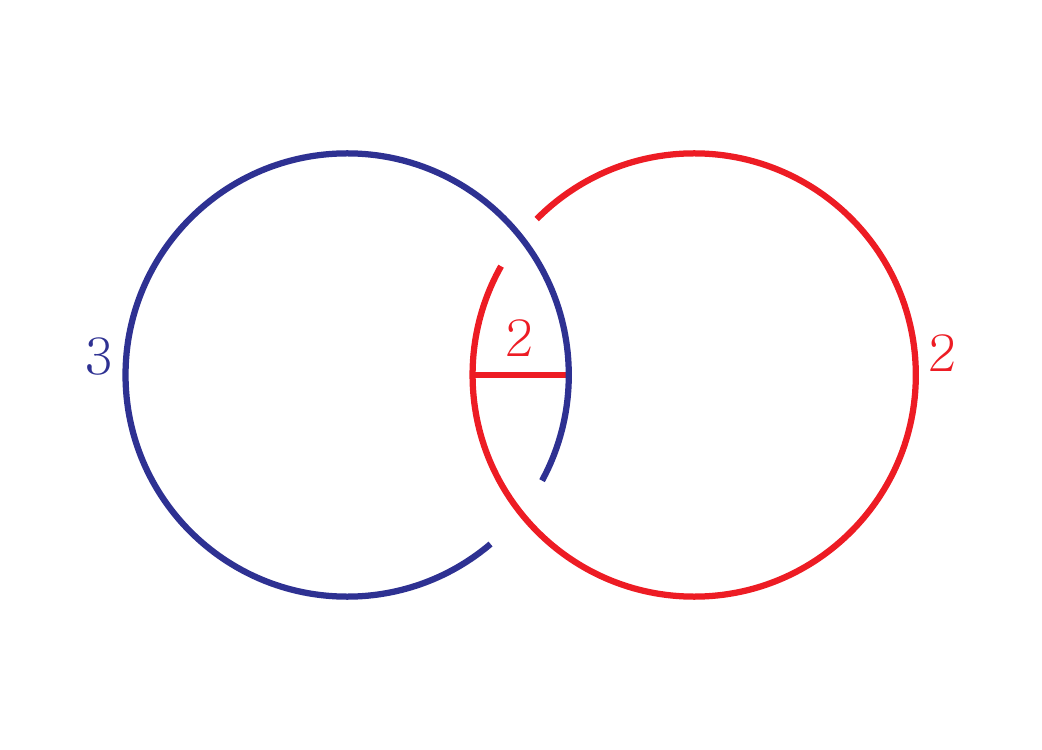} \\[-1mm]
{\small\textup{(a)}} & {\small\textup{(b)}} \\[2mm]
\multicolumn{2}{c}{\HandleDiagram{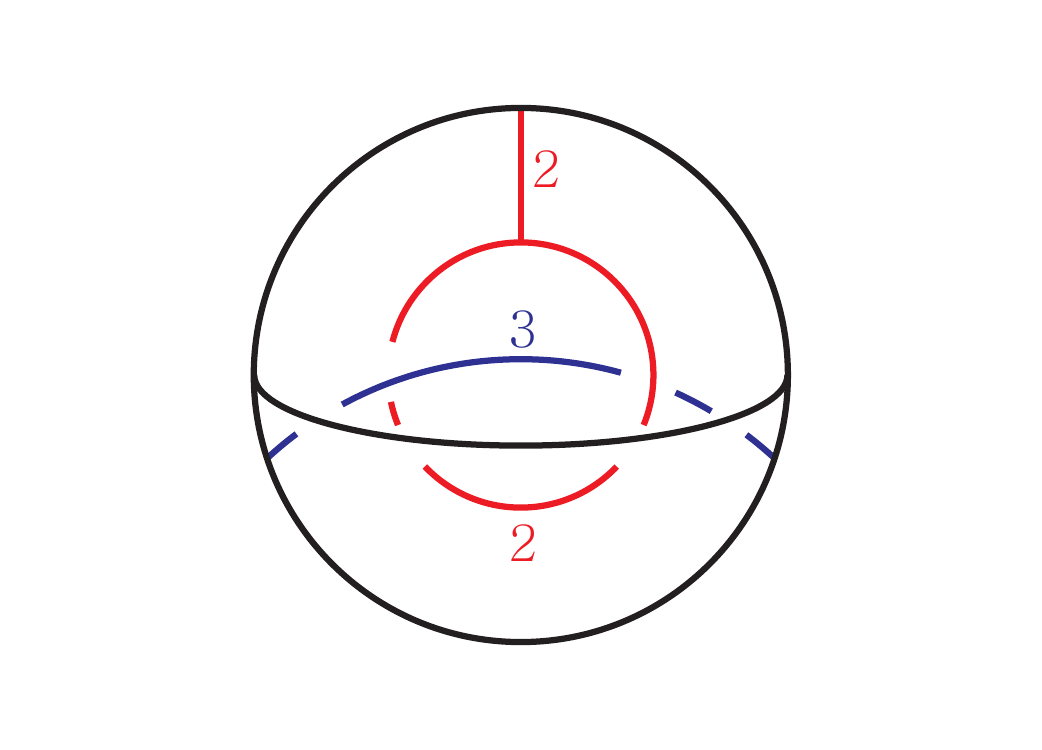}} \\[-1mm]
\multicolumn{2}{c}{{\small\textup{(c)}}}
\end{tabular}
\caption{\textup{(a)} \(N_{\mathbf T}\). \textup{(b)} The labeled singular graph of \(S^3\sslash G_{21}\). \textup{(c)} \(\mathcal M(\mathbf T,\alpha)\).}\label{fig:G21-1h}
\end{figure}

\begin{figure}[!htbp]
\centering
\begin{tabular}{@{}c@{\hspace{6mm}}c@{}}
\HandleDiagram{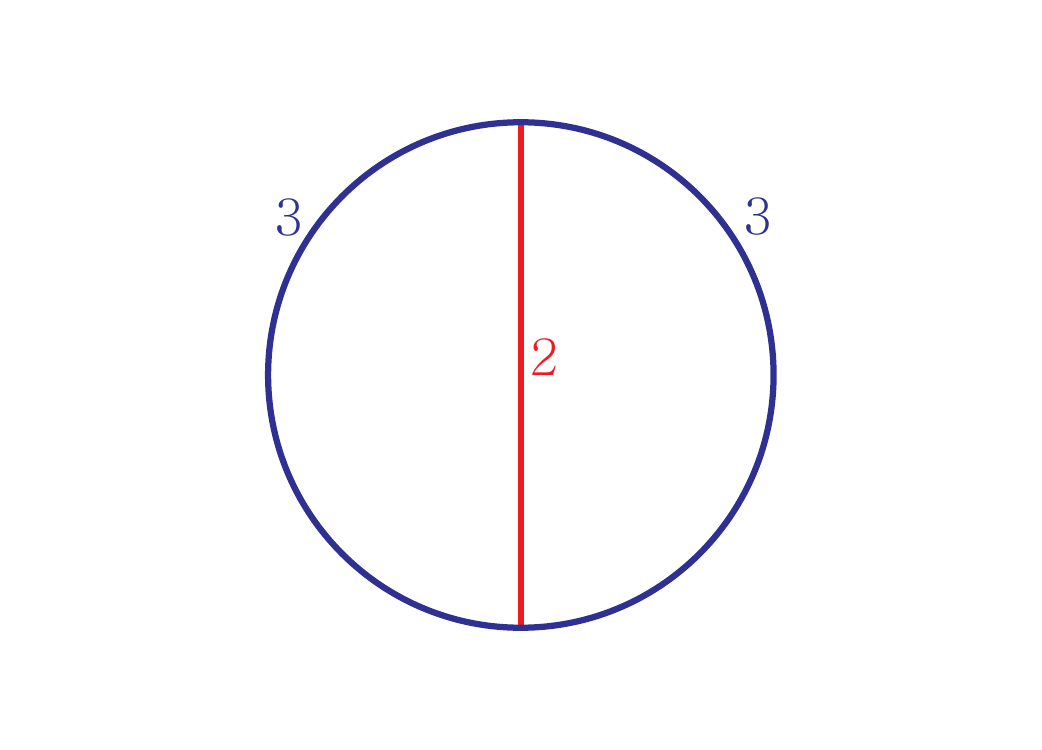} &
\HandleDiagram{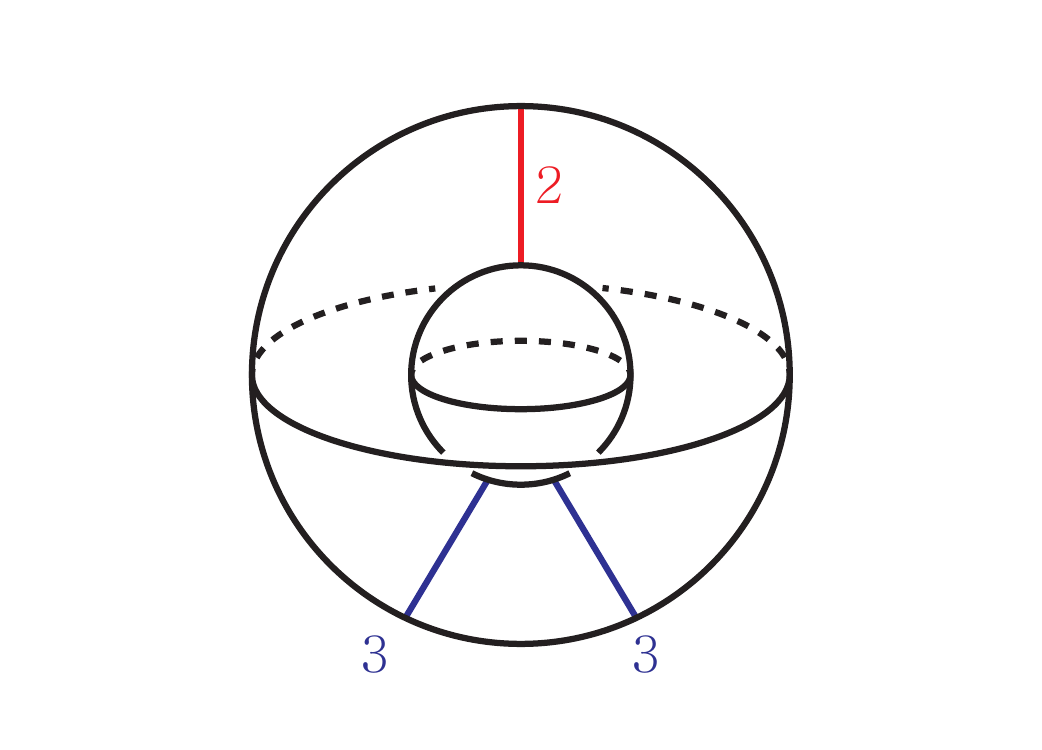} \\[-1mm]
{\small\textup{(a)}} & {\small\textup{(b)}}
\end{tabular}
\caption{\textup{(a)} The labeled singular graph of \(S^3\sslash G_{21'}\). \textup{(b)} \(I(\Sigma_{\mathbf T})\).}\label{fig:G21p-1h}
\end{figure}

\begin{figure}[!htbp]
\centering
\begin{tabular}{@{}c@{\hspace{6mm}}c@{}}
\HandleDiagram{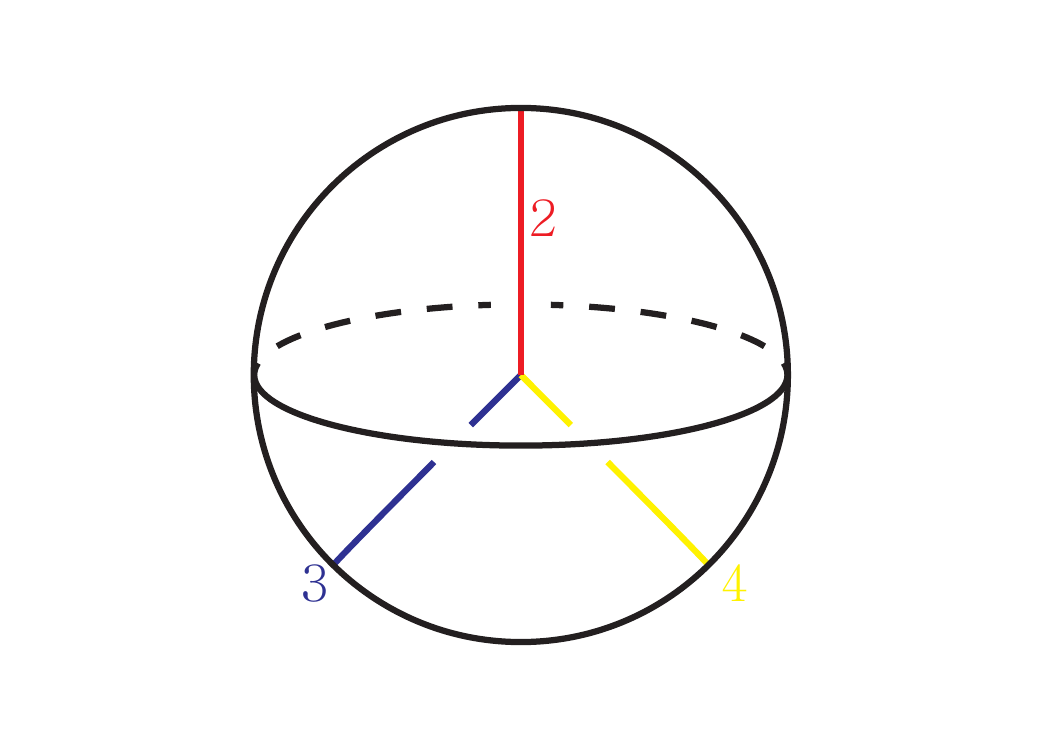} &
\HandleDiagram{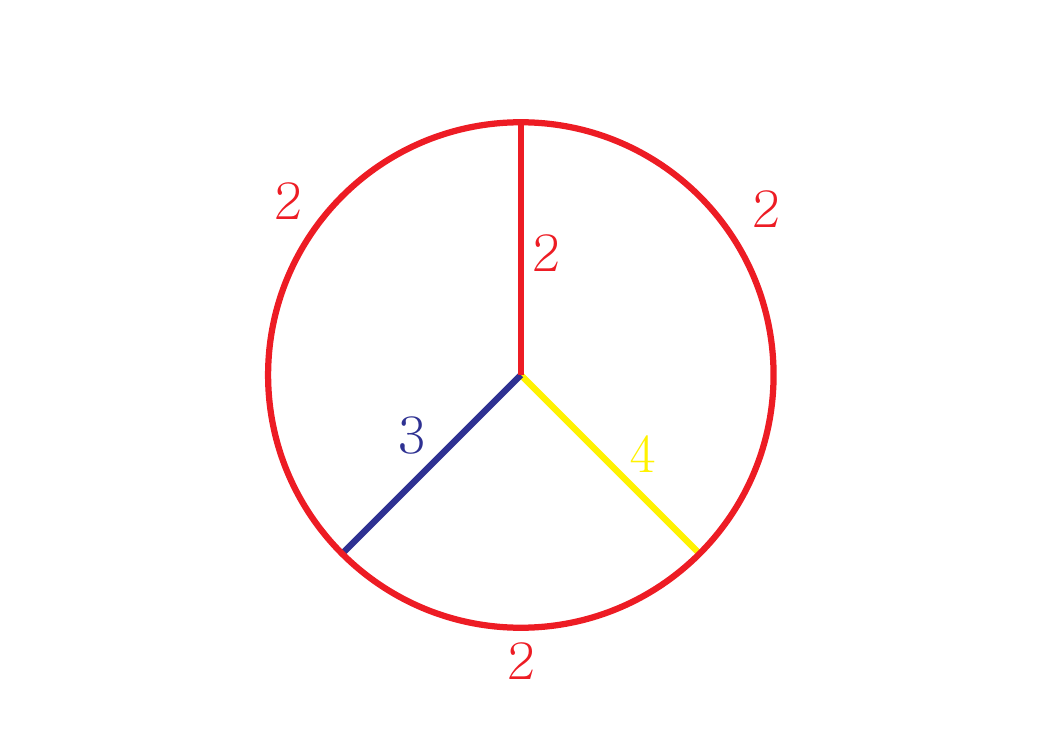} \\[-1mm]
{\small\textup{(a)}} & {\small\textup{(b)}} \\[2mm]
\multicolumn{2}{c}{\HandleDiagram{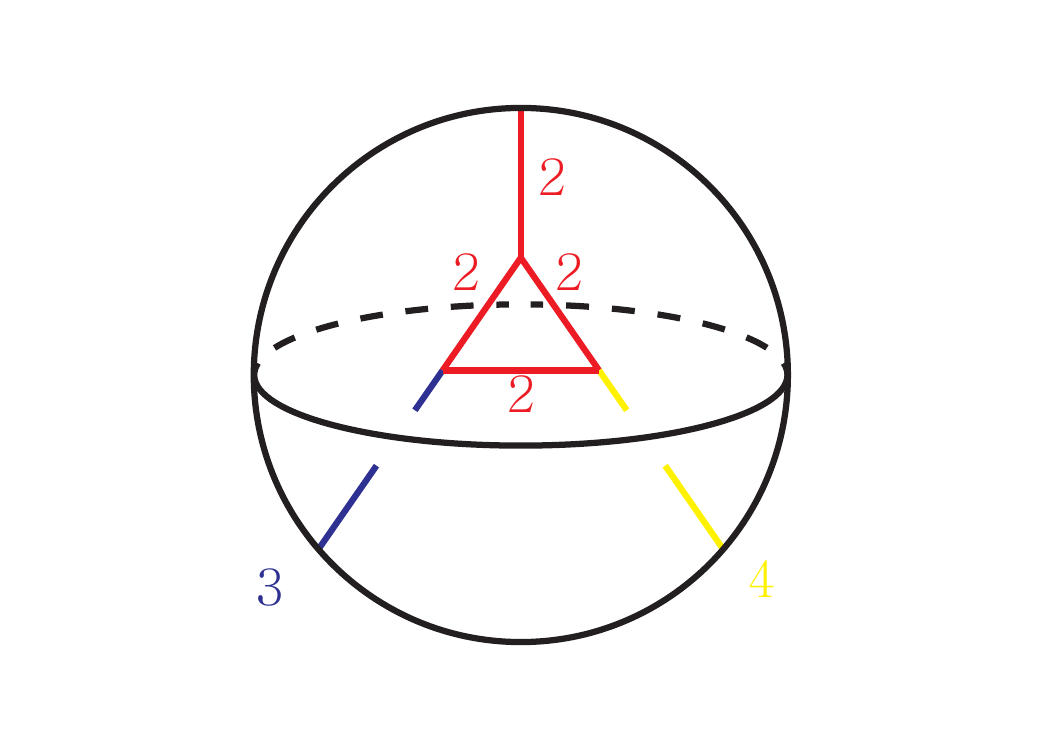}} \\[-1mm]
\multicolumn{2}{c}{{\small\textup{(c)}}}
\end{tabular}
\caption{\textup{(a)} \(N_{\mathbf O}\). \textup{(b)} The labeled singular graph of \(S^3\sslash G_{26}\). \textup{(c)} \(\mathcal M(\mathbf O,\alpha)\).}\label{fig:G26-1h}
\end{figure}

\begin{figure}[!htbp]
\centering
\begin{tabular}{@{}c@{\hspace{6mm}}c@{}}
\HandleDiagram{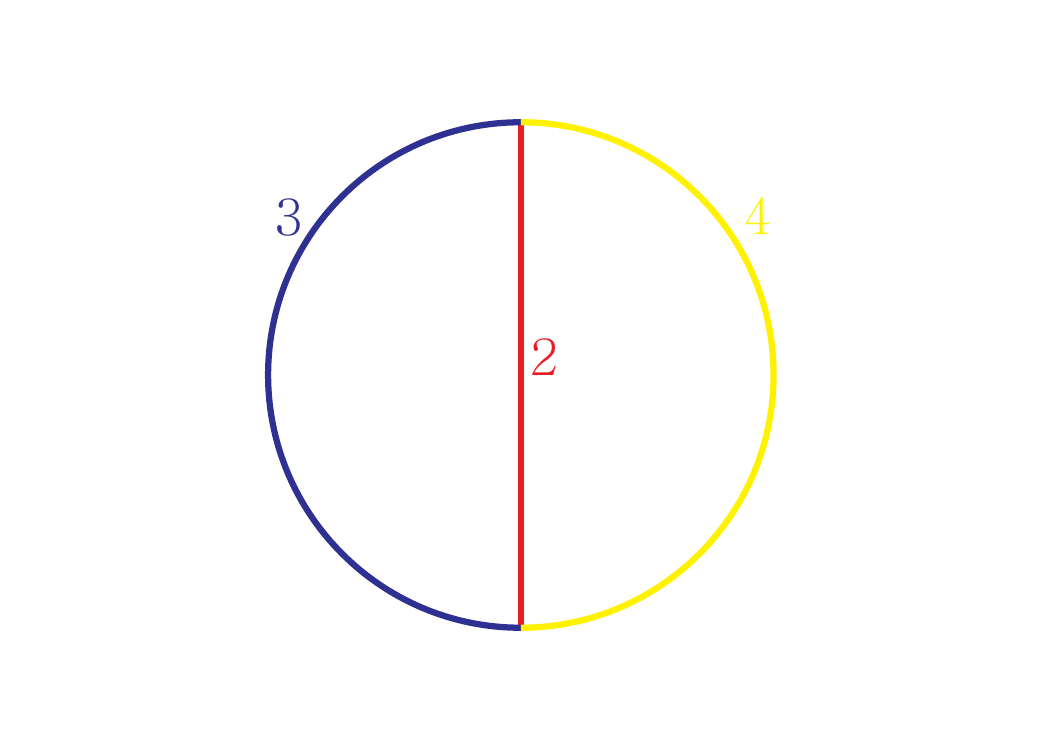} &
\HandleDiagram{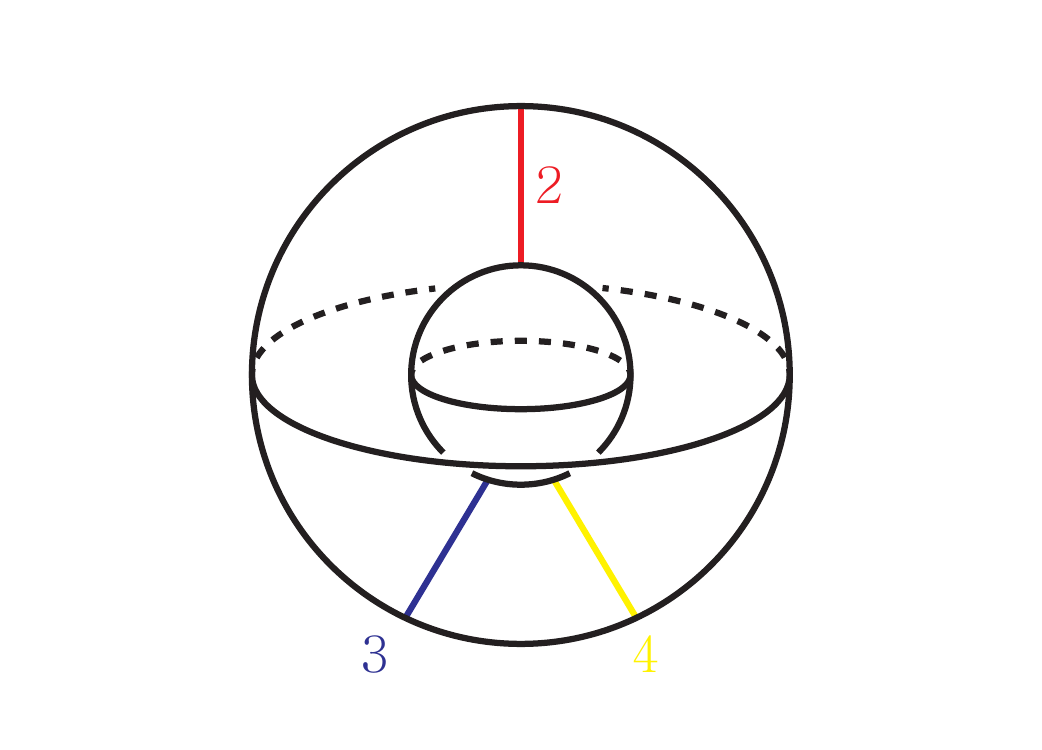} \\[-1mm]
{\small\textup{(a)}} & {\small\textup{(b)}}
\end{tabular}
\caption{\textup{(a)} The labeled singular graph of \(S^3\sslash G_{26'}\). \textup{(b)} \(I(\Sigma_{\mathbf O})\).}\label{fig:G26p-1h}
\end{figure}

\begin{figure}[!htbp]
\centering
\begin{tabular}{@{}c@{\hspace{6mm}}c@{}}
\HandleDiagram{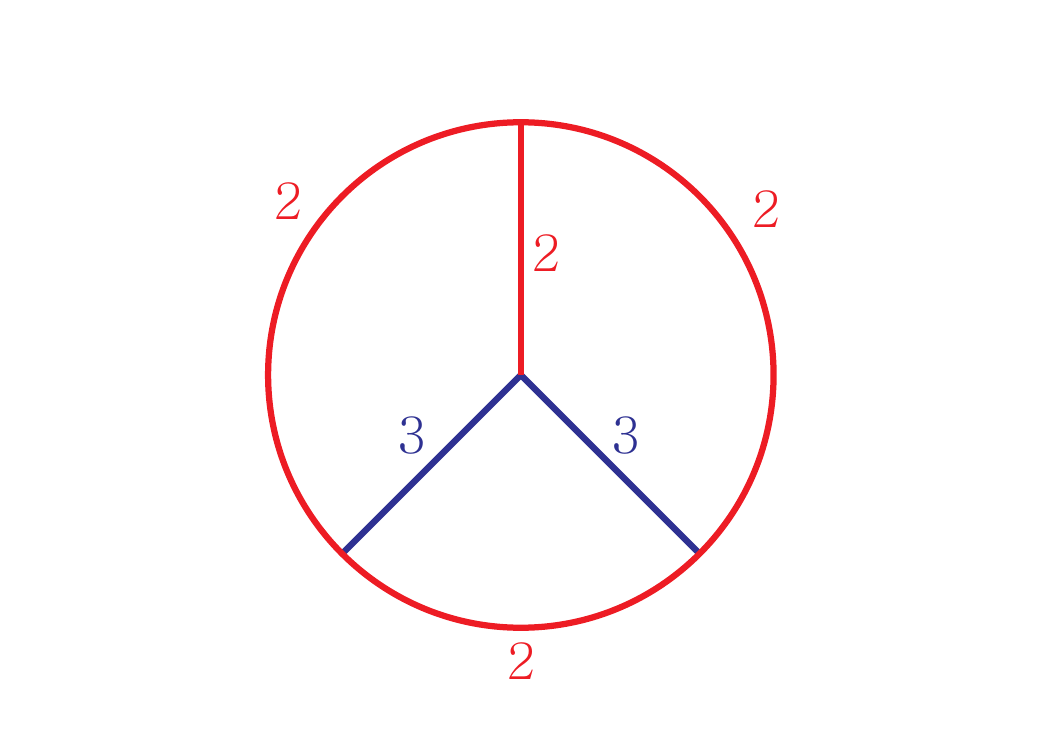} &
\HandleDiagram{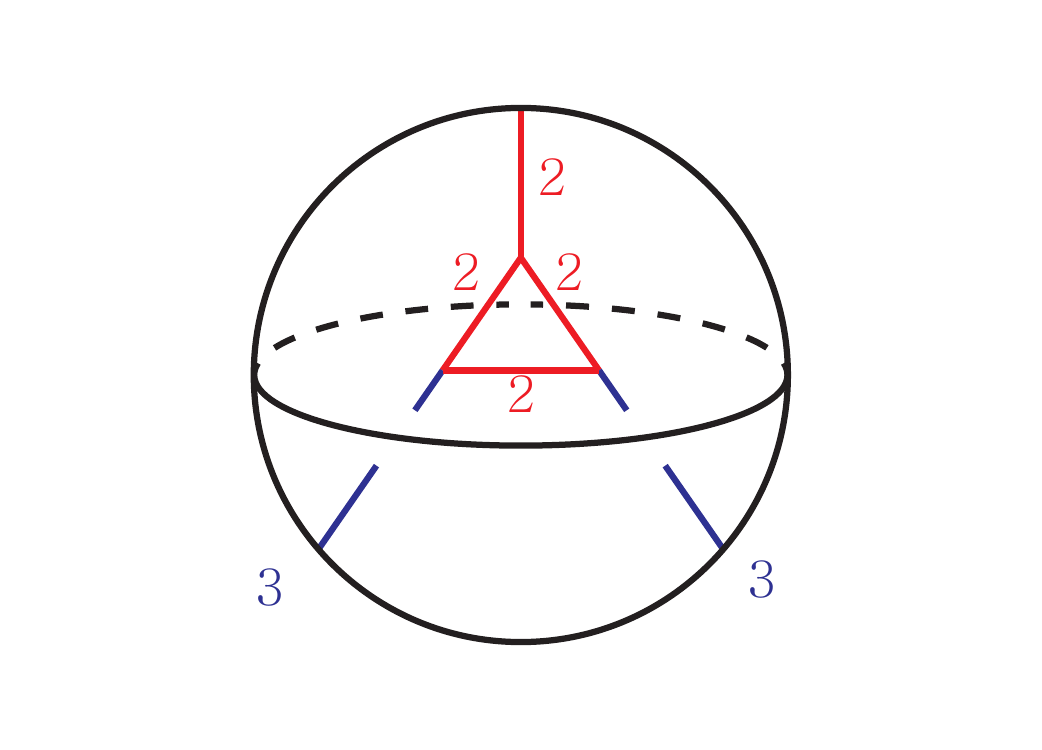} \\[-1mm]
{\small\textup{(a)}} & {\small\textup{(b)}}
\end{tabular}
\caption{\textup{(a)} The labeled singular graph of \(S^3\sslash G_{26''}\). \textup{(b)} \(\mathcal M(\mathbf T,\varepsilon_{\mathbf T})\).}\label{fig:G26pp-1h}
\end{figure}

\begin{figure}[!htbp]
\centering
\begin{tabular}{@{}c@{\hspace{6mm}}c@{}}
\HandleDiagram{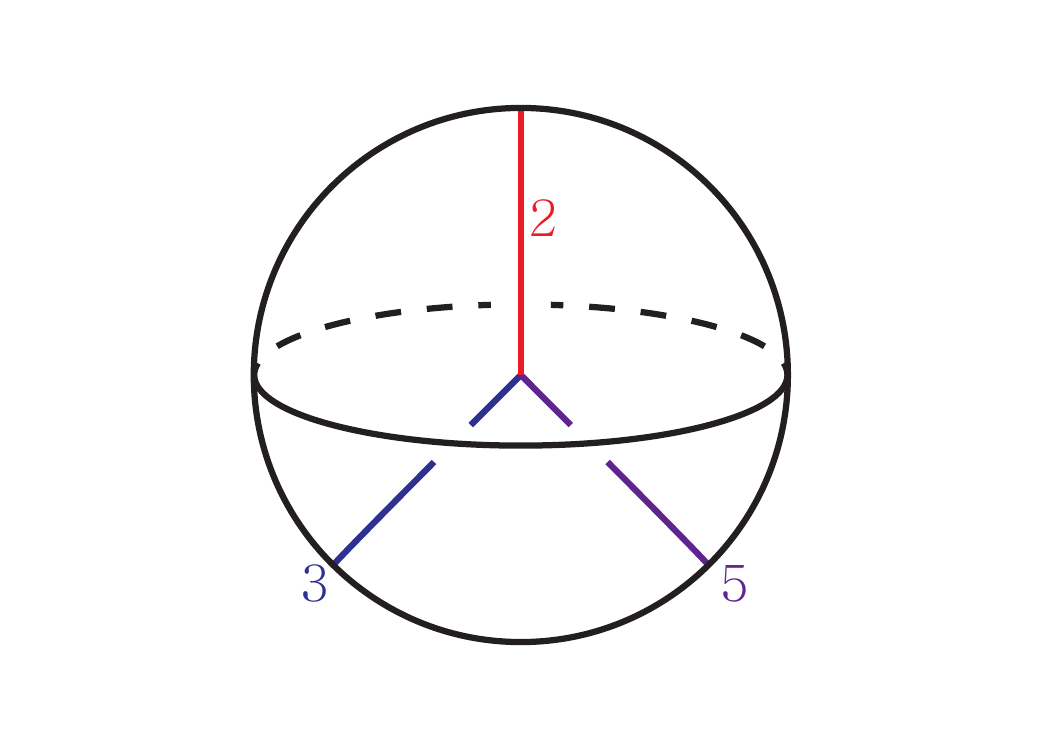} &
\HandleDiagram{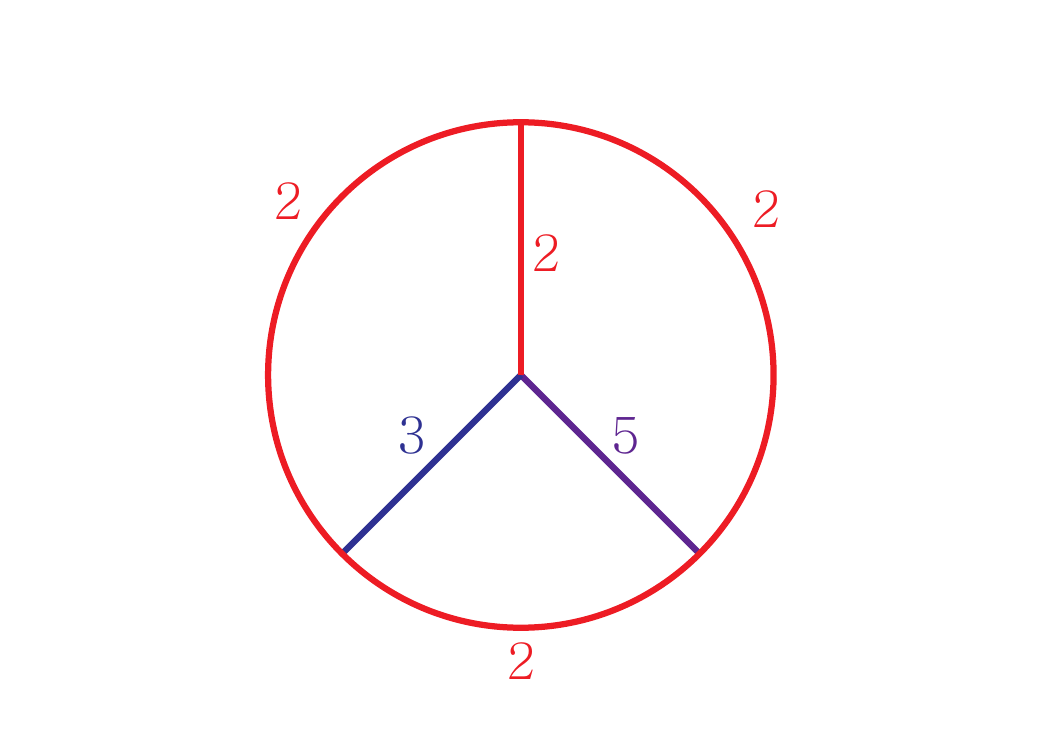} \\[-1mm]
{\small\textup{(a)}} & {\small\textup{(b)}} \\[2mm]
\multicolumn{2}{c}{\HandleDiagram{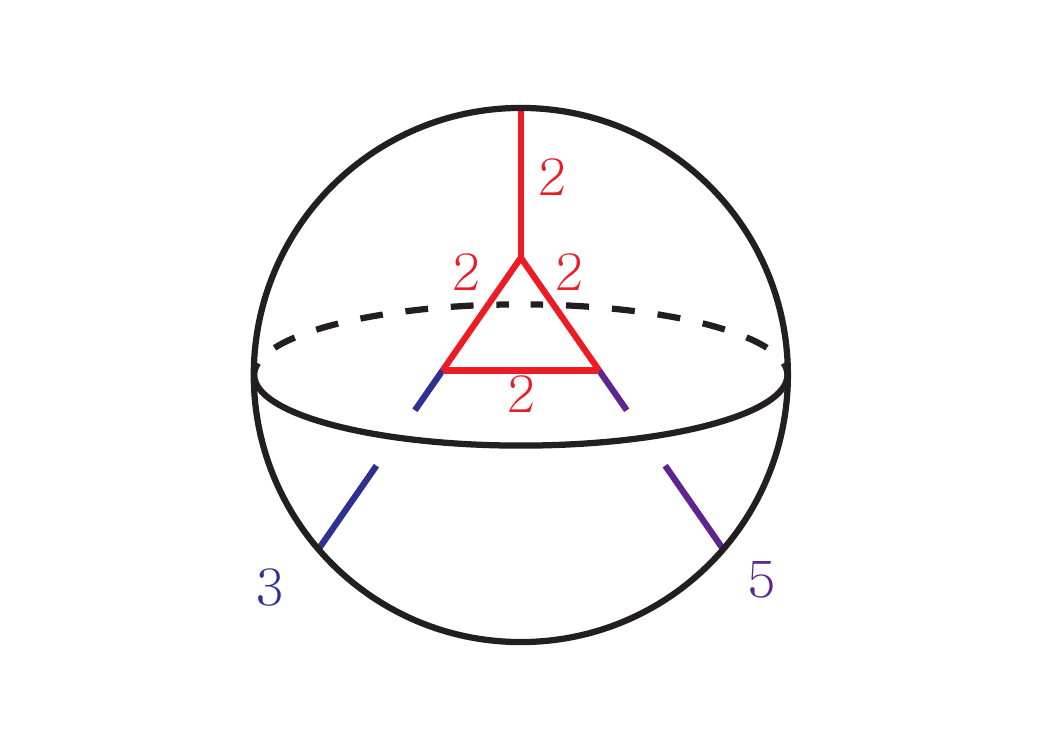}} \\[-1mm]
\multicolumn{2}{c}{{\small\textup{(c)}}}
\end{tabular}
\caption{\textup{(a)} \(N_{\mathbf I}\). \textup{(b)} The labeled singular graph of \(S^3\sslash G_{31}\). \textup{(c)} \(\mathcal M(\mathbf I,\alpha)\).}\label{fig:G31-1h}
\end{figure}

\begin{figure}[!htbp]
\centering
\begin{tabular}{@{}c@{\hspace{6mm}}c@{}}
\HandleDiagram{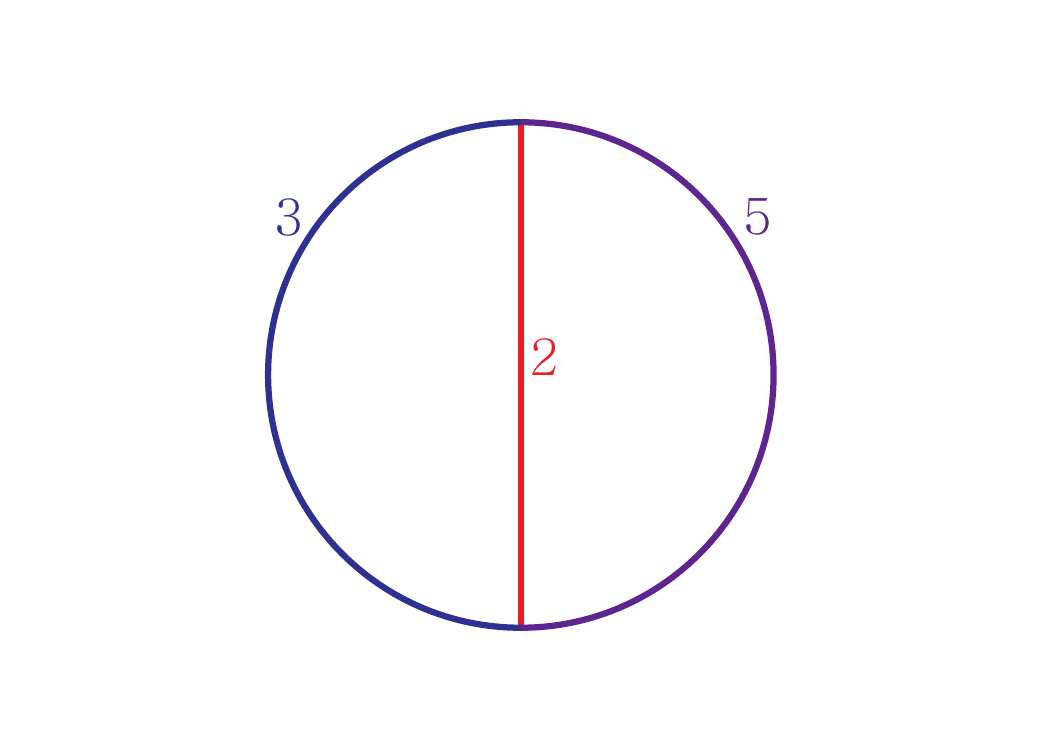} &
\HandleDiagram{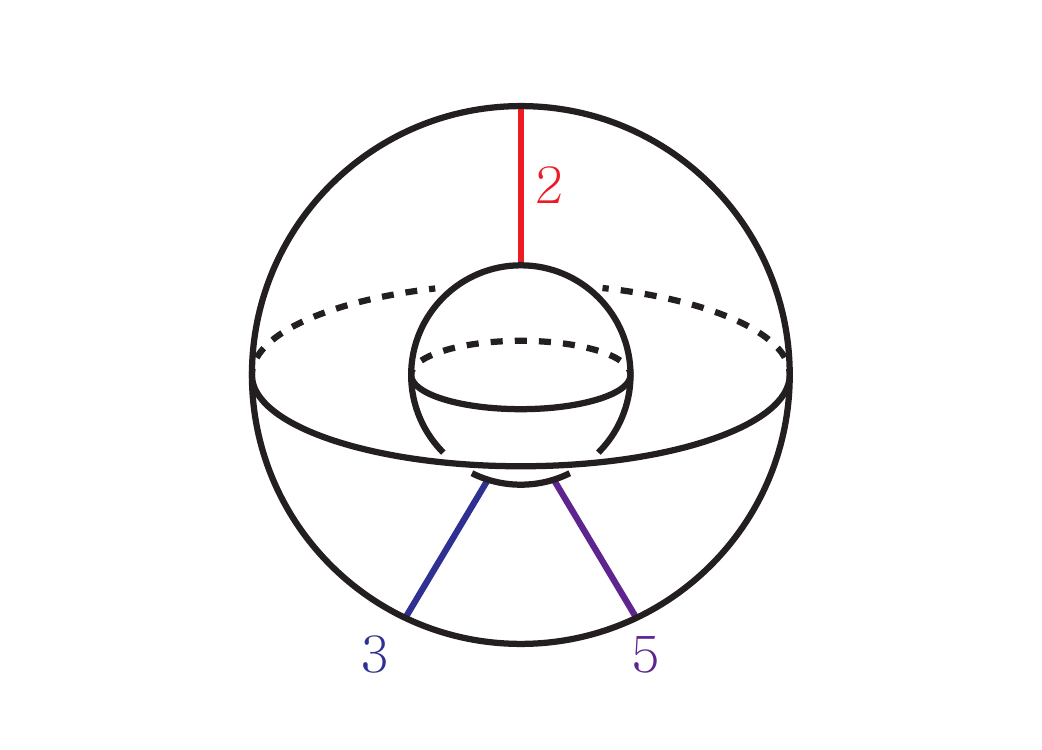} \\[-1mm]
{\small\textup{(a)}} & {\small\textup{(b)}}
\end{tabular}
\caption{\textup{(a)} The labeled singular graph of \(S^3\sslash G_{31'}\). \textup{(b)} \(I(\Sigma_{\mathbf I})\).}\label{fig:G31p-1h}
\end{figure}

\FloatBarrier

\section*{Declaration of AI use}

The second author used ChatGPT (OpenAI) for English-language editing, LaTeX assistance, reference checking, organization of the exposition, and assistance in the review and revision of selected mathematical arguments. The authors take full responsibility for all mathematical statements and proofs in this paper.

\end{document}